\documentclass[a4paper, 12pt]{article}
\usepackage{enumerate, theorem}
\usepackage{amsmath, amsfonts, amssymb}
\usepackage[height=22.5cm, width=15cm]{geometry}
\usepackage[all]{xy}
\usepackage{mathrsfs, yfonts}

\DeclareMathOperator{\id}{id}

\DeclareMathOperator{\C}{\mathbb{C}}

\newcommand{\A}{\tilde{\mathcal{A}}}

\newcommand{\parag}[1]{\paragraph{\sc{#1.}}}

\newtheorem{thm}{Theorem}[subsection]
\newtheorem{defn}[thm]{Definition}
\newtheorem{cor}[thm]{Corollary}
\newtheorem{prop}[thm]{Proposition}
\newtheorem{lemma}[thm]{Lemma}

\begin{document}

\title{ Contruction of quasi-invariant holomorphic \\ parameters for the Gauss-Manin connection of a holomorphic map to a curve.\\
(second version).}

\author{Daniel Barlet\footnote{Barlet Daniel, Institut Elie Cartan UMR 7502  \newline
Universit\'e de Lorraine, CNRS, INRIA  et  Institut Universitaire de France, \newline
BP 239 - F - 54506 Vandoeuvre-l\`es-Nancy Cedex.France. \newline
e-mail : Daniel.Barlet@iecn.u-nancy.fr}.}

\date{26/03/12}

\maketitle

\hfill {\it To my children who are my (quasi-)invariant parameters.}

\section*{Abstract.}

In this paper we consider holomorphic families of  frescos (i.e.  filtered differential equations with a regular singularity) and we construct a locally versal holomorphic family for every fixed Bernstein polynomial. We construct also several holomorphic parameters (a holomorphic parameter is a function defined on a set of isomorphism classes of frescos) which are quasi-invariant by changes of variable. This is motivated by the fact that a fresco is associated to a relative de Rham cohomology class on a one parameter  degeneration of compact complex manifolds, up to a change of variable in the  parameter. Then the value of a quasi-invariant holomorphic parameter on such data produces a holomorphic (quasi-)invariant of such a situation.

\parag{AMS Classification} 32 S 25, 32 S 40, 32 S 50.

\parag{Key words} Degenerating family of complex manifolds, Relative de Rham cohomology classes, filtered Gauss-Manin connection, asymptotic expansion, vanishing period, Theme, Fresco, geometric (a,b)-module.

\newpage

\tableofcontents

\section{Introduction.}

When we consider a proper holomorphic map \ $\tilde{f }: X \to C$ \ of a complex manifold \ $X$ \ on a smooth complex curve \ $C$ \  with a critical value at a point \ $0$ \ in \ $C$, the choice of a local coordinate near this point allows to dispose of an holomorphic function \ $f$. Then we may construct, using this function,  an (a,b)-module structure on the cohomology sheaves of the formal completion (in \ $f$) \ of the complex of sheaves \ $(Ker\, df^{\bullet},d^{\bullet})$ \ (see [B.12]). These (a,b)-modules represent a filtered version of the Gauss-Manin connection of \ $f$. The most simple example of this construction is the Brieskorn module (see [Br.70]) of a function with an isolated singular point. See also  [B.08] for the case of a germ of holomorphic function with  a 1-dimensional critical locus. \\
But it is clear that this construction depends seriously on the choice of the function \ $f$, that is to say on the choice of the local coordinate near the critical point \ $0$ \ in the complex curve \ $C$.\\
The aim of the present paper is to study the behaviour of such constructions when we make a change of local coordinate near the origin. We consider the case of frescos  which are monogenic geometric (a,b)-modules corresponding  to a minimal filtered differential equation associated to a relative de Rham cohomology class on \ $X$ \ (see [B.09] and [B.10]).  \\
 A holomorphic parameter is a function on the set of isomorphism classes of frescos which behave holomorphically in a holomorphic family of frescos. In general, a holomorphic parameter is not invariant by a change of variable, but we prove a theorem of stability of holomorphic families of frescos by a  change of variable (see  theorem \ref{chgt. var. hol.}) and this implies that a holomorphic parameter gives again a holomorphic parameter by a change of variable.\\
 Our main task is the construction of holomorphic parameters  which are (quasi)-invariant by all changes of variable. The first observation is the fact that the {\em principal Jordan-H{\"o}lder} sequence of a fresco is  covariant with the changes variable (see theorem  \ref{J-H. hol.}). That implies that any holomorphic parameter which is  (quasi-)invariant by changes of variable obtained for some quotient between  two terms of this sequence will produce a holomorphic parameter which is (quasi-)invariant by changes of variable for the fresco itself. Now the fact that the parameter of a rank \ $2$ \ fresco  is a holomorphic parameter which is (quasi-)invariant by changes of variable will produce a first collection of holomorphic parameters which are (quasi-)invariant by changes of variable. These holomorphic parameters are associated to Jordan blocks of size at least \ $2$ \ of the monodromy.\\
 Another construction of holomorphic parameters, which are (quasi-)invariant by changes of variable, is given. It  is non trivial even in the semi-simple case for which the previous construction does not give anything. In this case the holomorphic parameters look like some  "cross ratio" of eigenvectors of the monodromy.\\
In order to enlarge the use of this second construction, we prove  that for any holomorphic family of \ $[\lambda]-$primitive frescos parametrized by a reduced complex space \ $X$ \ there exists a dense Zariski open set \ $X'$ \ in \ $X$ \ on which the family of semi-simple parts of the frescos is again a holomorphic family. \\
 It is clear that combining these two constructions and this result we are able to produce many interesting  holomorphic parameters which are (quasi-)invariant by changes of variable for a general holomorphic family of  frescos.\\
 
 To finish this introduction let me decribe now the most significant results in this paper. After some preliminaries on the change of variable of a fresco, we introduce the notion of holomorphic family of frescos and we prove some basic results on this notion : theorem \ref{hol. dual.} which is the stability of holomorphic families by twisted duality and theorem \ref{deux sur trois} which says that in a short exact sequence of \ $\hat{A}_X-$sheaves, when two are holomorphic families the third is also a holomorphic family. The theorem \ref{J-H. hol.} shows that the quotient of two terms in the principal Jordan-H{\"o}lder sequence in a holomorphic families define also holomorphic families. Then, in section 4, we construct  a "standard" locally versal holomorphic family for any given fundamental invariants.\\
 In section 5 we prove that the collection of semi-simple parts of members  in a holomorphic family of \ $[\lambda]-$primitive frescos  is a holomorphic family on a Zariski dense open set in the parameter space,. This result is rather interesting because it is a tool to produce quasi-invariant holomorphic parameters (see above) which do not come from the principal J-H. sequence. \\
 
In order to construct more quasi-invariant holomorphic parameters than these we obtain from the rank \ $2$ \ quotients of the principal J-H. sequence, we study the rank \ $3$ \ case in detail. The classification of rank \ $3$ \  frescos is given in the appendix and the study of the changes of variable in rank \ $3$ \ leads to the construction of a new holomorphic parameter which behave in a very simple way by change of variable. This allows to produce a new quasi-invariant holomorphic parameter for each pair of rank \ $3$ \ quotients of terms in the principal J-H. sequence. These results are collected in theorems \ref{second.parm.} and \ref{inv. param.}.\\

\section{Change of variable in the ring \ $\hat{A}$.}

\subsection{Some facts on frescos.}

For basic properties of (a,b)-modules the reader may consult [B.93] [B.09] or [B.12].

\begin{defn}\label{fresco}
We call a geometric monogenic (a,b)-module a {\bf fresco}.
\end{defn}

\begin{defn}\label{primitive}
A geometric (a,b)-module is call {\bf \ $[\Lambda]-$primitive} for some subset\ $[\Lambda] \subset \mathbb{Q}\big/\mathbb{Z}$ \ if each root of its Bernstein polynomial is in \ $[-\Lambda]$.
\end{defn}

For any geometric (a,b)-module and any subset \ $[\Lambda] \in \mathbb{Q}\big/\mathbb{Z}$ \ there exists an unique maximal submodule \ $E[\Lambda]$ \ which is \ $[\Lambda]-$primitive. It is a normal submodule. The quotient \ $E\big/E[\Lambda]$ \ is the maximal \ $[^c\Lambda]-$primitive quotient of \ $E$, where \ $ [^c\Lambda]$ \ is the complement of \ $[\Lambda]$ \ in \ $\mathbb{Q}\big/\mathbb{Z}$.\\

Recall that is \ $F \subset E$ \ is a normal\footnote{A sub-module \ $F \subset E$ \ is normal if \ $E\big/F$ \ is again free on \ $\C[[b]]$. Note that for \ $G \subset F \subset E$ \ it is equivalent to ask that \ $G$ \ is normal in \ $F$ \ and \ $F$ \ normal in \ $E$ \ or to ask that \ $G$ \ and \ $F$ \ are normal is \ $E$.}  submodule of a fresco \ $E$, then \ $F$ \ and \ $E\big/F$ \ are frescos. Moreover, if \ $E$ \ is \ $[\Lambda]-$primitive, so are \ $F$ \ and  \ $E\big/F$.\\

As any regular (a,b)-module, a \ $[\Lambda]-$primitive fresco \ $E$ \ admits a Jordan-H{\"o}lder sequence, which is a sequence of normal sub-modules (which are frescos)
 $$0 = F_0 \subset F_1 \subset \cdots \subset F_{k-1} \subset F_k = E $$
 such that the quotients \ $F_{j}\big/F_{j-1}$ \ have rank 1.
  Then, for each \ $j \in [1,k]$ \ we have \ $F_{j}\big/F_{j-1} \simeq E_{\lambda_j} \simeq \hat{A}\big/\hat{A}.(a - \lambda_j.b)$.\\
  For any J-H. sequence of a \ $[\Lambda]-$primitive  fresco \ $E$ \ the sequence of  numbers  \ $\mu_1, \dots, \mu_k$ \ associated to the consecutive quotients are such that the (unordered) set  \ $\{\mu_1+1, \dots, \mu_k+k\}$ \ is independant of the choice of the J-H. sequence. \\
  Recall that for a  fresco \ $E$ \ the Bernstein polynomial \ $B_E$ \ is the {\em characteristic} polynomial of the action of \ $-b^{-1}.a$ \ on \ $E^{\sharp}\big/b.E^{\sharp}$. So its degree is equal to the rank of \ $E$ \ over \ $\C[[b]]$. The roots of \ $B_E$ \ are the opposite of the numbers  \ $\{\mu_1+1, \dots, \mu_k+k\}$.\\
  It is proved in  [B.09] proposition 3.5.2  that any \ $[\lambda]-$primitive fresco admits  a J-H. sequence such that 
  $$ \lambda_1+1 \leq \lambda_2+2 \leq \cdots \leq \lambda_k + k.$$
   Moreover  it is proved   in   [B.12]  proposition 1.2.5 such a J-H. sequence is unique (so the sub-modules \ $F_j$ \ are uniquely determined).\\
    It is called the {\bf principal J-H. sequence of \ $E$}. \\
  In {\em loc. cit.} we call {\bf fundamental invariants} of the rank \ $k$ \ $[\lambda]-$primitive fresco \ $E$ \ the ordered sequence of  numbers \ $(\lambda_1, \dots, \lambda_k)$ \ corresponding to the principal J-H. sequence. It is often convenient, in the \ $[\lambda]-$primitive case, to define the numbers \ $p_j$ \ by the formula \ $\lambda_{j+1} : = \lambda_j + p_j -1$ \ for \ $j \in [1,k-1]$ \ and to give the fundamental invariants via the numbers \ $\lambda_1, p_1, \dots, p_{k-1}$. Note that \ $p_j$ \ are natural integers and that \ $\lambda_1$ \ is rational number such that \ $\lambda_1 > k-1$. 
  \smallskip
 
Remark that the uniqueness of \ $F_1$ \ for the principal J-H. sequence of a \ $[\lambda]-$primitive fresco \ $E$ \ implies that the dimension of \ $Ker\, (a - \lambda_1.b)$ \ is a one dimensional vector space and that \ $F_1  = \C[[b]].e_1$ \ where \ $e_1$ \ is any non zero vector in this kernel.

\parag{Notation and convention} We want to extend the previous terminology to general frescos. So we fix now a total order on \ $\mathbb{Q}\big/\mathbb{Z}$. Then for each regular (a,b)-module  there exists an unique filtration of \ $E$ \ by normal submodules \ $(G_h)_{h \in H}$, that we shall call  the {\bf primitive filtration} of \ $E$ \ (associated to the order choosen above), such that
\begin{itemize}
\item If \ $E$ \ is \ $[\Lambda]-$primitive with \ $[\Lambda] = \{ [\mu_1]<  \dots < [\mu_H]\}$ \ then  we have
 \ $G_h : = E[\{[\mu_1], \dots,[\mu_h]\}]$.\\
 \end{itemize}
 
For a \ $[\Lambda]-$primitive fresco \ $E$ \  we define the {\bf principal Jordan-H{\"o}lder sequence} as the unique J-H.  sequence  of \ $E$ \  satisfying the following properties
\begin{enumerate}
\item This J-H. sequence is compatible with the primitive filtration \ $(G_h)_{h \in [1, H]}$ \ of \ $E$.
\item In each quotient \ $G_{h}\big/G_{h-1}$ \ of the primitive filtration of \ $E$ \ which is a \ $[\mu_h]-$primitive fresco, the J-H. sequence induces the principal  J-H. sequence.
\end{enumerate}
Note that existence and  uniqueness of the principal J-H. sequence of a \ $[\lambda]-$primitive fresco implies existence and  uniqueness (for a given order on \ $\mathbb{Q}\big/\mathbb{Z}$) in the general case. \\
When we consider a \ $[\Lambda]-$primitive fresco with \ $[\Lambda] =  \{ [\mu_1]<  \dots < [\mu_H]\}$ \ we shall call \ $(\lambda_1, \dots, \lambda_k)$ \ the {\bf fundamental invariants} of \ $E$ \ when the numbers \ $\lambda_1, \dots, \lambda_k$ \ are ordered in the same way that the successive rank \ $1$ \ quotients of the principal J-H. sequence of \ $E$. For instance \ $\lambda_1$ \ is the first fundamental invariant of \ $E[\mu_1]$ \ where \ $[\mu_1]$ \ is the smallest class in  \ $\mathbb{Q}\big/\mathbb{Z}$ \ of the opposite of the roots of the Bernstein polynomial of  $E$.

\subsection{Embedding frescos.}

\begin{defn}\label{asympt. 0}
Let \ $\lambda$ \ be a rational number in \ $]0,1]$ \ and \ $N$ \ an integer. We define the left \ $\hat{A}-$module \ $\Xi_{\lambda}^{(N)}$ \ as the free \ $\C[[b]]$ \ module generated by \ $e_0, \dots,e_N$ \ with the action of \ $a$ \ define by the following rules :
\begin{enumerate}
\item \ $ a.e_0 = \lambda.b.e_0 $ ;
\item For \ $j \in [1,N]$ \  \ $a.e_j = \lambda.b.e_j + b.e_{j-1} $ ;
\item The left action of \ $a$ \ is continuous for the \ $b-$adic topology of \ $\Xi_{\lambda}^{(N)}$. As \ $a.\Xi_{\lambda}^{(N)} \subset b.\Xi_{\lambda}^{(N)}$ \ the action of \ $a$ \ extends to an action of \ $\C[[a]]$\\
\end{enumerate}
\end{defn}

It is an easy computation to check that \ $a.b - b.a = b^2$ \ on \ $\Xi_{\lambda}^{(N)}$, so we have defined a simple pole (a,b)-module of rank \ $N+1$.\\
Remark that \ $\Xi_{\lambda}^{(N)}$ \ corresponds to the formal asymptotic expansions of the type
$$ \sum_{j = 0}^N \  S_j(b).s^{\lambda-1}.\frac{(Log\, s)^j}{j!} = \sum_{j = 0}^N \  T_j(s).s^{\lambda-1}.\frac{(Log\, s)^j}{j!} $$
where the action of \ $a$ \ is the multiplication with \ $s$ \ and \ $b$ \ is the "primitive without constant".
For \ $[\Lambda] \subset \mathbb{Q}\big/\mathbb{Z}$ \ a finite subset, we define \ $\Xi_{\Lambda}^{(N)}$ \ 
as the direct sum
$$ \Xi_{\Lambda}^{(N)} : = \oplus_{\lambda_i \in [\Lambda]} \quad \Xi_{\lambda_i}^{(N)} .$$
When \ $V$ \ is a finite dimensional  complex vector  space, we shall also consider vector-valued asymptotic expansions. Then we let  \ $a$ \ and \ $b$ \ acts on \ $\Xi_{\Lambda}^{(N)} \otimes V$ \ as \ $a \otimes \id_V$ \ and \ $b \otimes \id_V$.\\

Now the following embedding theorem is an elementary by-product of the theorem 4.2.1 in [B.09] (The  more precise version for \ $[\lambda]-$primitive frescos stated  in the remark below  is given in [B.12].)

\begin{thm}\label{embed.}
Let \ $E$ \ be any rank \ $k$ \ geometric (a,b)-module. Then there exists an integer \ $N \leq k-1$, a complex vector space \ $V$ \ of dimension \ $\leq k$ \ and an embedding of \ $E$ \ in the left \ $\hat{A}-$module \ $\Xi_{\Lambda}^{(N)} \otimes V$, where \ $[\Lambda]$ \ is the subset of \ $\mathbb{Q}\big/\mathbb{Z}$ \ of the classes containing an opposite of a root of the Bernstein polynomial of \ $E$.
\end{thm}

\parag{Remark} When \ $E$ \ is \ $[\lambda]-$primitive we may choose \ $[\Lambda] = [\lambda]$ \ and the minimal \ $N$ \ is equal to \ $d(E)-1$ \ where \ $d(E)$ \ is the ss-depth of \ $E$;   the minimal dimension for \ $V$ \ is \ $\delta(E)$ \ the rank of the semi-simple part of \ $E$ \ (see below for the corresponding definitions.) $\hfill \square$\\

Let me recall basic definitions and results  proved in [B.12]  for  frescos.

\begin{defn}\label{theme + ss.}
Let \ $E$ \ be a  fresco. We say that \ $E$ \ is a {\bf theme} if there exists an embedding of \ $E$ \ in \ $\Xi_{\Lambda}^{(N)}$ \ for some integer \ $N$.\\
We say that \ $E$ \ is {\bf semi-simple} when any \ $\hat{A}-$linear map \ $\varphi : E \to \Xi_{\lambda}^{(N)}$ \ has rank at most equal to\ $1$ \ for any choice of \ $[\lambda]$.
\end{defn}

\parag{Remarks} 
\begin{enumerate}
\item Of course any rank 1 fresco is at the same time a theme and semi-simple. But it is clear that in rank \ $\geq 2$ \ a \ $[\lambda]-$primitive fresco cannot be at the same time a theme and semi-simple.
\item Themes and semi-simple frescos have nice stability properties. See the [B.09] and [B.12].
\end{enumerate}

\begin{prop}[{See  [B.12]} ]\label{structure}
Let \ $E$ \ be a  fresco. There exists a biggest semi-simple submodule \ $S_1(E)$ \ in \ $E$. It is normal and has rank \ $\delta(E)$ \ at least equal to \ $1$. Then we define inductively the semi-simple filtration \ $(S_j(E))_{j\in [1,d]}$ \ by the formula
$$ S_{j+1}(E) : =  q_j^{-1}\left[ S_1(E\big/S_j(E))\right] $$
where \ $q_j : E \to E\big/S_j(E)$ \ is the quotient map. 
\end{prop}

Then \ $d(E) : = \inf\{j \in \mathbb{N}, S_j(E) = E \} $ \ is called the {\bf ss-depth} of \ $E$. 
For instance for a  fresco \ $E$ \ we have \ $d(E) =1$ \ if and only if \ $E$ \ is semi-simple, and \ $d(E) = rk(E)$ \ if and only if \ $E$ \ is a \ $[\lambda]-$primitive theme.\\
 Note that  \ $\delta(E) = 1$ \ if and only if \ $E$ \ is a \ $[\lambda]-$primitive theme, and \ $\delta(E) = rk(E)$ \ if and only if \ $E$ \ is semi-simple.

\smallskip

The isomorphism class of \ $E$ \ determines the isomorphism classes of the normal sub-modules\ $S_j(E)$, and of the quotients \ $ E\big/S_j(E)$ .

\subsection{Definition and first properties of a change of variable.}

We shall work with the \ $\C-$algebra
\begin{equation}
\hat{A} : = \sum_{\nu = 0} ^{+\infty} \  P_{\nu}(a).b^{\nu} \quad {\rm where} \quad P_{\nu} \in \C[[a]] \tag{1}
\end{equation}
with the product law defined by the following two conditions
\begin{enumerate}
\item \ $a.b - b.a = b^2$ .
\item The right and left  multiplications by any \ $T \in \C[[a]]$ \ is continuous for the \ $b-$adic filtration of \ $\hat{A}$.
\end{enumerate}
The first condition implies  the identities 
 $$a.b^n = b^n.a + n.b^{n+1} \quad {\rm and} \quad  a^n.b = b.a^n + n.b.a^{n-1}.b$$ 
  and with the second condition we obtain
\begin{align*}
& a.S(b) = S(b).a + b^2.S'(b) \quad {\rm for \ any} \quad S \in \C[[b]] \quad {\rm and} \tag{2a}\\
& T(a).b = b.T(a) + b.T'(a).b \quad {\rm for \ any} \quad T \in \C[[a]] \tag{2b}
\end{align*}

\begin{lemma}\label{debut}
We have the following properties of the \ $\C-$algebra \ $\hat{A}$.
\begin{enumerate}
\item For any \ $x,y$ \ in \ $\hat{A}$ \ we have \ $x.y - y. x \in \hat{A}.b^2$.
\item For each \ $n \in \mathbb{N}$ \ we have  \ $b^n.\hat{A} = \hat{A}.b^n$ \ and it is a two-sided ideal in \ $\hat{A}$.
\item We have \ $a.\hat{A} + b.\hat{A} = \hat{A}.a + \hat{A}.b$ \ and any element \ $x \in \hat{A}$ \ is invertible if and only if \ $x \not\in a.\hat{A} + b.\hat{A}$.
\end{enumerate}
\end{lemma}

The easy proof is left to the reader. $\hfill \blacksquare$

%\parag{Proof} We let the reader prove \ $1$ \ as an exercice. Let me prove \ $2$ \ by induction on \ $n \geq 0$. Assume that this is true for \ $n$ \ and consider \ $b^{n+1}.x$, for some \ $x \in \hat{A}$. We may write \ $b^n.x = y.b^n$ \ and then \ $b^{n+1}.x = b.y.b^n$. Now, using the fact that \ $b.y = y.b + z.b^2$  \ from \ $1$ \ we obtain \ $b^{n+1}.x = y.b^{n+1} + z.b^{n+2} = (y + z.b).b^{n+1}$.\\
%Let \ $a.x + b. y $ \ for \ $x,y \in \hat{A}$, and write \ $x = x_0(a) + x'.b$ \ and \ $b.y = y'.b$ \ using \ $2$ \ for \ $n = 1$. Then \ $a.x + b.y = x_0(a).a + (a.x' + y').b$ \ and \ $a.\hat{A} + b.\hat{A} \subset  \hat{A}.a + \hat{A}.b$. Consider now  \ $u.a + v.b$ \ with \ $u,v \in \hat{A}$. Put \ $u = u_0(a) + b.u' $ \ using \ $(2b)$ \ and \ $v.b = b.v'$ \ using again \ $2$ \ for \ $n =1$. Then we have \ $u.a + v.b = a.u_0(a) + b.(u'.a + v')$.\\
%Now consider  \ $x =  \sum_{\nu = 0} ^{+\infty} \  P_{\nu}(a).b^{\nu} $ \ with \ $P_0(0) \not= 0$. Then \ $P_0(a)$ \ is an invertible element of \ $\C[[a]]$. Now \ $P_0(a)^{-1}.x = 1 - b.y$ \ for some \ $y \in \hat{A}$. But now the serie  \ $ \sum_{n=0}^{+\infty} \ (b.y)^n $ \ converges for the \ $b-$adic filtration of \ $\hat{A}$ \ as \ $(b.y)^n$ \ lies in \ $b^n.\hat{A}$ \ thanks to \ $1$ \ and \ $2$. Then \ $1 -b.y$ \ is invertible in \ $\hat{A}$, and so is\ $x$.\\
%Conversely, if \ $x$ \ is invertible, then so is the image of \ $x$ \ in \ $\C[[a]] = \hat{A}\big/b.\hat{A}$, and  we have \ $x_0(0) \not= 0$ \ and \ $x = x_0(0) + a.y + b.z $.
% Then \ $x$ \ does not lies in \ $a.\hat{A} + b.\hat{A}$. $\hfill \blacksquare$

\bigskip

When \ $E$ \ is a regular (a,b)-module, it is in a natural way a \ $\hat{A}-$module, because it is complete for the \ $a-$adic filtration\footnote{in fact the existence of an integer \ $N$ \ such that \ $a^N.E \subset b.E$ \ is enough for that purpose.}. Then we may identify the category of  regular (a,b)-modules with the subcategory of left \ $\hat{A}-$module that they defined. The following easy lemma make this fact more explicit. Recall that \ $\A$ \ is the subalgebra of  \ $\hat{A}$ \ in which the coefficients of \ $b^{\nu}$ \ are polynomials in \ $a$. So an (a,b)-module is a left \ $\A-$module which is free and finite type on \ $\C[[b]] \subset \A \subset \hat{A}$.

\begin{lemma}\label{div.}
Let \ $E$ \ a regular (a,b)-module and let  \ $u \in \hat{A}$ \ and \ $x \in E$. Then there exists \ $v \in \A$ \ such that  \ $u.x = v.x$.
\end{lemma}

\parag{proof} As \ $\A.x \subset E$ \ is a regular (a,b)-module there exists (see [B.09] ) a monic polynomial in \ $a$ \ with coefficients in \ $\C[[b]]$ \ and initial form in (a,b) homogeneous of degree \ $k$, such that \ $P.x = 0$. Now write for each \ $n \geq k$
$$ a^n = Q_n.P + R_n $$
where \ $Q_n $ \ is in \ $\hat{A}$ \ and \ $R_n$ \ is a polynomial in \ $a$ \ of degree \ $\leq k-1$. Assume that we know that \ $val_{a,b}(Q_n) \geq n-k$ \ and \ $val_{a,b}(R_n) \geq n$. We shall prove that this implies \ $val_{a,b}(Q_{n+1}) \geq n+1-k$ \ and \ $val_{a,b}(R_{n+1}) \geq n+1$. Write
$$ R_n : = \sum_{j=1}^k \ a^{k-j}.R_n^j .$$
Then we have from our assumption \ $val_b(R_n^j) \geq n-k+j$. Now write
$$ a^k.R^1_n = R_n^1.P + S_n $$
where \ $S_n$ \ is a polynomial of degree \ $\leq k-1$ \ in \ $a$ \ and with \ $val_{a,b}(S_n) \geq n+1$ \ because \ $val_b(R^1_n) \geq n-k+1$. Now the equality
$$ a^{n+1} = (a.Q_n + R_n^1).P + S_n + \sum_{j=2}^k a^{k-j+1}.R_n^j $$
gives \ $Q_{n+1} = a.Q_n + R_n^1$ \ and \ $R_{n+1} = S_n + \sum_{j=2}^k a^{k-j+1}.R_n^j$. So we obtain \ $val_{a,b}(Q_{n+1}) \geq n+1-k$ \ as  \ $a.Q_n$ \ and \ $R^1_n$ \ have valuation in \ $(a,b)$ \ at least equal \ $n+1-k$, and also \ $val_{a,b}(R_{n+1}) \geq n+1-k$ \ because \ $S_n$ \ has valuation in (a,b) at least equal to \ $n+1$ \ and each term \ $a^{k-j+1}.R_n^j $ \ also for \ $j \in [2,k]$.\\
For \ $n = k$ \ we have \ $Q_k = 1$ \ and \ $R_k = P - a^k$ \ so our hypothesis is true and we may begin the induction. Now write
$$ u = \sum_{n = 0}^{\infty} \ u_n.a^n $$
where \ $u_n$ \ is in \ $\C[[b]]$. So we obtain
$$ u = \sum_{n=0}^{k-1} \ u_n.a^n  + \big(\sum_{n=k}^{\infty} u_n.Q_n \big).P + \sum_{n=k}^{\infty} u_n.R_n $$
where the convergence is obtained from the previous inequalities on valuations. Then we may take \ 
$v : =  \sum_{n=0}^{k-1} \ u_n.a^n + \sum_{n=k}^{\infty} u_n.R_n $ \ which is a polynomial in \ $a$ \ of degree \ $\leq k-1$ \ with coefficients in \ $\C[[b]]$. $\hfill \blacksquare$

\parag{Important remark}  In the sequel we shall always identify frescos with the corresponding left \ $\hat{A}-$modules that they define.\\
 Note that two frescos are isomorphic as (a,b)-modules ( that is to say as left \ $\A-$modules) if and only if the corresponding left \ $\hat{A}-$modules are isomorphic. $\hfill \square$

\begin{lemma}\label{chgt. variable 0}
Let \ $\theta \in \C[[a]]$ \ such that \ $\theta(0) = 0$ \ and \ $\chi(\theta) : = \theta'(0) \not= 0$. Then the elements \ $\alpha : = \theta(a)$ \ and \ $\beta : = b.\theta'(a)$ \ satisfy the commutation relation in \ $\hat{A}$ :
$$ \alpha.\beta - \beta.\alpha = \beta^2 .$$
\end{lemma}

\parag{Proof} This is a simple computation using the identity in \ $(2b)$
\begin{align*}
& \alpha.\beta - \beta.\alpha  = \theta(a).b.\theta'(a) - b.\theta'(a).\theta(a) = b.\theta'(a).b.\theta'(a) = \beta^2
\end{align*}$ \hfill \blacksquare$\\

So for any such \ $\theta$ \ there exists an unique automorphism of \ $\C-$algebra \ $\Theta : \hat{A} \to  \hat{A} $ \ such that \ $ \Theta(1) = 1, \Theta(a) = \alpha , \Theta(b) = \beta $ \ which is continuous for the \ $a-$adic and  \ $b-$adic filtrations. We have, when \ $x = \sum_{\nu = 0} ^{+\infty} \ P_{\nu}(a).b^{\nu} $
$$ \Theta(x) = \sum_{\nu = 0} ^{+\infty} \ P_{\nu}(\alpha).\beta^{\nu}. $$

\begin{defn}\label{chgt. variable 1}
We shall say that  \ $\theta \in \C[[a]]$ \ such that \  $\theta(0) = 0$ \ and \ $\chi(\theta) = \theta'(0) \not= 0$ \ is  a {\bf change of variable} in \ $\hat{A}$. For any left \ $\hat{A}-$module \ $E$, we define a new \ $\hat{A}-$module \ $\theta_*(E)$, called {\bf the change of variable of \ $E$ \ for \ $\theta$}, by letting \ $\hat{A}$ \ act on \ $E$ \ via the automorphism \ $\Theta$.\\
We shall say that \ $\theta$ \ is {\bf unimodular} when \ $\chi(\theta) : = \theta'(0) = 1$.
\end{defn}

Note that, as a \ $\C-$vector space, we have \ $\theta_*(E) = E$ \ and for \ $u \in \hat{A}$ \ and \ $e \in \theta_*(E)$ \ we put \ $u.e : = \Theta(u)._Ee$, where on the right handside \ $\hat{A}$ \ acts on \ $e$ \ as an element of \ $E$.

\parag{Notation} When we consider an element \ $x \in \theta_*(E)$ \ we write \ $a.x$ \ and \ $b.x$ \ for the action of \ $a$ \ and \ $b$ \ on the \ $\hat{A}-$module \ $\theta_*(E)$. Considering now \ $x$ \ as an element of \ $E$, this means that we look in fact at \ $\alpha.x$ \ and \ $\beta.x$ \ respectively. So we have to distinguish carefully if we are looking at \ $x$ \ as an element of \ $\theta_*(E)$ \ or as an element in \ $E$.

\begin{lemma}\label{chgt. variable 2}
For any \ $\hat{A}-$module \ $E$ \ and any change of variable \ $\theta$ \ we have
$$ b^n.\theta_*(E) = \beta^n.E =  b^n.E \quad \forall \ n \in \mathbb{N}.$$
\end{lemma}

\parag{Proof} We shall prove the statement (i. e. the second equality)  by induction on \ $n \geq 0$. As the case \ $n = 0$ \ is clear, assume that the equality \ $\beta^n.E = b^n.E$ \ is proved and we shall prove it for \ $n+1$. The inclusion \ $\beta^{n+1}.E \subset b^{n+1}.E$ \ is easy because \ $\beta = b.\theta'(a)$, $b^n.E$ \ is stable by \ $\hat{A}$ \ as \ $b^n.\hat{A} = \hat{A}.b^n$ \ and our induction hypothesis allows to conclude.\\
Assume that \ $x \in E$. As \ $\theta'(0)$ \ is not \ zero, $\theta'(a)$ \ is invertible in \ $\hat{A}$, the element  \ $\theta'(a)^{-1}.b^n.x$ \ is in \ $b^n.E = \beta^n.E$ \ and so there exists \ $y \in E$ \ such that \ $\beta^n.y = \theta'(a)^{-1}.b^n.x$. Then \ $\beta^{n+1}.y = b.\theta'(a).\theta'(a)^{-1}.b^n.x = b^{n+1}.x$, which gives the desired  inclusion \ $b^{n+1}.E \subset \beta^{n+1}.E. \hfill \blacksquare$\\

\begin{cor}\label{chgt. variable 3}
If we have a \ $\hat{A}-$linear map \ $f : E \to F$ \ between two left \ $\hat{A}-$modules, then the {\bf same map} induces a \ $\hat{A}-$linear map between \ $\theta_*(E)$ \ and \ $\theta_*(F)$. So \ $\theta_*$ \ is an exact functor.
\end{cor}

\parag{Proof} The map \ $f$ \ is \ $\C[[a]]-$linear and also \ $\C[[b]]-$linear. So it is \ $\C[[\alpha]]-$linear. It is also \ $\C[\beta]-$linear because we have
$$ f(\beta.x) = f(b.\theta'(a).x) = b.\theta'(a).f(x)= \beta.f(x) .$$
It is now enough to see that \ $f$ \ is continuous for the \ $\beta-$adic filtration to conclude that \ $f$ \ is \ $\C[[\beta]]-$linear. But as the previous lemma shows that the \ $\beta-$adic filtration is equal to the \ $b-$adic filtration for which \ $f$ \ is continuous, the map \ $f$ \ is \ $\hat{A}-$linear from \ $\theta_*(E)$ \ to \ $\theta_*(F)$. $\hfill \blacksquare$

\parag{Remarks}
\begin{enumerate}
\item Thanks to the previous corollary, it is obvious that any change of\ variable  $\theta$ \ defines a functor \ $\theta_*$ \  on the category of left \ $\hat{A}-$module into itself, which is the identity on morphisms, and it is an automorphism of this category. \\
Warning : it is not true in general that \ $\theta_*(E)$ \ is isomorphic to \ $E$, even when \ $\theta$ \  is unimodular (see section 6.3).
\item Remark also that for any left \ $\hat{A}-$module \ $E$ \  the sub$-\hat{A}-$modules of the \ $\hat{A}-$module \ $\theta_*(E)$ \ are the same sub-vector spaces than the sub$-\hat{A}-$modules of \ $E$ \ of the  given  vector space \ $E = \theta_*(E)$.
\item If the left \ $\hat{A}-$module \ $E$ \ is finitely generated, then so is \ $\theta_*(E)$. This a easy consequence of the fact that \ $\theta_*(\hat{A})$ \ is isomorphic to \ $\hat{A}$ \ as a left \ $\hat{A}-$module via the map \ $\Theta$.
\end{enumerate}

\begin{lemma}\label{chgt. variable 4}
Let \ $E$ \ be a left \ $\hat{A}-$module. Then \ $T_b(E)$, the \ $b-$torsion sub-$\C[[b]]-$module of \ $E$,  is a left  \ $\hat{A}-$sub-module. For any change of variable \ $\theta$ \ we have the equality
 $$ \theta_*(T_b(E)) = T_b(\theta_*(E)) = T_{\beta}(E).$$
 If \ $E$ \ is finitely generated on \ $\C[[b]]$ \ then \ $\theta_*(E)$ \ is finitely generated on \ $\C[[b]]$
 \end{lemma}
 
\parag{Proof} First remark that for any left \ $\hat{A}-$module \ $E$ \ and for each \ $n \in \mathbb{N}$ \ the sub-vector space \ $Ker \, b^n$ \ is an \ $\hat{A}-$submodule of \ $E$. For \ $n = 0$ \ this is clear. For \ $n = 1$ \ it is an immediate  consequence of the formula \ $b.T(a) = T(a).b + b.T'(a).b$. So assume that \ $n \geq 2$ \ and the assertion true for \ $Ker \, b^{n-1}$. Let \ $x \in Ker \, b^n$ \ and \ $T \in \C[[a]]$. Then we have\ $b^n.T(a).x = b^{n - 1}.T(a).b.x + b^n.T'(a).b.x$; but \ $b.x$ \ is in \ $Ker \, b^{n-1}$ \ and also \ $T(a).b.x$ \ and \ $T'(a).b.x$ \ because our inductive assumption. We conclude that \ $T(a).x$ \ is in \ $Ker\, b^n$.
Now let us show that we have \ $Ker \, \beta^n = Ker \, b^n$ \ for each \ $n \in \mathbb{N}$. Again the case \ $n = 0$ \ is obvious and for \ $n=1$ \ we have \ $Ker \, \beta \subset Ker \, b$ \ because \ $\theta'(a)$ \ is invertible in \ $\hat{A}$ \ and \ $Ker \, b \subset Ker \, \beta $ \ because \ $Ker \, b$ \ is stable by \ $\C[[a]]$.\\
 Assume \ $n \geq 2$ \ and that  the equality \ $Ker \, \beta^{n-1} = Ker \, b^{n-1}$ \ is proved. If \ $b^n.x = 0$ \ then \ $b.x $ \ is in  \ $Ker \, b^{n-1} = Ker \, \beta^{n-1}$. Now \ $\beta^n.x = \beta^{n-1}.b.\theta'(a).x$ \ and
 $$ b.\theta'(a).x = \theta'(a).b.x - b.\theta''(a).b.x $$
 which is also in \ $Ker \, \beta^{n-1} = Ker \, b^{n-1}$ \ because this kernel is a \ $\hat{A}-$submodule.
 So \ $Ker \, b^n \subset Ker \, \beta^n$.\\
Conversely, if \ $\beta^n.x = 0$ \ we have \ $\beta.x \in Ker \, \beta^{n-1} = Ker \, b^{n-1}$ \ and then 
$$b^{n-1}.\beta.x = b^n.\theta'(a).x = 0 .$$
 So \ $\theta'(a).x$ \ is in \ $Ker \, b^n$ \ which is \ $\C[[a]]-$stable, and \ $\theta'(a)$ \ is invertible in \ $\C[[a]]$. So \ $x$ \ is in \ $Ker \, b^n$, and we have the equality \ $Ker \, b^n = Ker \, \beta^n$.\\
Then the \ $b-$torsion of \ $E$ \ co{\"i}ncides with the \ $b-$torsion of \ $\theta_*(E)$.\\
As we know that \ $b.E = \beta.E$, the vector spaces \ $E\big/b.E$ \ and \ $  \theta_*(E)\big/b.\theta_*(E)$ \ co{\"i}ncide. So the finite dimension of \ $E\big/b.E $ \ implies that \ $\theta_*(E)$ \ is finitely generated as a \ $\C[[b]]-$module. $\hfill \blacksquare$\\

We introduce in [B.08] (see also the appendix of [B.12]) the notion of a small \ $\A-$module on the algebra
$$ \A : = \{ \sum_{\nu = 0}^{+\infty} \ P_{\nu}(a).b^{\nu} \} \quad {\rm where} \quad P_{\nu} \in \C[a] ,$$
with conditions for \ $\A$ \ analoguous to conditions 1 and 2  for \ $\hat{A}$.\\
In fact we were mainly interested in the case where these modules are naturally \ $\hat{A}-$modules. The following definition is the \ $\hat{A}-$analog. \\
First define for any \ $\hat{A}-$module \ $E$ 
$$  \hat{T}_a(E) : = \{ x \in E \ / \  \C[[b]].x \subset  T_a(E) \} $$
where \ $T_a(E)$ \ is the \ $a-$torsion of \ $E$. Note that \ $T_a(E)$ \ is not stable by \ $b$ \ in general, so it is not a sub$-\A-$module (and "a fortiori" not a sub$-\hat{A}-$module).

\begin{defn}\label{small}
We shall say that a \ $\hat{A}-$module\ $E$ \  is {\bf small} when it satisfies the following conditions :
\begin{enumerate}
\item  \ $E$ \ is a finitely generated \ $\C[[b]]-$module.
\item We have the inclusion \ $T_b(E) \subset \hat{T}_a(E)$.
\item  There exists \ $N \in \mathbb{N}$ such \ $a^N. \hat{T}_a(E) = 0 $.
\end{enumerate}
\end{defn}

\begin{prop}\label{change small}
The class of small \ $\hat{A}-$modules is stable by any change of variable.
\end{prop}

The proof of this proposition will be an easy consequence of the lemma \ref{chgt. variable 4} and the next lemma.

\begin{lemma}\label{torsion a}
For any \ $\hat{A}-$module \ $E$, $\hat{T}_a(E)$ \ is the maximal sub$-\hat{A}-$module contained in \ $T_a(E)$.
\end{lemma}

\parag{Proof} It is clear that any sub$-\hat{A}-$module of \ $E$ \ contained in \ $T_a(E)$ \ is contained in \ $\hat{T}_a(E)$. So the only point to prove is that \ $\hat{T}_a(E)$ \ is a sub$-\hat{A}-$module. Let \ $x$ \ be an element in \ $\hat{T}_a(E)$ \ and \ $z$ \ be an element of \ $\hat{A}$. We shall prove that \ $z.x$ \ is in \ $T_a(E)$ \ which is enough to conclude.\\
As \ $x$ \ is in \ $T_a(E)$, there exists \ $N_0 \in \mathbb{N}$ \ with \ $a^{N_0}.x = 0$. Now write
$$ z = \sum_{\nu = 0}^{+\infty} \ b^{\nu}.P_{\nu}(a) $$
where \ $P_{\nu}$ \ is in \ $\C[[a]]$. Put for each \ $\nu \geq 0 $
$$ P_{\nu}(a) = Q_{\nu}(a).a^{N_0} + R_{\nu}(a) $$
with \ $Q_{\nu} \in \C[[a]]$ \ and \ $R_{\nu} \in \C[a]$ \ a polynomial of degree at most \ $N_0 - 1$. Then we have 
 $$z.x = \sum_{\nu = 0}^{+\infty} \ b^{\nu}.R_{\nu}(a).x = \sum_{j=0}^{N_0-1} \ S_j(b).a^j.x  .$$
 But we may write in \ $\hat{A}$
 $$  \sum_{j=0}^{N_0-1} \ S_j(b).a^j =  \sum_{j=0}^{N_0-1} \ a^j.T_j(b) $$
 where each \ $T_j$ \ is in \ $\C[[b]]$. As \ $T_j(b).x$ \ is in \ $T_a(E)$, there exists an integer \ $N_1$ \ such that \ $a^{N_1}.T_j(b).x = 0 $ \ for each \ $j \in [0,N_0-1]$. Then we get
 $$a^{N_1}.z.x =  \sum_{j=0}^{N_0-1} \ a^{N_1+j}.T_j(b).x = 0 $$
 and \ $z.x$ \ is in \ $T_a(E)$. $\hfill \blacksquare$\\
 
 As \ $T_b(E)$ \ is a sub$-\hat{A}-$module of \ $E$, we have equivalence between the two inclusions
 \begin{enumerate}
 \item \  $ T_b(E) \subset T_a(E)$ 
 \item \ $ T_b(E) \subset \hat{T}_a(E) $.
 \end{enumerate}
 For any change of variable \ $\theta$ \ the equality \ $T_a(E) = T_{\alpha}(E) = T_a(\theta_*(E))$ \ is obvious. So the previous lemma implies the equality \ $ \hat{T}_a(E) = \hat{T}_a(\theta_*(E))$ \ as \ $E$ \ and \ $\theta_*(E)$ \ have the same sub$-\hat{A}-$modules, and so the conditions 1,  2 and 3 of the definition \ref{small} are stable by any  change of variable, thanks to the lemma \ref{chgt. variable 4}. This proves the proposition \ref{change small}. $\hfill \blacksquare$
 
 \bigskip

 \begin{defn}\label{regular}
 We shall say that a left \ $\hat{A}-$module \ $E$ \ is a {\bf simple pole $(\hat{a},b)-$module} when \ $E$ \ is free and finitely generated over \ $\C[[b]]$ \ and satisfies \ $a.E \subset b.E$.\\
 We shall say that a left \ $\hat{A}-$module \ $E$ \ is a {\bf regular $(\hat{a},b)-$module} when it is contained in a simple pole $(\hat{a},b)-$module. 
 \end{defn}
 
 \parag{Remarks}
  \begin{enumerate}
 \item A  regular $(\hat{a},b)-$module is free and finitely generated on \ $\C[[b]]$ \ by definition. 
 \item The previous terminology is compatible with the terminology on regular  (a,b)-modules (see [B.93]) because an (a,b)-module (which is not "a priori" a \ $\C[[a]]-$\\
 module) is complete for the \ $a-$adic filtration when it is regular. So the action of \ $\hat{A}$ \ is well defined on a regular (a,b)-module and a regular (a,b)-module is canonically a regular \ $(\hat{a},b)-$module (see the important remark before lemma \ref{chgt. variable 0}).$\hfill \square$
 \end{enumerate}
  
\begin{lemma}\label{chgt. variable 5}
Let \ $E$ \ a simple pole $(\hat{a},b)-$module (resp. a regular $(\hat{a},b)-$module). Then for any change of variable \ $\theta$ \ the left \ $\hat{A}-$module \ $\theta_*(E)$ \ is a simple pole $(\hat{a},b)-$module (resp. a regular $(\hat{a},b)-$module).
\end{lemma}

\parag{Proof} As \ $E$ \ is free and finitely generated on \ $\C[[b]]$ \ the lemma \ref{chgt. variable 4} mplies that \ $\theta_*(E)$ \ is free and finitely generated on \ $\C[[b]]$. Then \ $\alpha.E \subset a.E \subset b.E = \beta.E$ \ gives that \ $\theta_*(E)$ \ has a simple pole when \ $E$ \ has a simple pole. As \ $\theta_*$ \ preserves \ $\hat{A}-$linear injection, it preserves regularity. $\hfill \blacksquare$\\

\begin{defn}\label{Bernstein 0}
Let \ $E$ \ be a regular $(\hat{a},b)-$module. We define the {\bf Bernstein polynomial} of \ $E$ \ as the minimal polynomial of the action of \ $-b^{-1}.a$ \ on the finite dimensional vector space \ $E^{\sharp}\big/b.E^{\sharp}$.
\end{defn}

\begin{prop}\label{Bernstein 1}
Let \ $E$ \ be a regular $(\hat{a},b)-$module and \ $\theta$ \ a change of variable. Then we have a canonical \ $\hat{A}-$linear isomorphism  \ $ \theta_*(E^{\sharp}) \to \theta_*(E)^{\sharp} $ \ and the Bernstein polynomial of \ $E$ \ and \ $\theta_*(E)$ \ co{\"i}ncide.
\end{prop}

\parag{Proof} As \ $\theta_*(E^{\sharp})$ \ is a simple pole $(\hat{a},b)-$module and \ $\theta_*(E)$ \ is regular, thanks to the previous lemma, the universal property of \ $ i : \theta_*(E) \to \theta_*(E)^{\sharp}$ \ implies, as the change of variable gives an \ $\hat{A}-$linear injection \ $\theta_*(i) : \theta_*(E) \to \theta_*(E^{\sharp})$, that there is an \ $\hat{A}-$linear factorization for \ $\theta_*(i) $
$$ \theta_*(E) \overset{i}{\to}  \theta_*(E)^{\sharp} \overset{j}{\to}  \theta_*(E^{\sharp}) .$$
It is then easy to see (using \ $\theta^{-1}$) \ that the injective map \ $\theta_*(i)$ \ has the same universal property than \ $i$. So \ $j$ \ is an isomorphism.\\
Now we have a natural isomorphism of vector spaces 
 $$ \theta_*(E)^{\sharp}\big/b.\theta_*(E)^{\sharp} \simeq \theta_*(E^{\sharp})\big/b.\theta_*(E^{\sharp}).$$
 To conclude, it is now enough to prove that the endomorphisms \ $b^{-1}.a$ \ on both sides are compatible with this isomorphism. This is a obvious consequence of the following fact : if \ $F$ \ is a simple pole $(\hat{a},b)-$module then \ $b^{-1}.a$ \ and \ $\beta^{-1}.\alpha$ \ induce the same endomorphism in \ $F\big/b.F$, because of  the relation
 $$ \beta^{-1}.\alpha = \theta'(a)^{-1}.(b^{-1}.a).(\theta(a)\big/a) $$
 and the fact that \ $a$ induces the \ $0$ \ map on \ $F\big/b.F$. $\hfill \blacksquare$

\begin{cor}\label{chgt. variable 6}
Let \ $E$ \ be rank 1 regular $(\hat{a},b)-$module. Then for any change of variable\ $\theta$ \ we have \ $\theta_*(E) \simeq E$.
\end{cor}

\parag{Proof} As any regular rank 1 $(\hat{a},b)-$module is isomorphic to \ $\hat{A}\big/(a - \lambda.b).\hat{A}$ \ for some \ $\lambda \in \C$ \ and as the corresponding Bernstein polynomial is \ $z + \lambda$, the previous proposition implies this result. $\hfill \blacksquare$\\

Remark that the isomorphism between \ $E_{\lambda}$ \ and \ $\theta_*(E_{\lambda})$ \ is unique up to a non zero constant, because \ $Aut_{\hat{A}}(E_{\lambda}) \simeq \C^*$. A precise description of this isomorphism will be given later on (see lemma \ref{rk 1} ).

\subsection{Change of variable for frescos and  themes.}

Our next proposition shows that the $(\hat{a},b)-$module of \ $[\lambda]-$primitive asymptotic expansions \ $\Xi_{\lambda}^{(N)}$ \ is stable by change of variable. The will use the following easy lemma ; the proof is left to the reader.

\begin{lemma}\label{tech.lin.}
Let \ $N $ \ be the endomorphism of \ $\C^n$ \ define by \ $N(e_j) = e_{j+1}$ \ for \ $j \in [1,n-1]$ \ and \ $N(e_n) = 0$. Then for any \ $\nu \in \C^*$ \ the  linear map
$$ f :  End(\C^n) \to End(\C^n) $$
given by \ $f(S) : = S.N - N.S + \nu.S $ \ is an isomorphism.
\end{lemma}

\begin{cor}\label{uniq. dev.}
Let \ $E$ \ be a simple pole $(\hat{a},b)-$module of rank \ $n$ \ such that the map induced by \ $a$ \ on \ $E\big/b^2.E$ \ is given in a suitable basis by
$$ a.e = b.( \lambda.Id_n + N).e $$
where \ $N$ \ is the principal nilpotent endomorphism of the previous lemma. Then there exists a \ $\C[[b]]-$basis \ $\varepsilon : = (\varepsilon_1, \dots, \varepsilon_n)$ \  of \ $E$ \ in which we have
$$ a.\varepsilon = b.( \lambda.Id_n + N).\varepsilon .$$
\end{cor}

\parag{Proof} We begin with a \ $\C[[b]]-$basis \ $e : = (e_1, \dots, e_n)$ \  of \ $E$ \ such that we have, thanks to our hypothesis,
$$ a.e = b.(\lambda.Id_n + N).e + b^2.Z.e $$
where \ $Z$ \ is in \ $End(\C^n) \otimes \C[[b]]$. We look for \ $S \in End(\C^n) \otimes \C[[b]]$ \ such that \ $S(0) = Id_n$ \ and that the \ $\C[[b]]-$basis \ $\varepsilon : = S.e$ \ satisfies our requirement. This gives the equation
\begin{equation*}
a.S.e = S.a.e + b^2.S'.e = S.b.(\lambda.Id_n + N).e + S.b^2.Z.e + b^2.S'.e = b.( \lambda.Id_n + N).S.e .
\end{equation*}
So we want to solve the equation
\begin{equation*}
S.N - N.S + b.S.Z + b.S' = 0 \tag{1}
\end{equation*}
Writting \ $S : = \sum_{\nu \geq 0}^{+\infty} \ S_{\nu}.b^{\nu}$ \ and \ $Z : = \sum_{\nu \geq 0}^{+\infty} \ Z_{\nu}.b^{\nu}$ \ we want to solve the recursion system
\begin{equation*}
S_{\nu}.N - N.S_{\nu} + \nu.S_{\nu} = - \sum_{j = 0}^{\nu-1} \ S_{\nu -j-1}.Z_j \quad \forall \nu \geq 1 \quad {\rm with} \quad \ S_0 = Id_n. \tag{2}
\end{equation*}
The previous lemma gives existence and uniqueness of the solution \ $S$. $\hfill \blacksquare$\\

\begin{prop}\label{chang. dev.}
Fix  \ $\lambda \in ]0,1] \cap \mathbb{Q}$ \ and \ $p \in \mathbb{N}$. Then for any change of variable \ $\theta$ \ the $(\hat{a},b)-$module \ $\theta_*(\Xi^{(p)}_{\lambda})$ \ is isomorphic to \ $\Xi^{(p)}_{\lambda}$.
\end{prop}

\parag{Proof} As the $(\hat{a},b)-$module \ $E : = \Xi^{(p)}_{\lambda}$ \   is a simple pole (a,b)-module of rank \ $n : = p+1$\ such that the map  induced by \ $a$ \ on \ $E\big/b^2.E$ \ is given in a suitable basis by
$$ a.e = b.( \lambda.Id_n + N).e $$
where \ $N$ \ is the principal nilpotent endomorphism of the lemma above, it is enough to prove that the map induced by \ $\alpha$ \ on \ $E\big/b^2.E \simeq E\big/\beta^2.E$ \ is again given by \ $\beta.( \lambda.Id_n + N)$. But we have \ $\alpha = \chi(\theta).a + b^2.E$ \ and \ $\beta = \chi(\theta).b + b^2.E$, which allow to conclude. $\hfill \blacksquare$.\\

Note that the compatibility of  changes of variable with direct sum implies that for any finite subset \ $[\Lambda] \subset \mathbb{Q}\big/\mathbb{Z}$ \ we have also
$$ \theta_*(\Xi_{[\Lambda]}^{(N)}\otimes V) \simeq \Xi_{[\Lambda]}^{(N)}\otimes V $$
for any finite dimensional complex vector space \ $V$.\\

So we obtain the following easy corollary.

\begin{cor} If \ $E$ \ is a geometric $(\hat{a},b)-$module, so is \ $\theta_*(E)$ \ for any change of variable \ $\theta$. Moreover, if \ $E$ \ is semi-simple, \ $\theta_*(E)$ \ is also semi-simple\footnote{It is easy to show that a  geometric $(\hat{a},b)-$module is semi-simple if and only if it admits an embedding in some \ $\Xi^{(0)}_{\Lambda}\otimes V$.} . 
\end{cor}

Our next  result  gives the the stability of the notions of fresco,  theme and of semi-simple geometric \ $(\hat{a},b)-$module  by any change of variable.

\begin{prop}\label{stability}
Let \ $E$ \ be a fresco  (resp. a theme) and let \ $\theta$ \ be any change of variable. Then \ $\theta_*(E)$ \ is a fresco   (resp. a theme). Moreover if \ $E$ \ is  \ $[\Lambda]-$primitive then \ $\theta_*(E)$ \ is also  \ $[\Lambda]-$primitive.
\end{prop}

The proof of the first statement is a consequence of the following lemma

\begin{lemma}\label{stab. gene.}
Let \ $E$ \ be a rank \ $k$ \  fresco and let \ $\theta \in \C[[a]]$ \ be a change of variable. Then a generator  \ $e$ \ of \ $E$ \ on \ $\hat{A}$ \ is also  a generator of \ $\theta_*(E)$ \ on \ $\hat{A}$.
\end{lemma}

\parag{Proof} Denote \ $\alpha_{\theta} : = \theta(a)$ \ and \ $\beta_{\theta}: = b.\theta'(a)$ \ the action of \ $a$ \ and \ $b$ \ on \ $E$ \ which define the \ $\hat{A}-$module \ $\theta_*(E)$. It is enough to show that \ $e, \alpha_{\theta}.e, \dots, \alpha_{\theta}^{k-1}.e$ \ is a \ $\C[[\beta_{\theta}]]-$basis of \ $\theta_*(E)$. And so it is enough to show that it induces a \ $\C-$basis of \ $\theta_*(E)\big/\beta_{\theta}.\theta_*(E) \simeq E\big/b.E$.\\
Assume that we have complex numbers \ $\mu_1, \dots, \mu_k$ \ and an element \ $x \in E$ \ such that
$$ \sum_{j=1}^k \ \mu_j.\alpha_{\theta}^j.e = \beta_{\theta}.x .$$
Then, as we have \ $\alpha_{\theta}^i = \chi(\theta)^i.a^i + a^{i+1}.\C[[a]] $ \ for each \ $i$, and \ $\alpha_{\theta}^k.E \subset \beta_{\theta}.E$, the matrix expressing the \ $(\alpha_{\theta}^i.e)_{i \in [0,k-1]}$ \ in the basis \ $(a^i.e)_{i \in [0,k-1]}$ \ is lower  triangular with the \ $\chi(\theta)^i \not= 0$ \ on the diagonal.  $\hfill \blacksquare$

\parag{Proof of the proposition \ref{stab. gene.}} In the case of a  theme, the result is an easy consequence of the fact that we have an isomorphism (see the remark after proposition \ref{chang. dev.}) \ $ \theta_*(\Xi_{\Lambda}^{(N)}) \simeq \Xi_{\Lambda}^{(N)}$. \\
The last  statement is consequence of the invariance of the Bernstein polynomial by change of variable. $\hfill \blacksquare$

\parag{Important Remark} The {\em primitive filtration}  of a fresco and its {\em principal Jordan-H{\"o}lder sequence} are covariant by any change of variable.  $\hfill \square$\\

\begin{cor}\label{chgt. variable struct.}
Let \ $E$ \ be a fresco and \ $\theta$ \ a change of variable. Then we have natural isomorphisms
$$ u_j : \theta_*(S_j(E)) \to S_j(\theta_*(E))\quad \forall j \geq 0 .$$
So the numbers \ $d(E)$ \ and \ $\delta(E)$ \ are invariant by any change of variable.
\end{cor}

\parag{Proof} As we know that any change of variable of a  semi-simple fresco is a   semi-simple fresco we have a natural maps \ $u_1 $. But extremality properties allow immediately to conclude that it is an isomorphism. The induction on \ $j$ \ is immediate as the functor \ $\theta_*$ \ is exact. $\hfill \blacksquare$\\

\subsection{Change variable for the dual.}

\begin{prop}\label{dual-change}
Let \ $E$ \ be a regular (a,b)-module and \ $\theta \in \C[[a]]$ \ a change of variable. Then we have a functorial isomorphism
$$ \theta_*(E^*) \to \theta(E)^* .$$
\end{prop}

\parag{proof} Let recall first that for any rank \ $p$ \ (a,b)-module \ $E$ \ we have a resolution
\begin{equation*}
 0 \to \A^p \overset{\square.M}{\longrightarrow} \A^p \to E \to  0  \tag{1}
 \end{equation*}
 where \ $M : = \id.a - N$, deduces from a \ $\C[[b]]-$basis \ $e$ \ and \ $a.e = N.e$ \ with \ $N$ \ a \ $(p,p)$ \ matrix with coefficients in \ $\C[[b]]$.
Then \ $E^*$ \ which is, by definition \ $Ext^1_{\A}(E,\A)$ \ endowed with the left \ $\A-$module structure associated with the right  \ $\A-$structure on \ $\A$ \ and the anti-automorphism of \ $\A$
$$ \eta : \A \to \A $$
defined by the conditions
\begin{enumerate}[i)]
\item \  $\eta(1) = 1$ \ and \ $\eta(x.y) = \eta(y).\eta(x) \quad \forall x, y \in \A$,
\item \ $\eta(a) = a$ \ and \ $\eta(b) = -b$.
\end{enumerate}
Now when \ $E$ \ is regular, \ $E^*$ \ is also regular and we may replace \ $\A$ \ by \ $\hat{A}$ \ in the previous considerations.\\
For any change of variable \ $\theta$ \ which is associated with the automorphism \ $\Theta$ \ of the algebra \ $\hat{A}$ \ defined by the conditions
\begin{enumerate}[i)]
\item \ $\Theta(1) = 1$ \ and \ $\Theta(x.y) = \Theta(x).\Theta(y) \quad \forall x, y \in \hat{A}$,
\item \ $\Theta(a) : = \theta(a)$ \ and \ $ \Theta(b) : = b.\theta'(a) $.
\end{enumerate}
Extending \ $\eta$ \ to \ $\hat{A}$ \  we have the relation \ $\eta\circ\Theta = \Theta\circ\eta $.\\
Now the exact sequence \ $(1)$ \ gives, after applying the functor \ $Hom_{\A}(- , \A)$ \ the exact sequence of right \ $\A-$modules
\begin{equation*}
0 \to \A^p \overset{\,^tM.\square}{\longrightarrow} \A^p \to  E^* \to 0  \tag{2}
\end{equation*}
and then the exact sequence of left \ $\A-$modules
\begin{equation*}
0 \to \A^p \overset{\square.\eta(^tM)}{\longrightarrow} \A^p \to E^* \to 0  \tag{3}
\end{equation*}
and after completion
\begin{equation*}
0 \to \hat{A}^p \overset{\square.\eta(^tM)}{\longrightarrow} \hat{A}^p \to E^* \to 0  \tag{$\hat{3}$}
\end{equation*}
Applying the functor \ $\Theta$ \ leads to the exact sequence
\begin{equation*}
0 \to \hat{A}^p \overset{\square.\Theta\circ\eta(^tM)}{\longrightarrow} \hat{A}^p \to \theta_*(E^*) \to 0  \tag{4}
\end{equation*}
Now applying the functor \ $\theta_*$ \ to the exact sequence \ $(\hat{1})$ \ gives the exact sequence of left \ $\hat{A}-$modules
\begin{equation*}
 0 \to \hat{A}^p \overset{\square.\Theta(M)}{\longrightarrow} \hat{A}^p \to \theta_*(E) \to  0  \tag{1bis}
 \end{equation*}
 and so the exact sequence of right \ $\hat{A}-$modules
 \begin{equation*}
0 \to \hat{A}^p \overset{\,^t\Theta(M).\square}{\longrightarrow} \hat{A}^p \to  \theta_*(E)^* \to 0  \tag{2bis}
\end{equation*}
and then the exact sequence of left \ $\hat{A}-$modules
\begin{equation*}
0 \to \A^p \overset{\square.\eta\circ\Theta(^tM)}{\longrightarrow} \A^p \to \theta_*(E)^* \to 0  \tag{3bis}
\end{equation*}
The relation \ \ $\eta\circ\Theta = \Theta\circ\eta $ \ allows to conclude. $\hfill \blacksquare$

\section{Holomorphic families of  frescos.}

\subsection{Holomorphic families and the principal  J-H. sequence.}

\begin{defn}\label{Hol. 1}
Let \ $X$ \ be a reduced complex space. We define the sheaf \ $\Xi_{\lambda, X}^{(N)}$ \ as the free sheaf of \ $\mathcal{O}_X[[b]]-$module with basis \ $e_0, \dots, e_N$. We define the action of \ $\hat{A}$ \ on this sheaf by the condition 1. 2. and 3. of the definition \ref{asympt. 0}. We obtain in this way a structure of left \ $\hat{A}-$module on this sheaf. We shall say that a sheaf on \ $X$ \ with these two compatible structures,  a \ $\mathcal{O}_X[[b]]-$module and a \ $\hat{A}-$left module, is a {\bf \ $\hat{A}_X-$sheaf}. 
\end{defn}

\parag{Fundamental example} For  any finite subset \ $\Lambda \subset \mathbb{Q}\big/\mathbb{Z}$ \ we define on 
 $$\Xi_{\Lambda, X}^{(N)} : = \oplus_{[\lambda_i] \in [\Lambda]} \quad \Xi_{\lambda_i,X}^{(N)}$$ 
 and for any finite dimensional complex vector space \ $V$ \ we define
 the actions of \ $a$ \ and \ $b$ \ on \ $\Xi_{\Lambda, X}^{(N)}\otimes V$ \ 
 via \ $ \oplus_{[\lambda_i] \in [\Lambda]} \big(a \otimes \id_V\big)$ \ and \ $ \oplus_{[\lambda_i] \in [\Lambda]} \big(b\otimes \id_V\big)$. This gives a structure of  \ $\hat{A}_X-$sheaf on \ $\Xi_{\Lambda, X}^{(N)}\otimes V$. $\hfill \square$

\parag{Remarks}
\begin{enumerate}
\item Of course the action of \ $\mathcal{O}_X$ \ and of \ $\hat{A}$ \ commutes on \ $\Xi_{\Lambda, X}^{(N)}\otimes V$.
\item For each point in \ $X$ \ we have an evaluation map
$$ ev_x : \Xi_{\Lambda, X}^{(N)}\otimes V \to \Xi_{\Lambda}^{(N)}\otimes V$$
which is \ $\hat{A}-$linear and surjective. For any \ $\hat{A}_X-$sub-sheaf \ $\mathbb{E}$ \ of \ $ \Xi_{\lambda, X}^{(N)}\otimes V $ \ we shall denote by \ $\mathbb{E}(x)$ \ its {\em image} by this evaluation map. \\
Note that \ $\mathbb{E}(x)$ \ is the image of the fiber \ $\mathbb{E}\big/\frak{M}_x.\mathbb{E}$ \  at \ $x$ \ of \ $\mathbb{E}$ \ in the fiber at \ $x$ \ of \ $\Xi_{\lambda, X}^{(N)}\otimes V$. So it is a quotient of its fiber at \ $x$ \ and the sheaf \ $\mathbb{E}$ \ may be free of finite type on \ $\mathcal{O}_X[[b]]$ \ and this quotient map may not be injective.
\end{enumerate}

Because of this last remark, the following easy lemma will be required later on.

\begin{lemma}\label{easy coh.}
Let \ $X$ \ be a reduced complex space and \ $\mathcal{F} \subset \mathcal{O}_X^p$ \ be a coherent sub-sheaf such that for each \ $x$ \ the evaluation map \ $ev_x : \mathcal{F} \to \C^p \simeq \mathcal{O}_X^p\big/\frak{M}_x.\mathcal{O}_X^p$ \ has rank r. Then \ $\mathcal{F}$ \ is locally free with  rank \ $r$ \ on \ $\mathcal{O}_X$ \ and is locally a direct factor of \ $\mathcal{O}_X^p$.
\end{lemma}

\parag{proof} Locally on \ $X$ \ consider  \ $M : \mathcal{O}_X^q \to \mathcal{O}_X^p $ \  a \ $\mathcal{O}_X-$map with image equal to  \ $\mathcal{F}$. Considering \ $M$ \ as a holomorphic  map to \ $End(\C^q,\C^p)$ \ we see that our hypothesis implies that the rank of \ $M(x)$ \ is constant and equal to \ $r$. So the image of \ $M$ \ is a rank \ $r$ \ vector bundle and \ $\mathcal{F}$, which is the sheaf of holomorphic sections of this rank \ $r$ \  vector bundle, is then locally free of rank \ $r$ \ and  is locally a direct factor of \ $\mathcal{O}_X^p$. $\hfill \blacksquare$\\

\begin{defn}\label{Hol. 2}
A section \ $\varphi$ \ of the sheaf \ $\Xi_{\Lambda, X}^{(N)}\otimes V$ \ will be called {\bf \ $k-$admissible} when the sub$-\hat{A}_X-$sheaf generated by \ $\varphi$ \ in \ $\Xi_{\Lambda, X}^{(N)}\otimes V$ \ is a free rank \ $k$ \ $\mathcal{O}_X[[b]]-$module with basis \ $\varphi, a.\varphi, \dots, a^{k-1}.\varphi$ \ such that, for each \ $x \in X$,  $\mathbb{E}(x)$ \ is a rank \ $k$ \ fresco, where we denote  \ $\mathbb{E}(x)$ \ the image of \ $\mathbb{E}$ \ by the evalution map \ $ev_x$.\\
A \ $\hat{A}_X-$sheaf \ $\mathbb{E}$ \ on \ $X$ \ will be called an {\bf holomorphic family of rank \ $k$ \  frescos} parametrized by \ $X$, when the sheaf \ $\mathbb{E}$ \  is locally isomorphic on \ $X$ \  to a \ $\hat{A}_X-$sheaf generated by an admissible section of \ $ \Xi_{\Lambda, X}^{(N)}\otimes V $, for some \ $N$, some finite set \ $[\Lambda] \subset \mathbb{Q}\big/\mathbb{Z}$ \  and some finite dimensional complex vector space \ $V$.\\
\end{defn}

Of course the sheaf \ $\mathbb{E}$ \ on \ $X$ \ defines the family \ $(\mathbb{E}(x))_{ x \in X}$ \ of rank \ $k$ \ $[\Lambda]-$primitive frescos (see the remark 2 below). And conversely, such a family will be holomorphic when it is given  by the evaluations maps \ $ev_x$ \ from a \ $\hat{A}_X-$sheaf \ $\mathbb{E}$ \ as in the previous definition.

\parag{Remarks}
\begin{enumerate}
\item The sub$-\hat{A}_X-$sheaf generated by a \ $k-$admissible section \ $\varphi$ \ of \ $\Xi_{\Lambda, X}^{(N)}\otimes V$ \ is stable by \ $\C[[a]]$ \ by definition. So we may write locally
$$ a^k.\varphi = \sum_{j=0}^{k-1} \ T_j.a^j.\varphi $$
where \ $T_j$ \ are unique (local) section of \ $\mathcal{O}_X[[b]]$. So there is unique monic degree \ $k$ \ polynomial \ $P : = a^k - \sum_{j=0}^{k-1} \ T_j.a^j$ \ in \ $\mathcal{O}_X[[b]][a]$ \ such that \ $P.\varphi = 0$. An easy argument of division shows that the left ideal of \ $\hat{A}_X$ \ generated by \ $P$ \ is the annihilator of \ $\varphi$.
\item When \ $\mathbb{E}$ \ is locally free rank \ $k$ on \ $\mathcal{O}_X[[b]]$ \ its fiber at \ $x$ \ $\mathbb{E}\big/\frak{M}_x.\mathbb{E}$ \ is a free rank \ $k$ \ $\C[[b]]-$module. In the case \ $\mathbb{E}$ \ is a holomorphic family of rank \ $k$ \ $[\Lambda]-$primitive frescos parametrized by \ $X$, the evaluation map 
 $$\mathbb{E}\big/\frak{M}_x.\mathbb{E} \to \mathbb{E}(x)$$ 
  associated to a local  isomorphism \ $\mathbb{E} \simeq \hat{A}_X.\varphi$, where \ $\varphi$ \ is a $k-$admissible section of \ $ \Xi_{\Lambda, X}^{(N)}\otimes V $, for some \ $N$, and some finite dimensional complex vector space \ $V$ \ and to the corresponding evaluation map \ $ev_x$, is an isomorphism because it is a surjective map between two rank \ $k$ \ free \ $\C[[b]]-$modules. So we shall  identify \ $\mathbb{E}\big/\frak{M}_x.\mathbb{E}$ \ and \ $\mathbb{E}(x)$ \ in this case, and use the notation \ $\mathbb{E}(x)$ \ for the fiber at \ $x$ \ of \ $\mathbb{E}$.
\item As long as we work inside a simple pole \ $\hat{A}_X-$sheaf (i.e. \ $a.\mathcal{F}\subset b.\mathcal{F}$) \ which is a locally free finite type \ $\mathcal{O}_X[[b]]$ \ module (as \ $\Xi_{\Lambda,X}^{(N)}\otimes V$), we dont need to consider the sheaf algebra \ $\hat{A}_X$ \ to define the \ $\hat{A}_X-$structure :  the \ $\hat{A}_X-$structure is completely defined by the \ $\mathcal{O}_X[[b]]-$structure and the action of \ $a$ \ (i.e. the left module structure on the sub-algebra \ $\mathcal{O}_X[[b]][a]$) because the \ $\mathcal{O}_X[[b]]-$completion implies the \ $\mathcal{O}_X[[a]]$ \ completion.\\
\item Our definition of an \ $1-$admissible section  is not satisfied by the following :
$$\varphi : =  z.s^{\lambda_1-1} \in \Xi_{\lambda,\C}^{(0)}$$
where \ $z$ \ is the coordinate on \ $X : = \C$ \ and \ $\lambda_1 \in 1+ \mathbb{Q}^{+*}$ \ is in \ $[\lambda]$ \ because the evaluation at \ $z = 0$ \ of the \ $\hat{A}_{\C}-$sheaf generated by \ $\varphi$ \ is \ $0$. But of course the "abstract" sheaf \ $\hat{A}_{\C}.\varphi$ \ is isomorphic to the subsheaf of \ $\Xi_{\lambda,\C}^{(0)}$ \ generated by the \ $1-$admissible section \ $\psi : = s^{\lambda_1-1}$ \ and so it defines an holomorphic family of rank \ $1$ \ frescos parametrized by \ $\C$ \ via the isomorphism define by  \ $\varphi \mapsto \psi$.\\
\end{enumerate}

\begin{lemma}\label{Hol. 3}
Let \ $\mathbb{E}$ \ be a holomorphic family of rank \ $k$ \ $[\Lambda]-$primitive frescos. Then for each \ $x$ \ in \ $X$ \ the \ $\hat{A}-$module \ $\mathbb{E}(x)$ \ is a rank \ $k$ \ $[\Lambda]-$primitive fresco and  the Bernstein polynomial of \ $\mathbb{E}(x)$ \ is locally constant on \ $X$.
\end{lemma}

\parag{Proof} It is enough to prove the result for the \ $\hat{A}_X-$sheaf generated by an admissible section \ $\varphi$ \ of the sheaf \ $\Xi_{\lambda, X}^{(N)}\otimes V$. Then write
\begin{equation*}
 a^k.\varphi = \sum_{j=1}^k \ S_j(b,x).a^{k-j}.\varphi  \tag{@}
 \end{equation*}
 
where  \ $S_1, \dots, S_k$ \ are sections of the sheaf \ $\mathcal{O}_X[[b]]$. For each given \ $x$, we know that \ $\mathbb{E}(x)$ \ is a fresco of rank \ $ k$. This implies that the valuation in \ $b$ \ of \ $S_j$ \ is \ $\geq j$. Then we have \ $S_j(b,x) = b^j.s_j(x) + b^{j+1}.\mathcal{O}_X[[b]] $ \ where \ $s_j$ \ is a holomorphic function on \ $X$. As the Bernstein element\footnote{the Bernstein element  of a fresco \ $E$ \  with fundamental invariants \ $\lambda_1, \dots, \lambda_k$ \ is the element \ $P_E : = (a - \lambda_1.b).\dots (a - \lambda_k.b)$ \  in \ $\hat{A}$ ; see [B.09].} of \ $\mathbb{E}(x)$ \ is equal to \ $a^k + \sum_{j=1}^k \ s_j(x).b^j.a^{k-j} $ \ and has rational coefficients, we conclude that \ $s_j$ \ is locally constant on \ $X$. $\hfill \blacksquare$\\

\begin{prop}\label{primitive part}
Let \ $\mathbb{E}$ \ be a holomorphic family of \ $[\Lambda]-$primitive frescos parametrized by a reduced complex space \ $X$ \  and let \ $[M] \subset [\Lambda]$. Then there exists a natural \ $\hat{A}_X-$sheaf \ $\mathbb{E}[M] \subset \mathbb{E}$ \ which is normal in \ $\mathbb{E}$ \ and defined a holomorphic family of \ $[M]-$primitive frescos such for each \ $x \in X$ \ we have
$$ \mathbb{E}(x)[M] = \mathbb{E}[M](x) .$$
\end{prop}

This proposition is an easy consequence of the following lemma. $\hfill \blacksquare$ 

\begin{lemma}\label{easy}
Let \ $\varphi$ \ be a \ $k-$admissible section of some sheaf \ $ \Xi_{\Lambda,X}^{(N)} \otimes V$ \ and \ $P : = (a - \lambda_1.b).S_1^{-1} \dots S_{k-1}^{-1}.(a - \lambda_k.b) $ \ a generator of the annihilator of \ $\varphi$. Then \ $(a - \lambda_k.b).\varphi$ \ is \ $(k-1)-$admissible and defines a normal corank \ $1$ \ holomorphic sub-family of the holomorphic family defined by \ $\mathbb{E}$.
\end{lemma}

\parag{proof} As \ $(a - \lambda_k.b).\varphi(x)$ \ generates  for each \ $x \in X$ \ a rank \ $(k-1)-$fresco which is normal in \ $\A.\varphi(x)$, the assertion is obvious. $\hfill \blacksquare$\\

\begin{thm}\label{J-H. hol.}
Let \ $X$ \ be a reduced complex space and \ $\mathbb{E}$ \ an holomorphic family of rank \ $k$ \  frescos parametrized by \ $X$. Then for each \ $j$ \ in \ $[1,k]$ \ there exists a sub$-\hat{A}_X-$sheaf \ $\mathbb{F}_j \subset \mathbb{E}$ \ which is an holomorphic family of rank \ $j$ \ frescos parametrized by \ $X$ \ such that for each \ $x$ \ in \ $X$ \ the sub-modules \ $\mathbb{F}_j(x), j \in [1,k]$, give the principal J-H. sequence of the fresco \ $\mathbb{E}(x)$.\\
Moreover, for each \ $j \in [1,k]$ \ the quotient sheaf \ $\mathbb{E}\big/\mathbb{F}_j$ \ is an holomorphic family of frescos.
\end{thm}

\parag{proof} In fact we shall prove by induction on \ $k \geq 1$ \ that for any \ $k-$admissible section \ $\varphi$ \ of the sheaf \ $\Xi_{\Lambda, X}^{(N)}\otimes V$ \ corresponding to an holomorphic family of rank \ $k$ \ frescos with fundamental invariants \ $\lambda_1, \dots, \lambda_k$, there exists locally on \ $X$, invertible sections \ $S_1, \dots, S_k$ \ of the sheaf \ $\mathcal{O}_X[[b]]$ \ such that
\begin{equation*}
 (a - \lambda_1.b).S_1^{-1}.(a - \lambda_2.b) \dots S_{k-1}^{-1}.(a - \lambda_k.b).S_k^{-1}.\varphi \equiv 0. \tag{@@}
 \end{equation*}
Then it will not be difficult to prove that for each \ $j \in [0,k-1]$ \ the section
$$ \varphi_j : = (a - \lambda_{k-j+1}.b).S_{k-j+1}^{-1} \dots S_{k-1}^{-1}.(a - \lambda_k.b).S_k^{-1}.\varphi $$
is \ $(k-j)-$admissible and generates the sheaf \ $\mathbb{F}_{k-j}$ \ corresponding to the family \ $\mathbb{F}_{k-j}(x), x \in X$ \ locally.\\

The case \ $k = 1$ \  is contained in the following  lemma. 

\begin{lemma}\label{Hol. rank 1}
Let \ $X$ \ be a connected reduced complex space and let  \ $\varphi$ \ be an $1-$admissible section of the sheaf \ $\Xi^{(N)}_{\Lambda,X}\otimes V$. Then there exists locally on \ $X$ \ an invertible section \ $S$ \ of the sheaf \ $\mathcal{O}_X[[b]]$ \ and a non vanishing holomorphic map \ $v_0 : X \to  V$ \ such that 
$$ \varphi = S(b,x).s^{\lambda_1-1} \otimes v_0(x)$$
where \ $\lambda_1$ \ is a positive rational number in \ $[\Lambda] \subset \mathbb{Q}\big/\mathbb{Z}$.
\end{lemma}

\parag{proof} By definition the \ $\mathcal{O}_X[[b]]-$module \ $\mathbb{E}$ \  generated by \ $\varphi$ \ is free with rank \ $1$ \ and stable by \ $a$. Write \ $a.\varphi = T(b,x).\varphi$. As for each \ $x$ \ in \ $X$ \ the fresco \ $\mathbb{E}(x)$ \ has rank \ $1$, it has a simple pole and \ $T(b,x)$ \ has valuation \ $\geq 1$ \ in \ $b$. Moreover if we put \ $T(b,x) = t(x).b + b^2.T_1(b,x)$ \ the Bernstein element  of \ $\mathbb{E}(x)$ \ is \ $a - t(x).b$. So, as in the lemma \ref{Hol. 3}, this implies that the holomorphic function \ $t$ \ is constant on \ $X$ ; let \ $\lambda_1$ \ the positive  rational value of \ $t$. We may write 
 $$a.\varphi =  b.\big(\lambda_1 +  b.T_1(b,x)\big).\varphi .$$
 Now the section \ $U : = \lambda_1 +  b.T_1(b,x)$ \ of \ $\mathcal{O}_X[[b]]$ \ is invertible and  we may solve the differential equation (with holomorphic parameter in \ $X$)
 $$b.S' + S.(U - \lambda_1) = 0  \quad S(0,x) \equiv 1 $$
 where \ $S$ \ is a section of \ $\mathcal{O}_X[[b]]$, and obtain the new generator \ $\psi : = S.\varphi$ \ for \ $\mathbb{E}$ \ with the following property
 \begin{align*}
 & a.S.\varphi = S.a.\varphi + b^2.S'.\varphi = S.b.U.\varphi + b^2.S'.\varphi =   \lambda_1.b.S.\varphi \quad {\rm and \  so} \\
 &a.\psi = \lambda_1.b.\psi 
  \end{align*}
  and this implies that there exists an holomorphic function \ $v_0 : X \to   V $ \ such that \ $\psi = s^{\lambda_1-1}\otimes v_0$. To conclude the proof it is enough to remark that if \ $v_0(x) = 0$ \ then \ $\mathbb{E}(x) = \{0\}$ \ contradicting our rank \ $1$ \ assumption. So \ $v_0$ \ does not vanish. $\hfill \blacksquare$\\
  
  So now assume that we have proved our assertion \ $(@@)$ \ for \ $k-1 \geq 1$, and consider now a \ $k-$admissible section \ $\varphi$ \ of the sheaf \ $ \Xi_{\Lambda, X}^{(N)}\otimes V$ \ corresponding to a holomorphic family of rank \ $k$ \ frescos with fundamental invariants \ $\lambda_1, \dots, \lambda_k$. Consider the subsheaf \ $K : = \big(\mathcal{O}_X.s^{\lambda_1-1}\otimes V\big) \cap \mathbb{E}$ \ of \ $\mathbb{E}$. It is \ $\mathcal{O}_X-$coherent as it is contained in the \ $\mathcal{O}_X-$coherent subsheaf  of \ $\Xi^{(0)}_{\Lambda,X} \otimes V$ \ where the degree in \ $b$ \ is bounded by some  \ $n \gg 1$. As \ $K(x) = F_1(\mathbb{E}(x)) \cap Ker\,(a - \lambda_1.b)$ \ has dimension 1 for each \ $x$ \ in \ $X$, it is a locally free rank \ $1$ \ $\mathcal{O}_X-$module which is a direct factor, thanks to lemma \ref{easy coh.} ; then we may find, at least locally, an holomorphic non vanishing  section \ $v_0 : X \to V$ \ which is a local basis for it. The subsheaf \ $\mathcal{O}_X[[b]].\psi$, where \ $\psi : = s^{\lambda_1-1}\otimes v_0$,  has a  simple pole \ because $a.\psi = \lambda_1.b.\psi$. Then \ $\psi$ \ is admissible an define the holomorphic family \ $\mathbb{F}_1(x)_{x \in X}$. \\
  
To prove that \ $\mathbb{E}\big/\mathbb{F}_1$ \ is a holomorphic family, remark first that if \ $\chi$ \ is a section of \ $\mathbb{E}$ \ of the form \ $\chi : = S(b,x).s^{\mu-1}\otimes v_0 $ \ for some \ $\mu \in [\lambda_1]$, where we assume that the generic valuation of \ $S$ \ in \ $b$ \ is \ $0$, then \ $\chi$ \ is a section of \ $\mathcal{O}_X[[b]].\psi$ \ because  we have \ $\mu  \geq \lambda_1$ : indeed,  near a generic point \ $x$ \ the section \ $S$ \ will be invertible and we find \ $s^{\mu-1}\otimes v_0$ \ in \ $\mathbb{E}$. But for each \ $x$ \ in \ $X$ \ there is no non zero element in \ $\mathbb{E}(x)$ \ which is annihilated by\ $a - \mu.b$ \ with \ $\mu < \lambda_1$. This proves our claim.\\
Now write locally  \ $\mathcal{O}_X\otimes V = (\mathcal{O}_X \otimes W) \oplus \mathcal{O}_X.v_0 $ \ where \ $W$ \ is an hyperplane in \ $V$, and consider the map
$$ \sigma : \Xi_{\Lambda,X}^{(N)} \otimes V \to  \Xi_{\Lambda,X}^{(N)} \otimes V $$
given by the identity map on \ $\Xi_{\Lambda',X}^{(N)}\otimes V$, where we define \ $[\Lambda'] : = [\Lambda]\setminus [\lambda_1]$, and also on  \ $ \Xi_{\lambda_1,X}^{(N)} \otimes W$ \ and given on \ $\Xi_{\lambda_1,X}^{(N)} \otimes v_0$ \ by the quotient map by \ $ \Xi_{\lambda_1,X}^{(0)} \otimes v_0$ :
$$ \Xi_{\lambda_1,X}^{(N)} \otimes v_0 \to  \Xi_{\lambda_1,X}^{(N-1)} \otimes v_0$$
composed with the obvious inclusion 
$$  \Xi_{\lambda_1,X}^{(N-1)} \otimes v_0 \to  \Xi_{\lambda_1,X}^{(N)} \otimes v_0.$$
As we  proved above  that the intersection of the kernel of \ $\sigma$ \ with \ $\mathbb{E}$ \ is \ $\mathbb{F}_1$ \  the quotient sheaf \ $\mathbb{E}\big/\mathbb{F}_1$ \ is isomorphic to the subsheaf of \ $\Xi_{\Lambda,X}^{(N)} \otimes V$ \ generated by the section \ $ \sigma(\varphi)$.\\
Moreover we have\ $\hat{A}_X.\sigma(\varphi)(x) \simeq \mathbb{E}(x)\big/\mathbb{F}_1(x)$ \ for each \ $x \in X$. We shall prove that the section  \ $\sigma(\varphi)$ \ is \ $(k-1)-$admissible.\\
 First we shall prove  that the sections \ $\sigma(\varphi), a.\sigma(\varphi), \dots, a^{k-2}.\sigma(\varphi)$ \ are locally free on \ $\mathcal{O}_X[[b]]$ \ in the sheaf \ $ \Xi_{\Lambda,X}^{(N)} \otimes V $. This is equivalent to show that the corresponding classes are \ $\mathcal{O}_X-$free modulo \ $b$. We shall use the following lemma.

\begin{lemma}\label{tech}
In the situation above write
$$ \psi = s^{\lambda_1-1} \otimes v_0 = \sum_{j=0}^{k-1} \ U_j.a^j.\varphi \  ; $$
then the \ $b-$valuation of  \ $U_j $ \ is at least \ $k-j-1$, for \ $j \in [0,k-1]$.
\end{lemma}

\parag{proof} Write \ $(a- \lambda_1.b).\psi = 0$. It gives the equations
\begin{align*}
& U_{k-1}.S_j + b^2.U'_j - \lambda_1.b.U_j + U_{j-1} = 0 
\end{align*}
where we have written  \ $a^k.\varphi = \sum_{j=0}^{k-1} S_j.a^j.\varphi$ \ and used the convention \ $U_{-1} \equiv 0 $.\\
Recall that the \ $b-$valuation of \ $S_j$ \ is at least \ $k-j$. Then an easy induction conclude the proof.$\hfill \blacksquare$

\parag{End of the proof of the theorem \ref{J-H. hol.}} Now  \ $\psi = U_{k-1}.a^{k-1}.\varphi + b.\mathbb{E}$ \ does not belong to \ $b.\mathbb{E}$ \ because \ $\mathbb{F}_1(x)$ \ is normal in \ $\mathbb{E}(x)$ \ for each \ $x$. So \ $U_{k-1}$ \ is invertible in \ $\mathcal{O}_X[[b]]$ \ and we have
$$ U_{k-1}.a^{k-1}.\sigma(\varphi) = - \sum_{j=0}^{k-2} \ U_j.a^j.\sigma(\varphi).$$
So we easily deduce that \ $\sigma(\varphi), a.\sigma(\varphi), \dots, a^{k-2}.\sigma(\varphi)$ \ is locally free on \ $\mathcal{O}_X[[b]]$ \ in the sheaf \ $ \Xi_{\Lambda,X}^{(N)} \otimes V $ \ and generate a subsheaf stable by \ $a$. So \ $\sigma(\varphi)$ \ is \ $(k-1)-$admissible with corresponding fundamental invariants \ $\lambda_2, \dots, \lambda_k$ \ and the induction hypothesis gives \ $S_2, \dots, S_k$ \ invertible in \ $\mathcal{O}_X[[b]]$ \ such that
$$ (a - \lambda_2.b).S_2^{-1} \dots S_{k-1}^{-1}.(a - \lambda_k.b).S_k^{-1}.\sigma(\varphi) \equiv 0 $$
in \ $\mathbb{E}\big/\mathbb{F}_1$. So we get
$$ (a - \lambda_2.b).S_2^{-1} \dots S_{k-1}^{-1}.(a - \lambda_k.b).S_k^{-1}.\varphi = T.\psi $$ with \ $T \in \mathcal{O}_X[[b]].\psi .$
This equality modulo \ $b$ \ gives an invertible element \ $f$ \  of \ $\mathcal{O}_X$ \ such that \ $f.a^{k-1}.\varphi = T(0).\psi \quad modulo \ b.\mathbb{E} $. \\
But we already know that 
 $\psi = U_{k-1}.a^{k-1}.\varphi + b.\mathbb{E}$, with \ $U_{k-1}(0)$ \ invertible in \ $\mathcal{O}_X$, so we conclude that we have 
$$ (f - T(0).U_{k-1}(0))a^{k-1}.\varphi \in b.\mathbb{E} .$$
As we know that for each \ $x \in X$ \ we have \ $a^{k-1}.\varphi(x) \not\in b.\mathbb{E}(x)$ \ it implies  that \ $f = U_{k-1}(0).T(0)$. So \ $T$ \ is an invertible element in \ $\mathcal{O}_X[[b]]$ \ and we may define \ $S_1 : = T$ \ to conclude our induction.\\

Now to complete the proof of the theorem it is enough to prove that for each \\
 $j \in [1,k-1]$ \ the section   \ $\varphi_j$ \ is \ $(k-j)-$admissible. But it is clear that, modulo \ $b.\mathbb{E}$, the classes of \ $\varphi_j, a.\varphi_j, \dots, a^{k-j-1}.\varphi_j$ \ co{\"i}ncide with \ $a^j.\varphi, \dots, a^{k-1}.\varphi$. So they are independant on \ $\mathcal{O}_X[[b]]$. The identity
$$ (a - \lambda_1.b).S_1^{-1} \dots (a - \lambda_{k-j}.b)S_{k-j}^{-1}.\varphi_j \equiv 0 $$
shows that \ $a^{k-j}.\varphi_j$ \ is in the \ $\mathcal{O}_X[[b]]-$module with basis \ $\varphi_j, a.\varphi_j, \dots, a^{k-j-1}.\varphi_j$, and then \ $\varphi_j$ \ is \ $(k-j)-$admissible for each \ $j$.$\hfill \blacksquare$\\

\begin{thm}\label{chgt. var. hol.}
Let \ $\mathbb{E}$ \ be a holomorphic family of rank \ $k$ \ $[\Lambda]-$primitive frescos parametrized by a reduced  complex space \ $X$\  and let \ $\theta \in \C[[a]]$ \ be a change of variable. Then the family \ $\theta_*(\mathbb{E}(x)), x \in X$ \ is holomorphic.
\end{thm}

\parag{Proof} The problem is local on \ $X$ \ and we may assume that \ $\mathbb{E}$ \ is generated by a \ $k-$admissible section \ $\varphi$ \ of some sheaf \ $\Xi_{\Lambda, X}^{(N)}\otimes V $. As the \ $\hat{A}_X-$module \ $\theta_*(\Xi_{\Lambda, X}^{(N)}\otimes V)$ \ is isomorphic to \ $\Xi_{\Lambda, X}^{(N)}\otimes V $, thanks to the remark following the  proposition \ref{chang. dev.}, it is enough to prove that \ $\theta_*(\varphi)$ \ is a \ $k-$admissible section generating \ $\theta_*(\mathbb{E})$. Then  it is enough to prove that the images in 
$$ \mathbb{E}\big/b.\mathbb{E} = \mathbb{E}\big/\beta.\mathbb{E} $$
of  \ $\varphi, \alpha.\varphi, \dots, \alpha^{k-1}.\varphi$ \  are \ $\mathcal{O}_X-$free, where \ $\mathbb{E} : = \hat{A}_X.\varphi$. And this is an easy consequence of the fact that the images of \ $\varphi, a.\varphi, \dots, a^{k-1}.\varphi$ \ are  \ $\mathcal{O}_X-$free by the definition of a \ $k-$admissible section.  $\hfill \blacksquare$\\

\subsection{Basic theorems on holomorphic families.}

In this section we shall prove some general stability results for the notion of holomorphic family of frescos.

\begin{thm}\label{suite exacte}
Let \ $X$ \ be a reduced complex space and consider an exact sequence of \ $\hat{A}_X-$sheaves :
$$ 0 \to \mathbb{F} \to \mathbb{E} \to \mathbb{G} \to 0 .$$
Assume that \ $ \mathbb{F}$ \ and \ $ \mathbb{G}$ \ are holomorphic families of frescos parametrized by \ $X$. Assume also  that for each \ $x \in X$ \ we have an exact sequence of frescos given by the fibers at \ $x$
$$0 \to \mathbb{F}(x) \to \mathbb{E}(x) \to \mathbb{G}(x) \to 0 $$ Then the family \ $\mathbb{E}$ \ is holomorphic.
\end{thm}

\parag{Proof} We shall show that it is enough to prove the case where \ $\mathbb{F}$ \ has rank \ $1$.
Assume that \ $\mathbb{F}$ \ has rank \ $k-1 \geq 2$ \ and that the theorem has been proved for rank \ $\leq k-2$ \ in the first sheaf of the exact sequence. Let \ $\mathbb{F}_1$ \ be the sheaf given by the rank \ $1$ \ term in the principal J-H. sequence of \ $\mathbb{F}$. Then we have an exact sequence of \ $\hat{A}_X-$sheaves :
$$  0 \to \mathbb{F}\big/\mathbb{F}_1 \to \mathbb{E}\big/\mathbb{F}_1 \to \mathbb{G} \to 0 $$
and now the first sheaf has rank \ $k-2$ \ and is an holomorphic family thanks to \ref{J-H. hol.}. To apply the induction hypothesis we have to check that \ $(\mathbb{E}\big/\mathbb{F}_1)(x) = \mathbb{E}(x)\big/\mathbb{F}_1(x)$ \ for each \ $x \in X$. But this is a consequence of the exactness of the  sequence
$$  \mathbb{F}_1(x) \to \mathbb{E}(x) \to (\mathbb{E}\big/\mathbb{F}_1)(x) \to 0 $$
(the tensor product is right exact) and the fact that the first map has to be injective be cause the ranks force its kernel to be torsion. So \ $(\mathbb{F}\big/\mathbb{F}_1)(x) = \mathbb{F}(x)\big/\mathbb{F}_1(x)$ \ and we have an exact sequence
$$ 0 \to (\mathbb{F}\big/\mathbb{F}_1)(x) \to (\mathbb{E}\big/\mathbb{F}_1)(x) \to \mathbb{G}(x) \to 0 $$
for each \ $x \in X$ \ and  the induction hypothesis implies the holomorphy of the family \ $ \mathbb{E}\big/\mathbb{F}_1$. Then applying the rank \ $1$ \ case to the exact sequence of \ $\hat{A}_X-$sheaves :
$$ 0 \to \mathbb{F}_1 \to \mathbb{E} \to \mathbb{E}\big/\mathbb{F}_1 \to 0 $$
we conclude that \ $\mathbb{E}$ \ is a holomorphic family.\\

Now we consider the case where \ $\mathbb{F}$ \ is rank \ $1$.  Let \ $k$ \ be the rank of \ $\mathbb{E}$. The assertion is local on \ $X$ \ so we may assume that \ $\mathbb{G}$ \ is generated by a \ $(k-1)-$admissible section \ $\psi$ \ of some sheaf \ $\Xi_{\Lambda,X}^{(N)} \otimes V$. We may find (locally on \ $X$) \ a section \ $\varphi$ \ of \ $\mathbb{E}$ \ which lift \ $\psi$. Let \ $Q$ \ be a section of the sheaf \ $\sum_{j=0}^{k-1} \ \mathcal{O}_X[[b]].a^j $, monic of degree \ $k-1$ \ in \ $a$ \ such that \ $\mathbb{G} \simeq \hat{A}_X\big/\hat{A}_X.Q $. So \ $Q$ \ generates the annihilator\footnote{ see the remark i) following the definition \ref{Hol. 2}.} of \ $\psi$ \ in \ $\hat{A}_X$.
Then \ $Q.\varphi$ \ is a section of \ $\mathbb{F}$ \ which generates \ $\mathbb{F}(x)$ \ for each \ $x \in X$, and as we may assume, at least locally on \ $X$, that \ $\mathbb{F} \simeq \mathcal{O}_X[[b]].s^{\lambda-1}$, there exists an invertible section \ $S_1$ \ of the sheaf \ $\mathcal{O}_X[[b]]$ \ such that \ $P : = (a - \lambda_1.b).S_1^{-1}.Q$ \ annihilates \ $\varphi$. So we may (locally) assume that \ $\varphi$ \ is a section of some sheaf \ $\Xi_{\Lambda,X}^{(N')}\otimes V'$ \ such that \ $\mathbb{E}(x) \simeq \hat{A}.\varphi(x)$ \ for all \ $x$.\\
To conclude, it is then enough to show that \ $\varphi, a.\varphi, \dots, a^{k-1}.\varphi$ \ is a \ $\mathcal{O}_X[[b]]-$basis of \ $\mathbb{E}$. As this is true for each given \ $x$ \ in \ $X$, the only point to prove is the fact that they generate the sheaf \ $\mathbb{E}$. Let \ $\sigma$ \ be a local section of  \ $\mathbb{E}$. We may write, locally on \ $X$, the image of \ $\sigma$ \ by the quotient map \ $\pi :  \mathbb{E} \to \mathbb{G}$ \ as 
$$  \pi(\sigma) = \sum_{j=0}^{k-2} \ T_j.a^j.\psi $$
where \ $T_j$ \ are local sections of \ $\mathcal{O}_X[[b]]$. Then \ $\sigma - \sum_{j=0}^{k-2} \ T_j.a^j.\varphi$ \ is a section of \ $\mathbb{F}$ \ and it may be written as \ $T.Q.\varphi$. This gives the conclusion. $\hfill \blacksquare$\\

Let me recall that for a regular  (a,b)-module \ $E$ \  the dual \ $E^*$ \ may be also defined as follows (comparaison with the section 2.5 is not obvious ; see [K.11] ) :\\
as a \ $\C[[b]]-$module \ $E^*$ \ is \ $Hom_{\C[[b]]}(E,E_0)$ \ where \ $E_0 : = \A\big/\A.a$, and the action of \ $a$ \ on \ $E^*$ \ is given by
\begin{equation*}
 (a.l)(x) = a.l(x) - l(a.x) \tag{*}
 \end{equation*}
for \ $l \in E^*$ \ and \ $x \in E$. We may also define directly in this way the \ $\delta-$dual \ $E^*\otimes E_{\delta}$ \ for any complex number \ $\delta$ \ by the same formula for \ $a$ \ acting on the free finite type \ $\C[[b]]-$module \ $Hom_{\C[[b]]}(E,E_{\delta})$.\\
For a holomorphic family we define the \ $\hat{A}_X-$sheaf \ $\mathbb{E}^*\otimes E_{\delta}$ \ as the \ $\mathcal{O}_X[[b]]-$module \ $Hom_{\mathcal{O}_X[[b]]}(\mathbb{E},E_{\delta})$ \ where \ $E_{\delta}$ \ is now the rank \ $1$ \ free \ $\mathcal{O}_X[[b]]-$module generated by \ $e_{\delta}$ \ on which we define \ $a.e_{\delta} = \delta.b.e_{\delta}$. In order to have a \ $\hat{A}_X-$sheaf the action of \ $a$ \ is defined by the formula \ $(^*)$. 

\parag{Remark} When \ $X$ \ is connected, for \ $\delta$ \ big enough, \ $[\mathbb{E}^*\otimes E_{\delta}](x)$ \ is a fresco for each \ $x \in X$. More precisely, if \ $\mathbb{E}$ \ is  a holomorphic \ $[\lambda]-$primitive family the inequality \ $\delta > \lambda_k + k-1$ \ implies that \ $\mathbb{E}(x)^*\otimes E_{\delta}$ \ is a fresco for each \ $x \in X$.$\hfill \square$

\parag{Example} For \ $\delta >1 $ \ a rational number and \ $\lambda \in \mathbb{Q} \cap ]0,1]$ \ the \ $\delta-$dual of \ $\Xi^{(N)}_\lambda $ \ is \ $\Xi^{(N)}_{\delta - \lambda}$. $\hfill \square$

\begin{thm}\label{hol. dual.}
Let \ $X$ \ be a reduced complex space and \ $\mathbb{E}$ \ be an holomorphic family of rank \ $k$ \  frescos parametrized by \ $X$. Assume that for some \ $\delta \in \mathbb{Q}$ \ and for each \ $x \in X$ \ the \ $\hat{A}-$module \ $\mathbb{E}(x)^*\otimes E_{\delta}$ \ is a fresco. Then the family \ $\mathbb{E}^*\otimes E_{\delta}$ \ is holomorphic.
\end{thm}

\parag{Proof}  Remark that we may assume \ $X$ \ connected as the problem is clearly local on \ $X$, so the rank is globally  well defined. We shall make an induction on the rank of the frescos.\\
In rank \ $1$ \ we know that the family may be locally given, assuming that the Bernstein element is \ $a - \lambda_1.b$ \ with \ $\lambda_1 \in \mathbb{Q}^{+*}$,  by the sheaf \ $\mathcal{O}_X[[b]].s^{\lambda_1-1} \subset \Xi_{\lambda,X}^{(0)} $ \ and the (twisted) dual is then defined (locally) by the sheaf \ $\mathcal{O}_X[[b]].s^{\delta - \lambda_1 -1} $ \ where the condition on \ $\delta \in \mathbb{Q}$ \ is that \ $\delta - \lambda_1 > 0 $. So the rank \ $1$ \ is clear.\\
Assume now that for \ $k-1 \geq 1$ \ the theorem is proved and consider  an holomorphic family \ $\mathbb{E}$ \  of rank \ $k$ \  frescos parametrized by \ $X$. We have a exact sequence, thanks to \ref{J-H. hol.},
$$ 0 \to \mathbb{F}_{k-1} \to \mathbb{E} \to \mathbb{E}\big/\mathbb{F}_{k-1} \to 0 $$
where \ $\mathbb{F}_{k-1}(x)$ \ is the family of the \ $(k-1)$ \ terms of the principal J-H. sequence of \ $\mathbb{E}(x)$ \ for each \ $x \in X$. Then we have an exact sequence of \ $\hat{A}_X-$sheaves
$$ 0 \to ( \mathbb{E}\big/\mathbb{F}_{k-1})^*\otimes E_{\delta} \to \mathbb{E}^*\otimes E_{\delta} \to \mathbb{F}_{k-1} ^*\otimes E_{\delta} \to 0  $$
and we know that the families \ $( \mathbb{E}\big/\mathbb{F}_{k-1})^*\otimes E_{\delta} $ \ and \ $\mathbb{F}_{k-1} ^*\otimes E_{\delta}  $ \ are holomorphic, using the induction hypothesis. Then the theorem \ref{suite exacte} gives the holomorphy of \ $ \mathbb{E}^*\otimes E_{\delta}$. $\hfill \blacksquare$\\

Combining these two theorems we obtain the following result.

\begin{thm}\label{deux sur trois}
Let \ $X$ \ be a reduced complex space and consider an exact sequence of \ $\hat{A}_X-$sheaves :
$$ 0 \to \mathbb{F} \to \mathbb{E} \to \mathbb{G} \to 0 .$$
Assume that two of these are homorphic families of  frescos parametrized by \ $X$ \ and  that for each \ $x \in X$ \ we have an exact sequence of frescos given by the fibers at \ $x$
$$0 \to \mathbb{F}(x) \to \mathbb{E}(x) \to \mathbb{G}(x) \to 0 $$
 then the third  family is holomorphic.
\end{thm}

\parag{Proof} The  statement in the case where \ $\mathbb{F}$ \ and \ $\mathbb{E}$ \ are holomorphic is easily reduced to  the rank \ $1$ \ case for \ $\mathbb{F}$ : if \ $\mathbb{F}$ \ has rank \ $r \geq 2$ \ let \ $\mathbb{F}_1$ \ be the first term is the principal J-H. sequence for \ $\mathbb{F}$. Then we have  the exact sequence
$$ 0 \to \mathbb{F}\big/\mathbb{F}_1 \to \mathbb{E}\big/\mathbb{F}_1 \to \mathbb{G} \to 0 $$
where the first two sheaves are holomorphic families. As we know that \ $\mathbb{F}$ \ is a holomorphic family  applying the theorem \ref{J-H. hol.} we obtain the holomorphy of \ $\mathbb{F}\big/\mathbb{F}_1$.  But  $\mathbb{F}_1(x)$ \  is normal in \ $\mathbb{E}(x)$, so  the sequence
$$ 0 \to \mathbb{F}(x)\big/\mathbb{F}_1(x) \to (\mathbb{E}\big/\mathbb{F}_1)(x) \to \mathbb{G}(x) \to 0 $$
is exact for each \ $x$ \ and the induction on the rank of \ $\mathbb{F}$ \  implies the holomorphy of \ $\mathbb{E}\big/\mathbb{F}_1$. Then the theorem \ref{suite exacte} applies to the exact sequence
$$ 0 \to \mathbb{F}_1 \to \mathbb{E} \to \mathbb{E}\big/\mathbb{E}_1$$
and gives the holomorphy of \ $\mathbb{E}$.\\
So we have to prove the first case in the rank \ $1$ \ case for \ $\mathbb{F}$. But, as the statement in local on \ $X$ \ we may assume that \ $\mathbb{E}$ \ is generated by a \ $k-$admissible section \ $\varphi$ \ of some sheaf \ $\Xi_{\Lambda,X}^{(N)}\otimes V$ \ and that  \ $\mathbb{F}$ \ is locally generated by a  section \ $\psi$ \ of \ $ \Xi_{\Lambda,X}^{(N)}\otimes V$ \ satisfying \ $(a - \mu.b).\psi = 0$ \ where \ $-\mu$ \ is the root of the Bernstein polynomial of \ $\mathbb{F}$. If we show that \ $\psi(x)$ \ never vanishes then  we may argue as in the proof of the theorem \ref{J-H. hol.} when we proved that \ $\mathbb{E}\big/\mathbb{F}_1$ \ is holomorphic. But if \ $\psi(x) = 0$ \ for some \ $x \in X$, the image of \ $\mathbb{F}(x)$ \ in \ $\mathbb{E}(x)$ \ is \ $0$ \ contradicting our assumption.\\
The second case (i.e. where \ $\mathbb{E}$ \ and \ $\mathbb{G}$ \ are assumed to be holomorphic)  is consequence of the first one using duality and  the theorem \ref{hol. dual.}. \\
We proved already the third case in the theorem \ref{suite exacte}. $\hfill \blacksquare$

\newpage

\section{Versal families of  frescos.}

\subsection{Preliminaries.}

\parag{Notation} For \ $\lambda_1, \dots, \lambda_k$ \ the fundamental invariants of a  fresco, we shall denote \ $\mathcal{F}(\lambda_1, \dots, \lambda_k)$ \ the set of isomorphism class of  frescos with these fundamental invariants.\\
Define  the rational numbers \ $p_j : = \lambda_{j+1} - \lambda_j + 1$ \ for \ $j \in [1,k-1]$. Recall that we have fixed an order on \ $\mathbb{Q}\big/\mathbb{Z}$ \ and that in the principal J-H. sequence for a fresco the fundamental invariants \ $\lambda_1, \dots, \lambda_k$ \ are listed in such a way that if \ $[\mu_1] < \dots < [\mu_d]$ \ are their distinct classes in \ $\mathbb{Q}\big/\mathbb{Z}$ \ we begin by the principal J-H. sequence of \ $E[\mu_1]$ \ and then continue with the principal J-H. sequence of \ $\big(E\big/E[\mu_1]\big)[\mu_2]$ \ and so on.\\
 Now for \ $j \in [1,k]$ \ let \ $ l_j$ \ be the maximal integer ($0$ \ is possible) such that \ $\lambda_{j}, \dots, \lambda_{j+l_j}$ \ are equal modulo \ $\mathbb{Z}$. Note that the numbers \ $p_j, \dots, p_{j + l_j-1}$ \ are integers (for \ $l_j \geq 1$)  and that \ $p_{j+l_j}$ \ is not an integer. Now, with this notation define
 $$ Y(\lambda_j, \dots, \lambda_k) : = \sum_{h=0}^{k-j-1} \ \C.b^h \quad  \oplus \sum_{q = p_j+ \dots+p_{j+l} \geq k-j}^{ l \in [0,l_j-1]}\ \C.b^{q}\quad \subset \ \C[[b]].$$
 It will be convenient to adopt the  convention \ $Y(\lambda_k) = \C.1 \subset  \C[[b]]$.\\
We shall denote
$$ \mathcal{G}(\lambda_1, \dots, \lambda_k) : = \left\{ (S_1, \dots, S_{k-1}) \ / \  S_j(0) = 1 \quad {\rm and} \quad S_j \in Y(\lambda_j, \dots,\lambda_k)\quad  \forall j \in [1,k-1] \right\}.$$
We have a natural map
$$ [ \quad ] :  \mathcal{G}(\lambda_1, \dots, \lambda_k)  \longrightarrow \mathcal{F}(\lambda_1, \dots, \lambda_k)$$
which send the point  \ $s = (S_1, \dots, S_{k-1}) \in \mathcal{G}(\lambda_1, \dots, \lambda_k)$ \ to the isomorphism class  \ $[s] \in \mathcal{F}(\lambda_1, \dots, \lambda_k)$ \ of the fresco
$$ s : =  \A\big/\A.(a - \lambda_1.b).S_1^{-1} \dots S_{k-1}^{-1}.(a - \lambda_k.b) .$$

\parag{Remark} We have a natural projection
$$ \pi : \mathcal{G}(\lambda_1, \dots, \lambda_k) \to \mathcal{G}(\lambda_2, \dots, \lambda_k) $$
defined in forgetting \ $S_1$. It corresponds to the quotient map \ $E \to E\big/F_1$ \ where \ $F_1$ \ is the first term of the principal J-H. sequence of \ $E$. $\hfill \square$\\

\begin{lemma}\label{suppl.}
 Let \ $T_1, \dots, T_k$ \ be invertible elements in \ $\C[[b]]$ \ and consider the action of 
 $$ P_j : = T_j.(a - \lambda_{j+1}.b).T_{j+1} \dots T_{k-1}.(a - \lambda_k.b).T_k $$
 on \ $E_{\lambda_j}$. We have \ $P_j.E_{\lambda_j} + Y(\lambda_j, \dots, \lambda_k).e_{\lambda_j} = E_{\lambda_j}$, where \ $e_{\lambda_j}$ \ is a standard generateur of \ $E_{\lambda_j}$.
 \end{lemma}
 
 \parag{proof} First remark that \ $T_h$ \ acts bijectively on \ $E_{\lambda_j}$ \ and respect the \ $b-$filtration. Remark also that for any \ $\mu \in \C$ \ the operator \ $(a - \mu.b)$ \ on \ $E_{\lambda_j}$ \ is also closed for the \ $b-$adic filtration. So \ $P_j.E_{\lambda_j}$ \ is closed with finite codimension image and contains \ $b^N.E_{\lambda_j}$ \ for \ $N \gg 1$. Then, to complete the proof, it is enough to show that for each  \ $q \geq k-j, q \not= p_j + \dots + p_{j+l}$ \ for each \ $l \in [0,l_j-1]$, there exists an element in \ $P_j.E_{\lambda_j} $ \ with initial term equal to \ $b^q.e_{\lambda_j}$. But applying \ $P_j$ \ to \ $\gamma.b^{q-k+j}.e_{\lambda_j}$ \ will give  such an element for a good choice of \ $\gamma \in \C^*$, because for each \ $h \geq j+1$, the action of \  $T_h$ \  changes the initial term by a non zero constant and the action of  \ $(a- \lambda_h.b)$ \  changes the initial term by the action of \ $b$ \ and  by a non zero constant, as our choice of \ $q$ \ implies that at each step the initial term is not in the kernel of \ $(a - \lambda_h.b)$. Remark that for \ $h \geq j+l_j$ \ we have \ $\lambda_h - \lambda_j \not\in \mathbb{N}$. $\hfill \blacksquare$
 
 \begin{prop}\label{red. k-adm.}
 Let \ $X$ \ a  connected reduced complex space and let \ $\varphi$ \ a \ $k-$admissible section on \ $X$ \ of the sheaf \ $\Xi^{(N)}_{X,\Lambda} \otimes V$ \ where \ $V$ \ is a finite dimensional complex vector space. Then
 \begin{enumerate}[i)]
 \item For all \ $k-$admissible section \ $\varphi$ \ in \ $\Xi^{(N)}_{X,\Lambda}\otimes V$ \ such that the corresponding frescos have fundamental invariants \ $\lambda_1, \dots, \lambda_k$ \ there exists 
  $$ P_k : = (a - \lambda_1.b).S_1^{-1} \dots S_{k-1}^{-1}.(a - \lambda_k.b).S_k^{-1} $$
  where \ $S_j$ \ are global sections on \ $X$ \ of the sheaf \ $\mathcal{O}_X[[b]]$ \ satisfying \ $S_j(0) \equiv 1$, such that \ $\hat{A}.P_k$ \ is the annihilator ideal of \ $\varphi$ \ in the sheaf \ $\hat{A}_X$.
  \item Moreover,  we may choose locally on \ $X$ \ another local \ $k-$admissible generator \ $\psi$ \  of \ $\hat{A}.\varphi$ \  in order that the annihilator of \ $\psi$ \  is of the form \ $P_k$ \ where \ $S_1, \dots, S_k$ \  satisfy the extra condition 
  $$ S_j \in \mathcal{O}_X\otimes Y(\lambda_j, \dots, \lambda_k) \quad  \forall j \in [1,k-1] \quad  {\rm and \ } \quad S_k \equiv 1.$$
  \end{enumerate}
  \end{prop}
  
  \parag{proof}Let \ $\mathbb{E} : = \hat{A}.\varphi$ \ and let \ $\mathbb{F}_{k-1}$ \ its corank 1 term of its principal J-H. sequence. Thanks to the theorem \ref{J-H. hol.}  and lemma \ref{Hol. rank 1} we have  an isomorphism \ $\mathbb{E}\big/\mathbb{F}_{k-1} \simeq \mathcal{O}_X[[b]].s^{\lambda_k-1} $. Then the image of \ $\varphi$ \ is equal to \ $S_k.s^{\lambda_k-1}$ \ where \ $S_k$ \ is a global section of  the sheaf \ $\mathcal{O}_X[[b]]$ \ which is  locally (and so globally) invertible. Then \ $(a - \lambda_k.b).S_k^{-1}.\varphi$ \ is \ $(k-1)-$admissible and generates \ $\mathbb{F}_{k-1}$. \\
  The \ $(k-1)-$admissibility is consequence of the lemma \ref{easy}. To show that \ $\varphi$ \  generates \ $\mathbb{F}_{k-1}$ \ let \ $\sigma \in \mathbb{F}_{k-1}$ \ and write \ $\sigma = P.S_k^{-1}.\varphi$ \ where \ $P : = \sum_{j=0}^{k-1} P_j.a^j$, and divide \ $P$ \ by \ $a - \lambda_k.b$ :
  $$ P = Q.(a - \lambda_k.b) + U $$
  where \ $U$ \ is a section of the sheaf \ $\mathcal{O}_X[[b]]$ \ and \ $Q$ \ a section of the sheaf \ $\hat{A}_X$ \ which is a polynomial of degree \ $k-1$ \ in \ $a$. Then we have
  $$ \sigma = P.S_k^{-1}.\varphi = Q.(a - \lambda_k.b).S_k^{-1}.\varphi + S_k^{-1}.U.\varphi $$
  and as the image of \ $\sigma$ \ is \ $0$ \ in \ $\mathbb{E}\big/\mathbb{F}_{k-1}$ \ we obtain
  $$ S_k^{-1}.U.s^{\lambda_k-1} \equiv 0$$
  which implies that \ $U \equiv 0$, and this shows that \ $\sigma = Q.(a - \lambda_k.b).S_k^{-1}.\varphi$.\\
  Now by induction on \ $k$, the point i) is proved.
  
  \smallskip
  
  We shall prove the point ii) by induction on \ $k$. As the case \ $k = 1$ \ in consequence of the lemma \ref{Hol. rank 1}, we shall assume the result proved in rank \ $k-1 \geq 1$. Let \ $\mathbb{E} : = \hat{A}_X.\varphi$ \ and \ $\mathbb{F}_1 \simeq \mathbb{O}_X[[b]].s^{\lambda_1-1}$ \ the first term of the principal J-H. of \ $\mathbb{E}$. So applying the induction hypothesis, as the fundamental invariants of \ $\mathbb{E}\big/\mathbb{F}_1$ \ are \ $\lambda_2, \dots, \lambda_k$, we find a local  generator \ $[\varphi]$ \  of \ $\mathbb{E}\big/\mathbb{F}_1$, which is \ $(k-1)-$admissible, induced (locally) by a local \ $k-$admissible section \ $\varphi$  \ of \ $\mathbb{E}$, such that 
  $$ P_{k-1} : = (a - \lambda_2,b).S_2^{-1} \dots S_{k-1}^{-1}.(a - \lambda_k.b) $$
  generates the annihilator of \ $[\varphi]$ \ where \ $S_j, j \in [2,k-1]$ \ satisfy \ $S_j(0) \equiv 1$ \ and \ $ S_j \in Y(\lambda_j, \dots, \lambda_k) $ \ for \ $j \in [2,k]$. Define the (local) section \ $T $ \ of \ $\mathcal{O}_X[[b]]$ \ by the equality \ $P_{k-1}.\varphi = T.s^{\lambda_1-1}$. As it has to be a generator of \ $\mathbb{F}_1$ \ the section \ $ T(0)$ \ is  invertible in \ $\mathcal{O}_X$. Now the lemma \ref{suppl.}, as we are working on a Stein open set, allows to find a section \ $S_1$ \ of the sheaf \ $\mathcal{O}_X[[b]]$ \ such that
  $$ S_j(0) \equiv 1, \quad S_j \in \mathcal{O}_X\otimes Y(\lambda_1, \dots, \lambda_k) \quad {\rm and} \quad T(0).S_1 - T \in P_{k-1}.\mathbb{F}_1.$$
  Note that, as \ $P_{k-1}.\mathbb{F}_1$ \ contains \ $b^N.\mathbb{F}_1$ \ for \ $N \gg 1$, it is easy to show that \ $P_{k-1}.\mathbb{F}_1\big/b^N.\mathbb{F}_1$ \ is a \ $\mathcal{O}_X-$coherent subsheaf of a locally free \ $\mathcal{O}_X-$sheaf, and then to use the surjectivity of the \ $\mathcal{O}_X-$linear  maps
  $$ P_{k-1}.\mathbb{F}_1\big/b^N.\mathbb{F}_1 \oplus \big(Y(\lambda_1, \dots, \lambda_k)\otimes \mathcal{O}_X\big) \longrightarrow  \mathbb{F}_1\big/b^N.\mathbb{F}_1 $$
  and 
  $$ \mathbb{F}_1\big/b^{N}.\mathbb{F}_1 \overset{P_{k-1}}{\longrightarrow} P_{k-1}.\mathbb{F}_1\big/b^N.\mathbb{F}_1$$
  will allow to use H. Cartan theorem B  on a Stein open set.
  So, using again the lemma \ref{suppl.} and the fact that we work on an open Stein set of \ $X$, we may find with the argument above a section \ $\chi$ \ of \ $\mathbb{F}_1$ \ such that 
   $$P_{k-1}.\chi = (T(0).S_1 - T).s^{\lambda_1-1} .$$
   Then \ $\psi : = S_k^{-1}.\varphi - \chi $ \ is again \ $k-$admissible and generates \ $\mathbb{E}$. Moreover it satisfies \ $P_{k-1}.\psi = T(0).S_1.s^{\lambda_1-1}$ \ and then
   $$ (a - \lambda_1.b).S_1^{-1}.P_{k-1}.\psi \equiv 0 .$$
   Then it is easy to see that \ $P_k : = (a - \lambda_1.b),S_1^{-1}$ \ generates the annihilator of \ $\psi$. $\hfill \blacksquare$\\
   
   \subsection{The existence theorem.}

\begin{thm}\label{standard versal family}
The tautological family \ $\mathbb{E}^{taut}$ \ parametrized by the affine space \ $G_k = \mathcal{G}(\lambda_1, \dots, \lambda_k)$ \ whose value at \ $S_1, \dots,  S_{k-1}$ \ is given by the sheaf
  $$\hat{A}_{G_k}\big/\hat{A}_{G_k}.(a - \lambda_1.b).S_1^{-1} \dots S_{k-1}^{-1}.(a -\lambda_k.b)$$
  is holomorphic\footnote{in fact polynomial !}. Moreover this family is locally versal at each point of \ $G_k$.
\end{thm}

\parag{proof} To prove the holomorphy of \ $\mathbb{E}^{taut}$ \ we shall construct by induction on \ $k \geq 1$ \ a \ $k-$admissible global section \ $\varphi_k$ \  of the sheaf
$$ \Xi^{(k-1)}_{G_k,\Lambda} \otimes \C^k $$
with annihilator generates in \ $\hat{A}_{G_k}$ \ by the section \
$$P_k : = (a - \lambda_1.b).S_1^{-1} \dots S_{k-1}^{-1}.(a -\lambda_k.b) .$$
Then we shall show that the family \ $\mathbb{E}[G_k] : = \hat{A}_{G_k}.\varphi_k$ \ is locally versal near each point.\\
The rank \ $1$ \ case is obvious as \ $G_1$ \ is one point. Then \ $\mathbb{E}(G_1)$ \ is the rank \ $1$ \ fresco \ $E_{\lambda_1}$. Note that in this case, this family is locally  universal thanks to the lemma  \ref{Hol. rank 1}.\\
Now assume \ $k \geq 2$ \ and that we have already constructed \ $\varphi_{k-1}$ \ with annihilator generated by \ $P_{k-1}$ \ and such that we obtain a locally versal family. As we may increase \ $\Lambda$ \ without problem, we shall assume that \ $\Lambda$ \ contains all \ $[\lambda_j], j \in [1,k]$.  We have a natural projection \ $ \pi : G_k \to G_{k-1} : = \mathcal{G}(\lambda_2, \dots, \lambda_k)$ \ and we may lift \ $\varphi_{k-1}$ \ successively to a section of the sheaves \ $\Xi^{(k-2)}_{G_k,\Lambda}\otimes \C^{k-1} \subset \Xi^{(k-1)}_{G_k,\Lambda}\otimes \C^k$ \ where we define \ $\C^k : = \C^{k-1}\oplus \C.e_k$. Define now \ $\varphi_k$ \ as the solution of the equation
$$ (a - \lambda_k.b).\varphi_k = S_{k-1}.[\pi^*(\varphi_{k-1}) \oplus 0 ] $$
satisfying the following condition :
 \begin{itemize}
 \item The component of \ $ \varphi_k$ \   on \ $ \Xi^{(k-1)}_{G_k,\Lambda}\otimes \C.e_k$ \    is equal to \ $s^{\lambda_k-1}\otimes e_k$. 
 \end{itemize}
  It is easy to see that such a solution exists and is unique with this extra condition.\\
 To prove that \ $\varphi_k$ \ is admissible, as \ $\varphi_k, a.\varphi_k, \dots, a^{k-1}.\varphi_k$ \ clearly generates \ $\hat{A}_{G_k}.\varphi_k$ \ as a \ $\mathcal{O}_{G_k}[[b]]-$modules thanks to the \ $k-$admissibility of \ $\varphi_{k-1}$, it is enough to prove that a relation like
 $$ \sum_{j=0}^{k-1} \ T_j.a^j.\varphi_k \equiv 0 $$
 where \ $T_j, j \in [0,k-1]$ \ are local sections of the sheaf \ $\mathcal{O}_{G_k}[[b]]$ \ implies that \ $T_j \equiv 0$ \ for each \ $j$. But a simple computation shows that it is equivalent to prove that a relation like
 $$ \tilde{T}_0.\varphi_k + \sum_{j=1}^{k-1} \ \tilde{T}_j.a^{j-1}.(a - \lambda_k.b).\varphi_k \equiv 0 $$
 also implies the vanishing of the local sections \ $\tilde{T}_j$ \ for \ $j \in [0,k-1]$. Looking at the component on \ $\Xi^{(k-1)}_{G_k,\Lambda}\otimes \C.e_k$ \ we see that  \ $\tilde{T}_0 \equiv 0$ \ and the \ $(k-1)-$admissibilty of \ $\varphi_{k-1}$ \ allows to conclude.\\
 As it is clear that \ $P_k$ \ generates the annihilator of \ $\varphi_k$, to conclude the proof we have only to show the local versality of the holomorphic family defined by \ $\varphi_k$.\\
 So let \ $\mathbb{E}$ \ be a holomorphic family of  frescos parametrized by a connected complex space \ $X$ \ with fundamental invariants \ $\lambda_1, \dots, \lambda_k$, and let \ $x_0$ \ be a point in \ $X$.  Up to replace \ $X$ \ by a small Stein connected open set \ $X_0$ \ containing \ $x_0$ \ we may assume that \ $\mathbb{E}\big|_{X_0}$ \ is generated by a \ $k-$admissible section \ $\varphi$ \ of the sheaf \ $\Xi^{(N)}_{X_0,\Lambda}\otimes V$ \ where \ $V$ \ is a finite dimensional complex vector space. The proposition \ref{red. k-adm.} allows, up to skrink \ $X_0$, to replace \ $\varphi$ \ by a generating section \ $\psi$ \ with annihilator \ $P_k : = (a - \lambda_1.b).S_1^{-1}\dots S_{k-1}^{-1}.(a - \lambda_k)$ \ with the conditions \ $S_j(0) \equiv 1$ \ and \ $S_j \in \mathcal{O}_X\otimes Y(\lambda_j, \dots, \lambda_k) $. Then this defines a  holomorphic map  \ $f : X_0 \to \mathcal{G}(\lambda_1, \dots, \lambda_k)$, and it is clear that \ $f^*(\mathbb{E}^{taut})$ \ is then isomorphic on \ $X_0$ \ with the sheaf \ $\hat{A}.\psi = \hat{A}.\varphi \simeq \mathbb{E}\big|_{X_0}$. $\hfill \blacksquare$\\

We study in detail the classification of rank \ $3$ \  frescos in the Appendix. This gives many examples of locally universal holomorphic families of such frescos (see also [B.10] for more such examples.)

\section{The theorem of the semi-simple part.}

\subsection{The statement.}

The aim of this section is to prove the following result

\begin{thm}\label{hol. ss part}
Let \ $\mathbb{E}$ \ be a holomorphic family of \ $[\lambda]-$primitive frescos parametrized by a reduced complex space \ $X$. Then there exists a dense Zariski open set \ $X'$ \ in \ $X$ \ on which the family of semi-simple frescos \ $(S_1(\mathbb{E}(x)), x \in X')$ \ is holomorphic.
\end{thm}

Note that we may not avoid to restrict our family to a dense Zariski open set in general  because of the example of the rank \ $2$ \ family
$$ \mathbb{E}(z) : = \hat{A}\big/\hat{A}.P_z $$
where \ $P_z : = (a - \lambda_1.b).(1 + z.b^p)^{-1}.(a - (\lambda_1+p-1).b) $ \ where \ $z \in \C, \lambda_1>1$ \ is in \ $\mathbb{Q}$ \ and \ $p \in \mathbb{N}^*$. This holomorphic family satisfies \ $S_1(\mathbb{E}(z)) \simeq E_{\lambda_1}$ \ for \ $z \in \C^*$ \ because \ $\mathbb{E}(z)$ \ is a theme in these cases, and \ $\mathbb{E}(0)$ \ is semi-simple.\\

\parag{Remark}
Some more results will be obtained in the proof of this theorem
\begin{enumerate}[i)]
\item There is natural \ $\hat{A}_X-$subsheaf \ $\mathcal{S}_1(\mathbb{E})\subset \mathbb{E}$ \ on \ $X$ \ which is a  finetely generated \ $\mathcal{O}_X[[b]]-$module such that on a dense Zariski open set \ $X'$ \ in \ $X$ \ it induces the holomorphic family of the semi-simple parts.
\item The subsheaf \ $\mathcal{S}_1(\mathbb{E})$ \ is normal and the quotient \ $\mathcal{S}_1(\mathbb{E})\big/b.\mathcal{S}_1(\mathbb{E})$ \ is a \ $\mathcal{O}_X-$coherent subsheaf of the locally free \ $\mathcal{O}_X-$module \ $\mathbb{E}\big/b.\mathbb{E}$.
\item It is probably not difficult to extend the above theorem to a general holomorphic family of frescos using the same lines that in our proof.
\end{enumerate}

\subsection{Preparations.}

The proof will need some preparations. First we have to recall some facts which are proved in  [B. 12]. 
Let \ $E$ \ be a fresco.
\begin{enumerate}[i)]
\item A necessary and sufficient condition for a fresco \ $E$ \ to be semi-simple is that \ $E$ \ may be embedded in \ $\Xi^{(0)}_{\Lambda} \otimes V$ \ for some finite set \ $\Lambda \subset ]0,1] \cap \mathbb{Q}$ \ and some finite dimensional complex vector space \ $V$. We may choose \ $V$ \ in order that \ $ \dim_{\C}(V) \leq rk(E)$.
\item For any \ $E$ \ there is a maximal semi-simple submodule \ $S_1(E)$ \ (see section 2.3). This submodule is normal, so it is a  semi-simple fresco. Its rank will be noted \ $\delta(E)$. If \ $E$ \ is embedded in \ $\Xi^{(N)}_{\Lambda} \otimes V$ \ then we have
$$ S_1(E) \simeq (\Xi^{(0)}_{\Lambda} \otimes V) \cap E .$$
\end{enumerate}

The following proposition is the key of the description of semi-simple part of a \ $[\lambda]-$primitive fresco which will be used in presence of holomorphic parameters.

\begin{prop}\label{tool1}
Let \ $E$ \ be a \ $[\lambda]-$primitive fresco with rank \ $k$ \ and  fundamental invariants \ $\lambda_1, \dots, \lambda_k$. For \ $\mu \in [\lambda]$ \ define \ $K_{\mu}$ \ as the kernel of \ $(a - \mu.b)$ \ acting on \ $E$. Then the following properties hold:
\begin{enumerate}[i)]
\item The vector space \ $K_{\mu}\big/b.K_{\mu-1} $ \ has dimension \ $1$ \ if  \ $-(\mu - \delta(E)+1)$ \ is a root of the Bernstein polynomial of \ $S_1(E)$, the semi-simple part of \ $E$.
\item The vector space \ $K_{\mu}\big/b.K_{\mu-1} $ \ is \ $\{0\}$ \ in all other cases.
\item For \ $\mu \geq \lambda_k + k - 1$ \ the dimension of \ $K_{\mu}$ \ is  equal to \ $\delta(E)$ \ the rank of \ $S_1(E)$.
\item For \ $\mu \geq \lambda_k+k-1$ \ and \ $q \geq \mu - \lambda_1 $ \ we have the equality
$$ S_1(E) = \left(b^{-q}.\C[[b]].K_{\mu}\right) \cap \  E $$
in \ $E[b^{-1}] : = E \otimes_{\C[[b]]} \C[[b]][b^{-1}]$.
\end{enumerate}
\end{prop}

The proof will use some lemma.

\begin{lemma}\label{une clef}
Let \ $E$ \ be a \ $[\lambda]-$primitive fresco with rank \ $k$ \ and  fundamental invariants \ $\lambda_1, \dots, \lambda_k$. Define the non negative integer \ $p$ \  by the formula
$$ \lambda + p = \lambda_1 - k + 1 .$$
If we have an embedding \ $E \hookrightarrow \Xi^{(N)}_{\lambda} \otimes V$ \ then, in fact, \ $E$ \ is contained in \ $b^p.\Xi^{(N)}_{\lambda} \otimes V$.
\end{lemma}

\parag{proof} First recall that in our situation we have \ $\lambda_1 > k-1$ \ and \ $\lambda$ \ is in \ $]0,1] \cap \mathbb{Q}$. So the integer \ $p$ \ we defined above is non negative.\\
We shall prove the lemma by induction on the rank \ $k$ \ of \ $E$. Note first that for \ $k = 1$ \ the result is obvious. Assume that the result is proved for \ $k-1$ \ for \ $k \geq 2$ \ and consider a rank \ $k$ \ fresco \ $E$. Note \ $F_{k-1}$ \ the rank \ $k-1$ \ submodule of the principal J-H. of \ $E$. Then the fundamental invariants of \ $F_{k-1}$ \ are \ $\lambda_1, \dots, \lambda_{k-1}$. Consider an embedding of \ $E$ \ in \ $ \Xi^{(N)}_{\lambda} \otimes V$ \ given by the image \ $\varphi $ \ of a  generator of \ $E$. Up to change \ $\varphi$ \ by an invertible element in \ $\C[[b]]$ \ we may assume that \ $ \psi : = (a - \lambda_k.b).\varphi$ \ generates the image of \ $F_{k-1}$. The inductive hypothesis implies that \ $\psi$ \ lies in \ $b^{p+1}.\Xi^{(N)}_{\lambda} \otimes V$. Then \ $\varphi$ \ and also \ $E$ \ is contained in \ $b^p.\Xi^{(N)}_{\lambda} \otimes V$. $\hfill \blacksquare$

\begin{cor}\label{cas absolu}
Let \ $E$ \ be a \ $[\lambda]-$primitive fresco with rank \ $k$ \ and fundamental invariants \ $\lambda_1, \dots, \lambda_k$. Let \ $\mu \in [\lambda]$ \ and \ $q \in \mathbb{N}$ \ such that \ $\mu - \lambda_1 \leq q - k + 1$ \ we have
$$ (b^{-q-1}.\C[[b]].K_{\mu}) \cap E = (b^{-q}.\C[[b]].K_{\mu}) \cap E  $$
in \ $E \otimes_{\C[[b]]} \C[[b]][b^{-1}]$ \ where we put as above \ $ K_{\mu} : = Ker[ (a - \mu.b) : E \to E ]$.
\end{cor}

\parag{proof} We want to prove that if \ $x$ \ is in $ \C[[b]].K_{\mu} \cap b^{q+1}.E$ \ then \ $x$ \ is in \ $b.\C[[b]].K_{\mu}  $. Choose an embedding\footnote{such an embedding always exists for a \ $[\lambda]-$primitive fresco thanks to the theorem \ref{embed.}.}  of \ $E$ \ in \ $\Xi^{(N)}_{\lambda} \otimes V$ \ and define \ $p \in \mathbb{N}$ \  as in the lemma \ref{une clef}. Then \ $x$ \ is in \ $b^{q+1}.E \subset b^{p+q+1}.\Xi^{(N)}_{\lambda} \otimes V$. If we write \ $x = y + z$ \ where \ $y \in K_{\mu}$ \ and \ $z \in b.\C[[b]].K_{\mu}$ \ we will have \ $y = \rho.s^{\mu-1}$ \ for some \ $\rho \in \C$. So, if \ $\rho \not= 0$, we must  have
$$ \mu - 1 \geq \lambda + p + q = \lambda_1 + q - k +1$$
contradicting our assumption;  so\ $\rho = 0$ \ and  \ $x = z$ \ is in \  
$b.\C[[b]].K_{\mu}$. $\hfill \blacksquare$\\

 \begin{lemma}\label{rang S1}
 Let \ $E$ \ be a (a,b)-module and let \ $x_1, \dots, x_p$ \ be \ $\C-$linearly independant elements in \ $E$ \ satisfying \ $(a - \mu.b).x_j = 0, \quad \forall j \in [1,p]$. Then they are independant over the ring \ $\C[[b]]$.
 \end{lemma}
 
 \parag{proof} Consider \ $S_1, \dots, S_p$ \ in \ $\C[[b]]$ \ such that at least one of them is not zero and satisfying
 $$ \sum_{j=1}^p \ S_j(b).x_j = 0 \quad {\rm in }\ \quad E.$$
 Assume, as \ $b$ \ is injective on \ $E$, that not all \ $S_j(0)$ \ vanishes.
 Applying \ $(a - \mu.b)$ \ to the previous equality gives, using again the injectivity of \ $b$,
 $$ \sum_{j=1}^p \ S'_j(b).x_j = 0 .$$
 So we have the same relation for any derivative of the \ $S_j$. This gives, using \ $b-$completion to prove Taylor formula in \ $\C[[b]]$,
 $$ \sum_{j=1}^p \ \left(\sum_{\nu = 0}^{\infty} \ \frac{1}{\nu !}.S^{(\nu)}_j(b).(-b)^{\nu}\right).x_j = 0 $$
 which implies \ $\sum_{j=1}^p \ S_j(0).x_j = 0$. And this contradicts our assumption. $\hfill \blacksquare$
 
 \parag{Consequence} We obtain from the previous lemma, observing that a rank \ $1$ \ sub-module is always contained in \ $S_1(E)$,  that \ $\dim_{\C} K_{\mu} \leq rank(S_1(E))$, with equality for \ $\mu$ \ large enough, because in a  semi-simple module \ $F$ \  the rank \ $1$ sub-modules of \ $F$ \  generates \ $F^{\flat}$, the maximal simple pole  sub-module in \ $F$,  which has same rank than \ $F$ \ (see [B.12] section 2). $\hfill \square$

\parag{proof of proposition \ref{tool1}} Consider an element \ $x$ \ in \ $E \setminus b.E$ \ such that \ $a.x = \mu.b.x$. Then \ $\C[[b]].x$ \ is a normal  rank \ $1$ \ sub-(a,b)-module of \ $E$, so it is contained in \ $S_1(E)$ \ and so \ $\mu - \delta(E)+1$ \ is the opposite of a root of the Bernstein polynomial of \ $S_1(E)$. \\
This prove that, in the case where \ $\mu - \delta(E)+1$ \ is  not the opposite of a root of the Bernstein polynomial of \ $S_1(E)$, an element \ $x \in E$ \ such \ $a.x = \mu.b.x$ \ is in \ $b.E$. But if \ $x = b.y$ \ we have, as \ $b$ \ is injective, \ $y \in K_{\mu-1}$ \ and so \ $K_{\mu} = b.K_{\mu-1}$. This proves i) and ii).\\
Note that in the case i) \ $\mu - k+1$ \ is also the opposite of a root of the Bernstein polynomial of \ $E$, so there exists \ $j \in [1,k]$ \ such that \ $\mu - \delta(E)+1 = \lambda_j+j - k$. Then \ $\mu = \lambda_j + j -1 + \delta(E) - k$ \ for some \ $j \in [1,k]$. So we have \ $\mu \leq \lambda_k+k-1$. So for \ $\mu \geq \lambda_k+k-1$ \ we have, using the consequence of the lemma \ref{rang S1}, \ $\dim_{\C} K_{\mu} = \delta(E)$, and point iii) is proved.\\
To prove point \ $iv)$ \ remark that \ $b^{-q}.\C[[b]].K_{\mu}$ \ is a semi-simple (regular)  (a,b)-module, as direct sum of rank \ $1$ \ regular (a,b)-modules. As it has rank  equal to \ $\dim_{\C} K_{\mu} = \delta(E) $, it is enough  to prove that \ $\left(b^{-q}.\C[[b]].K_{\mu}\right) \cap \  E $ \ is normal in \ $E$ : it is contained  in \ $S_1(E)$ \ and with the same rank, so  normality will give the equality with \ $S_1(E)$. \\
To prove normality  let \ $x \in E$ \ such that  \ $b.x \in \left(b^{-q}.\C[[b]].K_{\mu}\right) \cap \  E $. Then  \ $x \in \left(b^{-q-1}.\C[[b]].K_{\mu}\right) \cap E$ \   and the equality proved in the  corollary \ref{cas absolu}
 $$\left(b^{-q}.\C[[b]].K_{\mu}\right) \cap \  E  = \left(b^{-q-1}.\C[[b]].K_{\mu}\right) \cap \  E $$
 for  \ $\mu \geq \lambda_k+k-1$ \ and \ $q \geq \mu - \lambda_1+ k - 1 $ \ allows to conclude.$\hfill \blacksquare$

 \subsection{Coherence of \ $\mathcal{K}_{\mu}$.}
 
 We begin by a very simple remark.
 
\parag{Remark} Let \ $\mathbb{E}$ \ be a locally free finite type \ $\mathcal{O}_X[[b]]-$module. If \ $U$ \ is an open set in \ $X$ \ and if \ $\sigma \in \Gamma(U, \mathbb{E})$ \ satisfies \ $\sigma(x) \in b^q.\mathbb{E}(x), \forall x \in U$ \ (resp. \ $\sigma(x) = 0, \forall x \in U$) then \ $ \sigma$ \ is in \ $\Gamma(U, b^p.\mathbb{E})$, (resp. \ $\sigma = 0$).\\
 
Let us give now a  version  with parameter of the corollary \ref{cas absolu} and of a part of  the proposition \ref{tool1}.

\begin{cor}\label{cas relatif}
Let \ $\mathbb{E}$ \ an holomorphic family of \ $[\lambda]-$primitive frescos parametrized by a connected reduced complex space \ $X$. Note \ $\lambda_1, \dots, \lambda_k$ \ the fundamental invariants of the \ $\mathbb{E}(x)$ \ and define
$$ \mathcal{K}_{\mu} : Ker[(a - \mu.b) : \mathbb{E} \to \mathbb{E} ] .$$
Then for  \ $q \geq \mu - \lambda_1 + k - 1$ \ we have
$$ (b^{-q-1}.\mathcal{O}_X[[b]].\mathcal{K}_{\mu})\, \cap \, \mathbb{E} = (b^{-q}.\mathcal{O}_X[[b]].\mathcal{K}_{\mu})\, \cap \,\mathbb{E}  $$
in \ $b^{-q-1}.\mathbb{E}$.\\
Moreover for \ $\mu \geq \lambda_k + k -1$ \ we have \ $ b.\mathcal{K}_{\mu} = \mathcal{K}_{\mu+1}$.
\end{cor}

\parag{proof} Let \ $\alpha$ \ a germ of section of the sheaf \ $(\mathcal{O}_X[[b]].\mathcal{K}_{\mu})\, \cap \, b^{q+1}.\mathbb{E}$ \ and write \ $\alpha = \beta + \gamma$ \ where \ $\beta$ \ is a germ of section of \ $\mathcal{K}_{\mu}$ \ and \ $\gamma$ \ a germ of section of \ $b.\mathcal{O}_X[[b]].\mathcal{K}_{\mu}$. Now embed (locally) \ $\mathbb{E}$ \ in \ $\Xi^{(N)}_{\lambda,X} \otimes V $. Then the lemma \ref{une clef} implies that \ $\mathbb{E}$ \ is contained in \ $b^p.\Xi^{(N)}_{\lambda,X} \otimes V $. Then write \ $\beta = \rho.s^{\mu-1}$ \ where \ $\rho$ \ is a germ of section of \ $\mathcal{O}_X$. If \ $\rho \not \equiv 0$,  we have
$$ \mu - 1 \geq \lambda - 1 + p + q + 1$$
which contradicts our assumption. So \ $\alpha = \gamma$ \ is in \ $b.\mathcal{O}_X[[b]].\mathcal{K}_{\mu}$ \ and we conclude using the injectivity of \ $b$.\\
To prove the second statement let \ $\sigma$ \ be a germ of section of \ $\mathcal{K}_{\mu+1}$ \ where we assume \ $\mu \geq \lambda_k + k - 1$. For each \ $x$ \ we have \ $\sigma(x) \in \mathbb{E}(x) $ \ which is in \ $K_{\mu+1}(x) = b.K_{\mu}(x)$ \ where \ we note \ $K_{\mu}(x)$ \ the kernel of \ $a - \mu.b$ \ acting on \ $\mathbb{E}(x)$ \ and we apply the proposition \ref{tool1}. Using the remark which begins the section 4.3, we conclude that \ $\sigma$ \ is in \ $b.\mathbb{E}$ \ and so in \ $b.\mathcal{K}_{\mu}$. $\hfill \blacksquare$\\

\begin{lemma}\label{basis of the saturation}
Let \ $E$ \ be a rank \ $k$ \ fresco with generator \ $e$. Then its saturation \ $E^{\sharp}$ \ by \ $b^{-1}.a$ \ is a free rank \ $k$ \ $\C[[b]]-$module with basis \ $\{b^{-j}.a^j.e, \, j \in [0,k-1]\}$. Moreover, the images in \ $E^{\sharp}\big/E$ \ of \ $b^{-j}.a^j.e, j\in [1,k-1]$ \ form a \ $\C-$basis of \ $E^{\sharp}\big/E$. Assuming that \ $(a - \lambda_1.b).S_1^{-1} \dots S_{k-1}^{-1}.(a - \lambda_k.b) $ \ generates the annihilator of \ $e$ \ and defining \ $e_1, \dots, e_k$ \ by the relations \ $e_k : = e$ \ and \ $(a - \lambda_j.b).e_j = S_{j-1}.e_{j-1}$ \ for \ $j \in [1,k]$ \ (with the convention \ $e_0 : = 0$), then the set \ $\{b^{-k+j}.e_j, j \in [1,k]\}$ \ is a \ $\C[[b]]-$basis of \ $E^{\sharp}$ \ and the images of the \ $b^{-k+j}.e_j, j \in [1,k-1]$ \ give a \ $\C-$basis of \ $E^{\sharp}\big/E$.
\end{lemma}

\parag{proof} Let \ $F_{k-1}$ \ be the sub-$\A-$module generated by \ $e_{k-1}$. Then it is easy to deduce from the relation \ $(b^{-1}.a - \lambda_k).e_k = S_{k-1}.b^{-1}.e_{k-1} $ \ that the equality
\begin{equation*}
E^{\sharp} = E + b^{-1}.F_{k-1}^{\sharp} \tag{@} 
\end{equation*}
holds. Using induction on the rank we obtain that  \ $\{ b^{-j}.a^j.S_{k-1}.e_{k-1}, j \in [1,k-1]\ \}$ \ is a \ $\C[[b]]-$basis for \ $F^{\sharp}$.\\
 This implies that \ $e_k$ \ and  \ $\{b^{-j-1}.a^j.S_{k-1}.e_{k-1}, j\in [1,k-1]\}$ \ generate \ $E^{\sharp}$ \ as  a \ $\C[[b]]-$module. Then the equality
$$ b^{-j-1}.a^j.S_{k-1}.e_{k-1} = b^{-j-1}.a^j.(a - \lambda_k.b).e_k =  b^{-j-1}.a^{j+1}.e_k - \lambda_k.b^{-j-1}.a^j.b.e_k $$
implies our first statement using the identity in \ $\A$
$$ a^j.b = b.a^j + \sum_{h=0}^{j-1} T_{h,j}.b^{j-h+1}.a^h $$
as we know that \ $E$ \ is free of rank \ $k$ \ on \ $\C[[b]]$. Now the relation \ $(@)$ \ shows that \ $E^{\sharp}\big/E \simeq b^{-1}.F_{k-1}^{\sharp}\big/F_{k-1} $, as we have \ $b^{-1}.F^{\sharp} \cap E = F$ \ because  \ $F_{k-1}$ \ is normal in \ $E$. Using again induction on the rank, we may assume that the set \ $\{b^{j-(k-1)}.e_j, j \in [1,k-2]\}$ \ is a \ $\C-$basis of \ $F_{k-1}^{\sharp}\big/F_{k-1}$. Then, thanks to the isomorphism above, it is enough to add \ $b^{-1}.e_{k-1}$ \ to the \ $b^{j-k}.e_j, j \in [1,k-2]$ \ to generate the vector space \ $E^{\sharp}\big/E$. And it is easy to show that \ $b^{-1}.e_{k-1}$ \ is not a linear combination of the  \ $b^{j-k}.e_j, j \in [1,k-2]$.\\
We leave as an exercice for the reader the proof of the other statements. $\hfill \blacksquare$

\begin{cor}\label{sature relatif}
Let \ $\varphi \in \Xi^{(N)}_{\lambda,X} \otimes V$ \ be a k-admissible section and put \ $\mathbb{E} : = \hat{A}_X.\varphi$. Define \ $\mathbb{E}^{\sharp}$ \ as the \ $\mathcal{O}_X[[b]]-$module generated by the sections \ $b^{-j}.a^j.\varphi,\, j \in [0,k-1]$. Then \ $\mathbb{E}^{\sharp}$ \ is a free rank \ $k$ \ $\mathcal{O}_X[[b]]-$module, stable by \ $b^{-1}.a$ \ and for each \ $x \in X$ \ we have
$$ \mathbb{E}^{\sharp}(x) \simeq \mathbb{E}(x)^{\sharp}.$$
\end{cor}

\parag{proof} It is clear from the previous lemma  that the set \ $\{b^{-j}.a^j.\varphi, \, j \in [0,k-1]\}$ \ is a \ $\mathcal{O}_X[[b]]-$basis of \ $\mathbb{E}^{\sharp}$. Also the quotient \ $\mathbb{E}^{\sharp}\big/\mathbb{E}$ \ is a free \ $\mathcal{O}_X-$module with basis \ $\{b^{-j}.a^j.\varphi, j \in [1,k-1]\}$. So the exact sequence of \ $\mathcal{O}_X-$modules
$$ 0 \to \mathbb{E} \to \mathbb{E}^{\sharp} \to \mathbb{E}^{\sharp}\big/\mathbb{E} \to 0 $$
splits. Then for each \ $x$ \ the map \ $\mathbb{E}(x) \to \mathbb{E}^{\sharp}(x)$ \ is injective. As \ $\mathbb{E}^{\sharp}$ \ has a simple pole the lemma  \ref{basis of the saturation} allows to see that the map induced between \ $\mathbb{E}(x)^{\sharp}$ \ and \ $\mathbb{E}^{\sharp}(x)$ \ is bijective. $\hfill \blacksquare$

\begin{lemma}\label{monomes}
Let \ $\varphi \in \Xi^{(N)}_{\lambda,X} \otimes V$ \ be a k-admissible section and let \ $\lambda_1, \dots, \lambda_k$ \ be the fundamental invariants of the corresponding holomorphic family of frescos. Let \ $\mu \in [\lambda]$ \ and \ $N \in \mathbb{N}$ \ such that \ $ N \geq \mu - \lambda_1 + k $. Then if \ $ \sigma$ \ is a section of the sheaf \ $\mathbb{E}^{\sharp}$ \ such that \ $(a - \mu.b).\sigma \in b^{N+1}.\mathbb{E}^{\sharp}$ \ then there exists \ $\tau \in \mathbb{E}^{\sharp}$ \ such that
\begin{enumerate}[i)]
\item  \ $ (a - \mu.b).\tau = 0 $ \ ;
\item  \ $ \sigma - \tau $ \ is a section of the subsheaf \ $b^N.\mathbb{E}^{\sharp}$.
\end{enumerate}
\end{lemma}

\parag{proof} Write \ $(a - \mu.b).\sigma = b^{N+1}.z_{N+1}$ \ and look for a \ $x_N \in \mathbb{E}^{\sharp}$ \ such that \ $(a - \mu.b).(\sigma + b^N.x_N) $ \ is in \ $b^{N+2}.\mathbb{E}^{\sharp}$. The equation
$$ (a - \mu.b).(\sigma + b^N.x_N) = b^{N+2}.y $$
is then equivalent, as \ $b$ \ is injective with
$$ (b^{-1}.a -  (\mu-N)).x_N = - z_{N+1} + b.y .$$
But on a Stein open set \ $U$ \  in  \ $X$ \ the map \ $(b^{-1}.a - (\mu-N))$ \ is bijective\footnote{Recall that \ $\lambda_1+1-k$ \ is the smallest eigenvalue of \ $b^{-1}.a$ \ acting on \ $E^{\sharp}\big/b.E^{\sharp}$ \ by definition of the Bernstein polynomial.} on \ $\Gamma(U, \mathbb{E}^{\sharp}\big/b.\mathbb{E}^{\sharp})$ \ and we have also the exactness\footnote{recall that \ $\mathbb{E}^{\sharp}$ \ is a free rank \ $k$ \ $\mathcal{O}_X[[b]]-$module, and we have on the Stein open set \ $U$ \ $H^1(U, \mathcal{O}_X[[b]]) = 0.$} of the sequence
$$ 0 \to \Gamma(U, b.\mathbb{E}^{\sharp}) \to \Gamma(U, \mathbb{E}^{\sharp}) \to \Gamma(U, \mathbb{E}^{\sharp}\big/b.\mathbb{E}^{\sharp}) \to 0 $$
which allows to find a solution \ $x_N \in \Gamma(U, \mathbb{E}^{\sharp}) $ \ for our equation. Then an easy induction on \ $j \geq 0$ \ allows to define a sequence \ $(x_{N+j})_{j\geq 0}$ \ of elements in \ $\Gamma(U, \mathbb{E}^{\sharp})$ \ such that, if we define \ $\tau \in  \Gamma(U, \mathbb{E}^{\sharp})$ \ as 
$$ \tau : = \sigma + \sum_{j=0}^{\infty} \ b^{N+j}.x_{N+j} ,$$
it satisfies the conditions \ $i)$ \ and \ $ii)$. $\hfill \blacksquare$\\

\begin{cor}\label{monomes E}
In the situation of the previous lemma, assume that \ $\sigma$ \ is now a section of \ $\mathbb{E}$ \ such that \ $(a - \mu.b).\sigma \in b^{N+1}.\mathbb{E}$. Then there exists a section \ $\tau \in \mathbb{E}$ \ such that
\begin{enumerate}[i)]
\item  \ $ (a - (\mu+k).b).\tau = 0 $ \ ;
\item  \ $ b^k.\sigma - \tau $ \ is a section of the subsheaf \ $b^{N}.\mathbb{E}$
\end{enumerate}
\end{cor}

\parag{proof} As a consequence of the lemma \ref{basis of the saturation} the inclusion \ $b^k.\mathbb{E}^{\sharp} \subset \mathbb{E}$ \ holds. The corollary is now an easy consequence of the lemma above. $\hfill \blacksquare$\\

\begin{lemma}\label{ordre suff.}
Let \ $E$ \ be a \ $[\lambda]-$primitive fresco with fundamental invariants \ $\lambda_1, \cdots, \lambda_k$, and let \ $\mu \in [\lambda]$. Then for \ $N \geq \mu - \lambda_1 + k$ \ the linear map 
$$ Ker \big[(a - \mu.b) : E \to E\big] \to  Ker \big[(a - \mu.b) : E\big/b^N.E  \to E\big/b^N.E\big] $$
is injective.
 \end{lemma}

\parag{Proof} If \ $x \in Ker \, (a - \mu.b)\cap b^N.E$, write \ $x = b^N.y$. Then \ $(a - (\mu-N).b).y = 0$ \ and this implies \ $y = 0$ \ as \ $\mu - N < \lambda_1 +1 - k$. $\hfill \blacksquare$\\

\begin{cor}\label{injective}
Let \ $\mathbb{E}$ \ be a holomorphic family of \ $[\lambda]-$primitive frescos parametrized by a reduced complex space \ $X$. Assume \ $X$ \ connected and note \ $\lambda_1, \dots, \lambda_k$ \ the fundamental invariants of the frescos \ $\mathbb{E}(x)$. Then for \ $\mu \in [\lambda]$ \ and \ $N \geq \mu - \lambda_1 + k$ \ the map induced by the quotient \ $\mathbb{E} \to \mathbb{E}\big/b^N.\mathbb{E}$ \ is injective on the sheaf \ $\mathcal{K}_{\mu}$.
\end{cor}

\parag{proof} Obvious from the remark which begins the section 4.3. $\hfill \blacksquare$

\begin{prop}\label{coherence}
Let \ $\mathbb{E}$ \ be a holomorphic family of \ $[\lambda]-$primitive frescos parametrized by a reduced complex space \ $X$. Define for \ $\mu \in [\lambda]$ \ the \ $\mathcal{O}_X-$module
$$ \mathcal{K}_{\mu} : Ker\left[ (a - \mu.b) : \mathbb{E} \to \mathbb{E}\right].$$
Then this sheaf is a coherent \ $\mathcal{O}_X-$module.
\end{prop}

\parag{proof} Note \ $\mathcal{K}_{\mu}^N$ \ the kernel of \ $a - \mu.b$ \ acting on the sheaf \ $\mathbb{E}\big/b^N.\mathbb{E}$. Then we have the commutative diagram
$$\xymatrix{& \mathcal{K}_{\mu} \ar[d]^{q_{N+1}} \ar[r]^{b^k} & \mathcal{K}_{\mu+k} \ar[d]^{q_N} \\
& \mathcal{K}_{\mu}^{N+1} \ar[r]^{b^k} & \mathcal{K}_{\mu+k}^N } $$
where the horizontal maps are induced by \ $b^k$ \ and the vertical map are the obvious quotient maps. The  corollary \ref{monomes E} gives the equality in \ $\mathcal{K}_{\mu+k}^N $
$$ b^k(\mathcal{K}_{\mu}^{N+1}) = q_N(\mathcal{K}_{\mu+k}) .$$
Now the corollary \ref{injective} gives that the map \ $q_N$ \ is injective for \ $N \geq \mu+ k - \lambda_1$.
As the sheaf \ $\mathcal{K}_{\mu}^N$ \ is \ $\mathcal{O}_X-$coherent for any choice of \ $\mu$ \ and \ $N$ \ because it is the  kernel of a \ $\mathcal{O}_X-$linear map between coherent (in fact locally free) \ $\mathcal{O}_X-$modules, the conclusion follows from the injectivity of \ $q_N$ \ choosing \ $N$ \ large enough. $\hfill \blacksquare$\\

Recall that for a fresco \ $E$ \ we denote \ $\delta(E)$ \ the rank of \ $S_1(E)$ \ the maximal semi-simple sub-module in \ $E$.

\begin{prop}\label{delta constant}
Let \ $\mathbb{E}$ \ be a holomorphic family of \ $[\lambda]-$primitive frescos parametrized by a reduced complex space \ $X$.There exists a dense Zariski open set in \ $X$ \ on which the function
$$ x \mapsto \delta(\mathbb{E}(x)) $$
is locally constant.
\end{prop}

\parag{proof} We may assume that \ $X$ \ is irreducible. Fix \ $\mu $ \ and an integer \ $m \gg 1$ \ such that for any \ $x \in X$ \ the map \ $K_{\mu}(x) \to K_{\mu}^m(x) $ \ is an isomorphism, where \ $K_{\mu}$ \ (resp. \ $K_{\mu}^m(x)$) \ is the kernel of \ $a - \mu.b$ \ acting on \ $\mathbb{E}(x)$ \ (resp. on \ $\mathbb{E}(x)\big/b^m.\mathbb{E}(x)$). This is possible thanks to the proposition \ref{tool1} the lemma \ref{ordre suff.} and the corollary \ref{monomes E}. Choose now a dense Zariski open set \ $X'$ \ in \ $X$ \ on which the \ $\mathcal{O}_X-$linear map
$$ a - \mu.b : \mathbb{E}\big/b^m.\mathbb{E} \to \mathbb{E}\big/b^m.\mathbb{E} $$
has constant rank. Then its kernel \ $\mathcal{K}_{\mu}^m$ \ is a locally free finite type \ $\mathcal{O}_X-$module, and we have for each \ $x \in X'$
$$ \mathcal{K}_{\mu}^m(x) = Ker\left[ a - \mu.b : \left(\mathbb{E}\big/b^m.\mathbb{E}\right)(x) \to \left(\mathbb{E}\big/b^m.\mathbb{E}\right)(x)\right] .$$
But as \ $ \mathbb{E}\big/b^m.\mathbb{E} $ \ is locally free on \ $\mathcal{O}_X$ \ we have for any \ $x \in X$ \ an isomorphism  
$$\left(\mathbb{E}\big/b^m.\mathbb{E}\right)(x) \simeq \mathbb{E}(x)\big/b^m.\mathbb{E}(x) .$$
So for each \ $x$ \ in \ $X'$ \ we have
$$ K_{\mu}(x) \simeq K_{\mu}^m(x) \simeq \mathcal{K}_{\mu}^m(x) .$$
This implies that on \ $X'$ \ we have \ $\delta(\mathbb{E}(x)) = rank(\mathcal{K}_{\mu}^m) $. $\hfill \blacksquare$

\parag{Remark} If we choose \ $X'$ \ smaller in order that the sheaf \ $\mathcal{K}_{\mu}$ \ is also locally free and isomorphic to \ $\mathcal{K}_{\mu}^m$, we see that on this smaller dense Zariski open set in \ $X$ \ we have the free \ $\mathcal{O}_X[[b]]-$module \ $\mathcal{O}_X[[b]].\mathcal{K}_{\mu}$ \ which is contained in \ $\mathcal{S}_1(\mathbb{E})$ \ as it will be defined in the proof of the theorem \ref{hol. ss part} which is given in the next subsection.

\subsection{The proof.}

\begin{lemma}\label{coker S1}
Let \ $E$ \ be a \ $[\lambda]-$primitive fresco with fundamental invariants \ $\lambda_1, \dots, \lambda_k$. Assume that \ $\delta(E) = \delta(F_{k-1}) + 1$ \ where \ $F_{k-1}$ \ is the \ $(k-1)-$th term of the principal J-H. sequence of \ $E$. Then let \ $-\mu$ \ be the root of the Bernstein polynomial of \ $S_1(E)$ \ such that \ $-\mu -1$ \ is not a root of the Bernstein polynomial of \ $S_1(F_{k-1})$ \ and define \ $r : = \mu - \lambda_k$. Then \ $r$ \ is a non negative integer and  we have an exact sequence of semi-simple frescos
$$ 0 \to S_1(F_{k-1}) \to S_1(E) \to b^r.(E\big/F_{k-1}) \to 0 $$
where the first map is the obvious inclusion and the second is  induced by the quotient  map\ $E \to E\big/F_{k-1} \simeq E_{\lambda_k}$.
\end{lemma}

Note that, as \ $-\mu$ \ is a root of the Bernstein polynomial of \ $S_1(E)$ \ there exists \ $j \in [1,k]$ \ such that \ $\mu = \lambda_j+j - \delta(E) \leq \lambda_k +k - \delta(E) $ \ and so  \ $r \leq k - \delta(E)$. 

\parag{proof} As we have \ $S_1(F_{k-1}) = S_1(E) \cap F_{k-1}$ \ (see [B.12] section 2), the kernel of the restriction to \ $S_1(E)$ \ of the quotient map \ $E \to E\big/F_{k-1}$ \ is \ $S_1(F_{k-1})$. As we assume the equality  $\delta(F_{k-1}) = \delta(E)-1$ \ the image of the second map has rank \ $1$ \ and so  is equal to \ $b^r.E_{\lambda_k}$ \ for some non negative integer \ $r$. Now this implies that \ $\lambda_k+r$  \ is the opposite of a root of the Bernstein polynomial of \ $S_1(E)$. The roots of the Bernstein polynomial of \ $S_1(E)$ \ are pairwise distinct (see [B.12] section 2), and the roots of the Bernstein polynomial of \ $S_1(F_{k-1})$ \ are of the form \ $-\nu-1$ \ where \ $- \nu$ \ is a root of the Bernstein polynomial of \ $S_1(E)$. This shows that \ $\mu = \lambda_k + r$ \ is well defined by the condition that \ $-\mu$ \ is the root of the Bernstein polynomial of \ $S_1(E)$ \ such that \ $-\mu - 1$ \ is not a root of the Bernstein polynomial of \ $S_1(F_{k-1})$.$\hfill \blacksquare$\\

\begin{prop}\label{FF}
For each holomorphic family \ $\mathbb{E}$ \ of \ $[\lambda]-$primitive frescos parametrized by a reduced complex space \ $X$ \ there exists a \ $\hat{A}_X-$subsheaf \ $\mathcal{S}_1(\mathbb{E}) \subset \mathbb{E} $ \ with the following properties :
\begin{enumerate}[i)]
\item For a holomorphic map \ $f : Y \to X$ \ we have a natural \ $\hat{A}_Y-$map
$$  f^*(\mathcal{S}_1(\mathbb{E})) \to \mathcal{S}_1(f^*(\mathbb{E}))    $$
where \ $f^*$ \  is the pull-back for the \ $\hat{A}_X-$modules to \ $\hat{A}_Y-$modules\footnote{ which is the same than  the pull-back for the \ $\mathcal{O}_X[[b]]-$modules to \ $\mathcal{O}_Y[[b]]-$modules when we forget the action of \ $a$.}. \\
For \ $X = \{pt\}$ \ parametrizing the \ $[\lambda]-$primitive fresco \ $E$ \ we have \ $\mathcal{S}_1(\mathbb{E})$ \ equal to \ $S_1(E)$.
\item When we have an exact sequence
$$ 0 \to \mathbb{F} \to \mathbb{E} \to \mathbb{E}\big/\mathbb{F} \to 0 $$
of holomorphic families of\ $[\lambda]-$primitive frescos parametrized by \ $X$ \ we have
$$ \mathcal{S}_1(\mathbb{F}) = \mathcal{S}_1(\mathbb{E}) \cap \mathbb{F} .$$
Note the following special case :  if \ $\eta : \mathbb{F} \to \mathbb{E}$ \ is an isomorphism then the above condition says that \ $\eta$ \ induces an isomorphism of \ $\hat{A}_X-$sheaves between \ $\mathcal{S}_1(\mathbb{F})$ \ and \ $\mathcal{S}_1(\mathbb{E})$.
\item We have \ $\mathcal{S}_1(\mathbb{E}) \cap b.\mathbb{E} = b.\mathcal{S}_1(\mathbb{E}) $.
\item The sheaf \ $\mathcal{S}_1(\mathbb{E})\big/b.\mathcal{S}_1(\mathbb{E})$ \ is \ $\mathcal{O}_X-$coherent.
\item On an open set \ $X'$ \ in \ $X$ \ on which the sheaf \ $\mathcal{S}_1(\mathbb{E})\big/b.\mathcal{S}_1(\mathbb{E})$ \ is a locally free finite type \ $\mathcal{O}_X-$module  and on which its inclusion in \ $\mathbb{E}\big/b.\mathbb{E}$ \ is a constant rank map of holomorphic vector bundles,  the natural map 
 $$\mathcal{S}_1(\mathbb{E})(x) \to \mathbb{E}(x)$$ 
  induces  an isomorphism  on a normal semi-simple sub-module of  \ $S_1(\mathbb{E}(x))$,  for each point \ $x \in X'$.
\item There exists a dense Zariski open set \ $X' $ \ in \ $X$ \ such that v) holds.
\end{enumerate}
\end{prop}

\parag{proof} Of course we may assume that \ $X$ \ is connected. So the rank \ $k$ \ and the fundamental invariants \ $\lambda_1, \dots, \lambda_k$ \ of the frescos \ $\mathbb{E}(x)$ \ are independent of \ $x \in X$. Choose \ $\mu \in [\lambda]$ \ such \ that \ $\mu \geq \lambda_k+k-1$ \ and \ $q \in \mathbb{N}$ \ such that \ $q \geq \mu -\lambda_1 + k - 1$. Then define 
 $$\mathcal{S}_1(\mathbb{E}) : = (b^{-q}.\mathcal{O}_X[[b]].\mathcal{K}_{\mu})\, \cap \, \mathbb{E}$$
 First we shall prove that \ $\mathcal{S}_1(\mathbb{E})$ \ is independant of the choices of \ $\mu$ \ and \ $q$. An immediate consequence of the corollary \ref{cas relatif} is  that the sheaf we defined depends only on \ $\mu -q$. So fix now \ $\mu \geq \lambda_k+k-1$; we have for \ $ q \geq \mu - \lambda_1+ k - 1$ 
  $$( b^{-q-1}.\mathcal{O}_X[[b]].\mathcal{K}_{\mu+1}) \, \cap \, \mathbb{E} = (b^{-q}.\mathcal{O}_X[[b]].\mathcal{K}_{\mu}) \, \cap\, \mathbb{E}$$
  as \ $\mathcal{K}_{\mu + 1} = b. \mathcal{K}_{\mu} $ \ thanks to the corollary \ref{cas relatif}, and we have
  $$ (b^{-q-1}.\mathcal{O}_X[[b]].\mathcal{K}_{\mu})\, \cap \, \mathbb{E} = (b^{-q}.\mathcal{O}_X[[b]].\mathcal{K}_{\mu})\, \cap \,\mathbb{E}  $$
thanks again to the corollary \ref{cas relatif}. So we proved that the subsheaf of \ $\mathbb{E}$  \\
 $(b^{-q}.\mathcal{O}_X[[b]].\mathcal{K}_{\mu}) \,  \cap \, \mathbb{E}$ \ is independant of \ $\mu$ and \ $q$ \ satisfying \ $\mu \geq \lambda_k+k-1$ \ and \ $q \geq \mu - \lambda_1+k-1$. The property iii) is then consequence of the equalities 
$$ b.\mathcal{S}_1(\mathbb{E}) = b.\left[(b^{-q-1}.\mathcal{O}_X[[b]].\mathcal{K}_{\mu})\, \cap \, \mathbb{E} \right] =  (b^{-q}.\mathcal{O}_X[[b]].\mathcal{K}_{\mu})\, \cap \,b.\mathbb{E} =  \mathcal{S}_1(\mathbb{E}) \cap b.\mathbb{E}.$$
The proof of i) is standard and ii) is easy  as we have for any \ $\mu$ \  the  equality 
$$\mathcal{K}_{\mu}(\mathbb{F}) = \mathcal{K}_{\mu}(\mathbb{E}) \cap \mathbb{F}$$ 
\smallskip
For the proof of iv) we shall use the following lemma.

\begin{lemma}\label{finite}
In the situation of the previous proposition, the sheaf \ $\mathcal{S}_1(\mathbb{E})$ \ is finitely generated on \ $\mathcal{O}_X[[b]]$.
\end{lemma}

\parag{proof} Fix \ $\mu$ \ large enough and let \ $s_1, \dots, s_p$ \ be a local generator of the sheaf \ $\mathcal{K}_{\mu}$ \ as \ $\mathcal{O}_X-$module. This exists as we proved that this sheaf is  \ $\mathcal{O}_X-$coherent. For \ $q$ \ large enough we shall show that the \ $\mathcal{O}_X-$module \ $\mathcal{F}$ \ generated in \ $b^{-q}.\mathbb{E}$ \ by the sections \ $b^{-j}.s_i, i \in [1,p], j \in [1,q]$ \ is \ $\mathcal{O}_X-$coherent. This is a consequence of the fact that for \ $m \gg q$ \ we have \ $\mathcal{F} \cap b^m.\mathbb{E} = \{0\}$ \ and so the map
$$ \mathcal{F} \longrightarrow b^{-q}.\mathbb{E}\big/ b^m.\mathbb{E} $$
is an isomorphism on the coherent  $\mathcal{O}_X-$submodule generated in \ $b^{-q}.\mathbb{E}\big/ b^m.\mathbb{E}$ \ (which is a locally free finite type \ $\mathcal{O}_X-$module) by the images of the sections  \ $b^{-j}.s_i$ \ where\ $ (i,j) \in [1,p]\times [1,q]$. Then \ $\mathcal{F} \cap \mathbb{E}$ \ which is the kernel of the \ $\mathcal{O}_X-$linear map \ $\mathcal{F} \to b^{-q}.\mathbb{E}\big/\mathbb{E}$ \ between two \ $\mathcal{O}_X-$coherent modules is \ $\mathcal{O}_X-$coherent.\\
Now we have
$$ \mathcal{S}_1(\mathbb{E}) = \mathcal{O}_X[[b]].\mathcal{K}_{\mu} + \mathcal{F}\, \cap \, \mathbb{E} $$
from our definition of \ $\mathcal{F}$, as \ $\mathcal{K}_{\mu}$ \ is a subsheaf of \ $\mathbb{E}$. So the lemma is proved. $\hfill \blacksquare$\\

\parag{End of the proof of the proposition \ref{FF}} Now as \ $\mathcal{S}_1(\mathbb{E})$ \ is a subsheaf of \ $\mathcal{O}_X[[b]]-$modules of the sheaf  \ $\mathbb{E}$ \ which locally free and finite type  on \ $\mathcal{O}_X[[b]]$,  we shall deduce that iv) holds, using iii). We have, near each point in \ $x \in X$, \ $\sigma_1, \dots, \sigma_p$ \ a local generator of the sheaf \ $\mathcal{S}_1(\mathbb{E})$ \ on \ $\mathcal{O}_X[[b]]$. Then we have a surjective \ $\mathcal{O}_X[[b]]-$linear map \ $\Sigma : \mathcal{O}_X[[b]]^p \to \mathcal{S}_1(\mathbb{E}) \subset \mathbb{E}$. Then \ $\Sigma$ \ induced a \ $\mathcal{O}_X-$linear map between \ $\mathcal{O}_X-$coherent sheaves \ $(\mathcal{O}_X)^p \to \mathbb{E}\big/b.\mathbb{E}$ \ with image \ $\mathcal{S}_1(\mathbb{E})\big/b.\mathcal{S}_1(\mathbb{E}) \subset \mathbb{E}\big/b.\mathbb{E}$ \ thanks to property iii) already  proved above. \\
The proof of v) is again standard, using our hypothesis which allows to locally split on \ $\mathcal{O}_X$ \ of the exact sequence
$$ 0 \to \mathcal{S}_1(\mathbb{E})\big/b.\mathcal{S}_1(\mathbb{E}) \to \mathbb{E}\big/b.\mathbb{E} \to \mathbb{E}\big/b\mathbb{E} +  \mathcal{S}_1(\mathbb{E}) \to 0 ,$$
because it implies that  \ $\mathcal{S}_1(\mathbb{E})$ \ is locally  a direct factor of \ $\mathbb{E}$ \ as a \ $\mathcal{O}_X[[b]]-$module. So the tensor product giving the evaluation at \ $x \in X'$ \  gives an isomorphism of the fiber at \ $x$ \ of \ $\mathcal{S}_1(\mathbb{E})$ \ into a normal\footnote{because direct factor over \ $\C[[b]]$.} submodule of \ $\left(b^{-q}.\C[[b]].K_{\mu}(x)\right) \cap \mathbb{E}(x)$ \ where \ $K_{\mu}(x)$ \ denotes the kernel of \ $(a - \mu.b)$ \ acting on \ $\mathbb{E}(x)$. But we have the equality  \ $\left(b^{-q}.\C[[b]].K_{\mu}(x)\right) \cap \mathbb{E}(x) = S_1(\mathbb{E}(x))$ \ thanks to the proposition \ref{tool1}.\\
To prove vi) it is enough to remark that for any morphism of \ $\mathcal{O}_X-$coherent modules there exists a dense Zariski open set \ $X'$ \ of \ $X$\ on which the coherent modules are locally free and the morphism has constant rank. $\hfill \blacksquare$\\

\parag{proof of the theorem \ref{hol. ss part}} We may  assume that \ $X$ \ is irreductible. We shall make an induction on the rank \ $k$ \ of the frescos in the holomorphic  family \ $\mathbb{E}$. In the rank \ $1$ \ case we have a family which is locally isomorphic to a constant family, so we have \ $\mathcal{S}_1(\mathbb{E}) = \mathbb{E}$ \  and the result is clear.\\
 Assume now \ $k \geq 2$ \ and the theorem proved for \ $ k' \leq k-1$. Consider the family \ $\mathbb{F}_{k-1}$ \ given by the rank \ $(k-1)-$th  term in the principal Jordan-H{\"o}lder sequence of \ $\mathbb{E}$. It is holomorphic, thanks to the theorem \ref{J-H. hol.}. Consider now the exact sequence of \ $\hat{A}_X-$sheaves, consequence of properties ii) of the proposition \ref{FF} 
$$ 0 \to \mathcal{S}_1(\mathbb{F}_{k-1}) \to \mathcal{S}_1(\mathbb{E}) \to \mathbb{E}\big/\mathbb{F}_{k-1} .$$
This induced an exact   sequence of  coherent \ $\mathcal{O}_X-$modules (see property ii) iii)  and iv) in \ref{FF})
$$ 0 \to  \mathcal{S}_1(\mathbb{F}_{k-1})\big/b. \mathcal{S}_1(\mathbb{F}_{k-1})  \overset{i}{\to} \mathcal{S}_1(\mathbb{E})\big/b.\mathcal{S}_1(\mathbb{E}) \to \mathbb{E}\big/(\mathbb{F}_{k-1}+b.\mathbb{E}) .$$
Now from the property vi) of \ref{FF} there exists a dense Zariski open set \ $X'$ \ on \ $X$ \ such that on it the family \ $\mathcal{S}_1(\mathbb{F}_{k-1})$ \ is holomorphic, the sheaf \ $\mathcal{S}_1(\mathbb{E})\big/b.\mathcal{S}_1(\mathbb{E})$ \ is locally free and the map \ $i$ \ has constant rank. Now cutting \ $X'$ \ with the dense Zariski open set given by the proposition \ref{delta constant}, we may assume that the rank of \ $\mathcal{S}_1(\mathbb{E}(x))$ \ is constant ; denote this rank \ $\delta(\mathbb{E})$. Then, using the lemma \ref{coker S1} the sheaf \ $ \mathbb{E}\big/(\mathbb{F}_{k-1}+b.\mathbb{E})$ \ is locally free of  rank  equal to \ $0$ \ or \ $1$. Now we shall prove that  the family given by the sheaf \ $\mathcal{S}_1(\mathbb{E})$ \ is a holomorphic family.\\
The first case is when the rank of \ \ $ \mathbb{E}\big/(\mathbb{F}_{k-1}+b.\mathbb{E})$ \  is \ $0$; the conclusion is obvious from the inductive hypothesis.\\
When the rank is \ $1$, in the exact sequence of the lemma  \ref{coker S1} corresponding to the evaluation at \ $x \in X'_1$, the integer  \ $r(x)$ \ is constant because the image of \ $\mathcal{S}_1(\mathbb{E})$ \ in  \ $\mathbb{E}\big/\mathbb{F}_{k-1} \simeq \mathbb{E}_{\lambda_k}$ \ is a \ $\hat{A}_X-$module which is  locally free rank \ $1$  \ $\mathcal{O}_X[[b]]-$module. Then on \ $X'_1$ \ we have an exact sequence of sheaves
$$ 0 \to \mathcal{S}_1(\mathbb{F}_{k-1}) \to \mathcal{S}_1(\mathbb{E}) \to \mathbb{E}_{\lambda_k + r} \to 0 $$
where \ $ \mathbb{E}_{\lambda_k + r} \simeq b^r. (\mathbb{E}\big/\mathbb{F}_{k-1})$ \ is a holomorphic family, thanks to the theorem \ref{J-H. hol.}. Now we conclude that \ $\mathcal{S}_1(\mathbb{E})$ \ induces a holomorphic family on \ $X'_1$ \ using the theorem \ref{deux sur trois}. \\
To finish the proof we have to show that, up to a restriction to a smaller  dense Zariski open set in \ $X$, the map \ $\mathcal{S}_1(\mathbb{E})(x) \to S_1(\mathbb{E}(x))$ \ is an isomorphism for each \ $x$. As we already know that this map is injective and that the image is a normal (semi-simple) submodule of \ $S_1(\mathbb{E}(x))$, it is enough to show that the rank of this image is equal to \ $\delta(\mathbb{E})$. A sufficient condition for that purpose is to know that for \ $\mu \gg 1$ \ the dimension of the image of \ $\mathcal{K}_{\mu}(x)$ \ in \ $K_{\mu}(x)$ \ is at least equal \ $\delta(\mathbb{E})$, thanks to the lemma \ref{rang S1}. But there exists a dense Zariski open set in \ $X$ \ on which the coherent sheaf \ $\mathcal{K}_{\mu}$ \ is locally free  of rank \ $\tilde{\delta}$ \  and locally direct factor of \ $\mathbb{E}$ \ as \ $\mathcal{O}_X-$module. The assertion \ $\tilde{\delta} = \delta(\mathbb{E})$ \ is consequence of the fact that \ $\mathcal{K}_{\mu} \cap \mathbb{F}_{k-1}$ \ has  rank \ $\delta(\mathbb{E}) - 1$ \ from our assumption and the fact that in this second case \ $\mathcal{K}_{\mu} \cap \mathbb{F}$ \ is strictly contained in \ $\mathcal{K}_{\mu}$ \ and the quotient cannot be a torsion sheaf on \ $\mathcal{O}_X $ \ as this quotient is a subsheaf of the locally free \ $\mathcal{O}_X[[b]]-$module \ $b^r. (\mathbb{E}\big/\mathbb{F}_{k-1})$. $\hfill \blacksquare$

  \section{ Change of variable for holomorphic \\ parameters.}

\subsection{Definition and statement of  the theorem.}

%All \ $\C-$algebras have a unit.\\

When we consider a sequence of variables \ $\rho : = (\rho_i)_{i \in \mathbb{N}}$ \ we shall denote by \ $\C[\rho]$ \ the \ $\C-$algebra generated by these variables. Then \ $\C[\rho][[b]]$ \ will be the commutative \ $\C-$algebra of formal power series
$$ \sum_{\nu = 0}^{\infty} \ P_{\nu}(\rho).b^{\nu} $$
where \ $P_{\nu}(\rho)$ \ is an element in \ $\C[\rho]$ \ so a polynomial in \ $\rho_0, \dots, \rho_{N(\nu)} $ \ where \ $N(\nu)$ \ is an integer depending on \ $\nu$. So each coefficient in the formal power serie in \ $b$ \ depends only on a finite number of the variables \ $\rho_i$.

\parag{example} When we shall consider unimodular changes of variable
$$ \theta(a) : = a + \sum_{i=2}^{\infty} \ \theta_i.a^i $$
we shall use the notation \ $\C[\theta]$ \ for the algebra associated to the algebraically indpendant variables \ $\theta_i, i \geq 2$.\\
%For  general changes of variable \ $\theta$ \ we shall use the algebra \ $\C[\theta][\chi(\theta),\chi(\theta)^{-1}]$ \ which which has non constant invertibles elements.

%the subset of \ $\C[[a]]$ \ defined by the conditions \ $\theta_0 = 0$ \ and \ $\chi(\theta) = \theta_1 = 1$ \ and 
% we shall denote by \ $\C_0[\theta]$ \ the subalgebra of \ $\C[\theta]$ \  generated by \ $(\theta_i)_{i \geq 2} $. $\hfill \square$\\

\begin{defn}\label{polyn. depend.}
Let \ $E$ \ be a (a,b)-module and let \ $e(\rho)$ \ be a family of elements in  \ $E$ \ depending on a family of variables \ $(\rho_i)_{i\in \mathbb{N}}$. We say that \ $e(\rho)$ \ {\bf\em depends polynomially on \ $\rho$} \ if there exists a fixed \ $\C[[b]]-$basis \ $e_1, \dots, e_k$ \ of \ $E$ \ such that
$$ e(\rho) = \sum_{j=1}^k \ S_j(\rho).e_j $$
where \ $S_j$ \ is for each \ $j \in [1,k]$ \ an element in the algebra  \ $\C[\rho][[b]]$.
\end{defn}

\parag{Remarks} \begin{enumerate}
\item It is important to remark that when \ $e(\rho)$ \ depends polynomially of \ $\rho$, then \ $a.e(\rho)$ \ also. Then for any \ $u \in \A$, again \ $u.e(\rho)$ \ depends polynomially on \ $\rho$. Note that this uses the stability of the algebra \ $\C[\rho][[b]]$ \ by derivation in \ $b$.
\item It is easy to see that we obtain an equivalent condition on the family \ $e(\rho)$ \ by asking that the coefficient of \ $e(\rho)$ \ are in \ $\C[\rho][[b]]$ \ in a \ $\C[[b]]-$basis  of \ $E$ \  basis depending polynomially of \ $\rho$.
\item The invertible elements in the algebra \ $\C[\rho][[b]]$ \ are exactly those elements with an invertible  constant term in \ $b$ \ in the algebra \ $\C[\rho]$.  When the variables \ $(\rho)_{i\in \mathbb{N}}$ \ are algebraically independant, the invertible elements are  those  with a constant term in \ $b$ \ in \ $\C^*$. This will be the case when we consider unimodular changes of variable : the variables \ $\theta_i, i \geq 2$ \ are algebraically independant. But this is not the case for  general changes of variables.$\hfill \square$
\end{enumerate}

\begin{prop}\label{gen. dep.}
Let \ $E$ \ be a rank \ $k$ \  fresco and let \ $\lambda_1, \dots, \lambda_k$ \ its fundamental invariants. Let \ $e(\rho)$ \ be a family of generators of \ $E$ \ depending polynomially on a family of algebraically independant variables \ $(\rho_i)_{i \in \mathbb{N}}$. Then there exists \ $S_1, \dots, S_k$ \ in \ $\C[\rho][[b]]$ \ such that \ $S_j(\rho)[0] = 1$ \ for each \ $j \in [1,k]$ \ and such that the annihilator of \ $e(\rho)$ \ in \ $E$ \ is generated by the element of \ $\A$
$$ P(\rho) : = (a - \lambda_1.b).S_1(\rho)^{-1} \dots S_{k-1}(\rho)^{-1}.(a - \lambda_k.b).S_k(\rho)^{-1} .$$
\end{prop}

\parag{Proof} The key result to prove this proposition is the rank \ $1$ \ case. In this case we may consider a standard generator \ $e_1$ \ of \ $E$ \ which is a \ $\C[[b]]-$basis of \ $E$ \ and satisfies
$$ (a - \lambda_1.b).e_1 = 0 .$$
Then, by definition, we may write
$$ e(\rho) = S_1(\rho).e_1 $$
where \ $S_1 $ \ is in \ $\C[\rho][[b]]$ \ is invertible in this algebra, so has a constant  non zero constant term \ in \ $b$. So, up to normalizing \ $e_1$, we may assume that \ $S_1(\rho)[0] = 1 $ \ and then define 
 $$P(\rho) : = (a - \lambda_1.b).S_1(\rho)^{-1} .$$
It clearly generates the annihilator of \ $e(\rho)$ \ for each \ $\rho$.\\
Assume now that the result is already proved for the rank \ $k-1 \geq 1$. Then consider the family \ $[e(\rho)]$ \ in the quotient \ $E\big/F_{k-1}$ \ where \ $F_{k-1}$ \ is the rank \ $k-1$ \ submodule of \ $E$ \ in its principal J-H. sequence. Remark first that \ $[e(\rho)]$ \ is a family of generators of
 \ $E\big/F_{k-1}$ \ which depends polynomially on \ $\rho$. This is a trivial consequence of the fact that we may choose a \ $\C[[b]]-$basis \ $e_1, \dots, e_k$ \  in \ $E$ \ such that \ $e_1, \dots, e_{k-1}$ \ is a \ $\C[[b]]$ \ basis of \ $F_{k-1}$ \ and \ $e_k$ \ maps to a standard generator of \ $E\big/F_{k-1}\simeq E_{\lambda_k}$. So the rank \ $1$ \  case gives \ $S_k \in \C[\rho][[b]]$ \ with \ $S_k(\rho)[0] = 1$ \  and such that \ $(a - \lambda_k.b).S_k(\rho)^{-1}.e(\rho)$ \ is in \ $F_{k-1}$ \ for each \ $\rho$. But then it is a family of generators of \ $F_{k-1}$ \ which depends polynomially on \ $\rho$ \ and the inductive assumption allows to conclude. $\hfill \blacksquare$\\

We shall apply this result to the change of variable of a given fresco.

\begin{thm}\label{dependance}
Let \ $E$ \ be a rank \ $k$ \  fresco with fundamental invariants \ $\lambda_1, \dots, \lambda_k$,  and let \ $\theta \in \C[[a]]$ \ be an unimodular change of variable. Let \ $e$ \ be a generator of \ $E$ \ such that 
$$ P^0 : = (a - \lambda_1.b).(S^0_1)^{-1} \dots (S^0_{k-1})^{-1}.(a - \lambda_k.b).(S^0_k)^{-1} $$
where \ $S^0_j \in \C[[b]]$ \ satisfies \ $S^0_j(0) = 1 \quad \forall j \in [1,k]$, generates the annihilator of \ $e$ \ in \ $E$. Denote by \ $s$ \ the family of variables associated to the non constant  coefficients of \ $k$ \ elements \ $S_1, \dots, S_k$ \ in \ $\C[[b]]$ \ and \ $s^0$ \ the value of these variables for our given \ $P^0$.
Then there exist \ $S_1(\theta, s), \dots, S_k(\theta, s) $ \ elements in \ $\C[\theta, s][[b]]$ \  satisfying for all \ $j \in [1,k]$ \  \ $S_j(\id, s^0) = S^0_j, S_j(\theta, s)[0] = 1$ \ and such that the generator \ $e$ \ of \ $\theta_*(E)$ \ has its annihilator in \ $\theta_*(E)$ \ generated by the action on  \ $E$ \ of
$$\Pi(\theta, s) : = (\alpha_{\theta} - \lambda_1.\beta_{\theta}).S_1(\theta, s^0)[\beta_{\theta}]^{-1}.(\alpha_{\theta} - \lambda_2.\beta_{\theta}) \dots (\alpha_{\theta} - \lambda_k.\beta_{\theta}).S_k(\theta, s^0)[\beta_{\theta}]^{-1}. $$
Moreover, for each integer \ $q$ \  there exists an integer \ $m(q)$ \ depending only on \ $q$ \ and  \ $\lambda_1, \dots, \lambda_k$ \ such that for any \ $j \in [1,k]$ \  the coefficient of \ $b^q$ \ in \ $S_j(\theta, s)$ \ is independant of \ $\theta_i$ \  and \ $s_{k.i}$ \ for \ $i \geq m(q)$.
\end{thm}

\parag{Remarks}\begin{enumerate}[i)]
\item As the action of a change of variable of the form \ $\theta(a) = \xi.a$ \ on the annihilator of a generator is obvious, the previous theorem implies the analog result for a general change of variable where we have to replace the algebra \\
 $\C[\theta] = \C[\theta_2, \theta_3, \dots]$ \ by the algebra \ $\C[\theta_2. \theta_3 \dots][\chi(\theta),\chi(\theta)^{-1}]$. 
\item Note that the actions on \ $a$ \ and \ $b$ \ on \ $\theta_*(E)$ \ are given by \ $\alpha_{\theta} : = \theta(a)$ \ and \ $\beta_{\theta} = b.\theta'(a)$ \ so  in \ $\theta_*(E)$ \ the annihilator of the generator \ $\theta(e)$\footnote{this notation is used to recall that we consider here  \ $e$ \ as an element in \ $\theta_*(E)$.} in \ $\hat{A}$ \ (for these  actions of \ $a$ \ and \ $b$ \ on \ $\theta_*(E))$ \  is generated by
$$P(\theta, s^0) : = (a - \lambda_1,b).S_1(\theta,s^0)[b]^{-1} \dots (a - \lambda_k.b).S_k(\theta, s^0)[b]^{-1}.$$
\item Let \ $e(\rho)$ \ be a family of generators of a rank \ $k$ \ fresco \ $E$ \  depending polynomially on a family of variables \ $\rho$, and assume that  we have \ $S_1(\rho), \dots, S_k(\rho)$ \ depending polynomially on \ $\rho$ \  with \ $S_j(\rho)[0] = 1\quad \forall j \in [1,k]$ \ and such that 
 $$P(\rho): = (a - \lambda_1.b).S_1(\rho)^{-1} \dots S_{k-1}(\rho)^{-1}.(a - \lambda_k.b).S_k(\rho)^{-1} $$ 
 satisfies \ $P(\rho).e(\rho) \equiv 0 $ \ in \ $E$. Then we shall have in \ $E$
  $$\Pi(\theta, s^{\rho}).e(\rho) \equiv 0  \quad \forall (\theta,\rho) $$
  where \ $s^{\rho}$ \ are the polynomials in \ $\rho$ \ defining \ $S_1(\rho), \dots, S_k(\rho)$ \ for each \ $\rho$. $\hfill \square$\\
  \end{enumerate}
 
 \subsection{Proof of the theorem \ref{dependance}.}

 We shall prove the theorem by induction on the rank \ $k$ \ of the fresco \ $E$. Let first give  the key of the rank \ $1$  \ case.

\begin{prop}\label{rk 1 polyn. depend.}
We consider unimodular  changes of variable. Let\ $\lambda$ \ be a new variable. There exists an unique (invertible) element \ $S \in \C[\theta,\lambda][[b]]$\  such that, for each complex value of \ $\lambda$ \ the element \ $\varepsilon(\theta,\lambda)$ \ in \ $E_{\lambda} : = \hat{A}\big/\hat{A}.(a - \lambda.b)$ \ defined as \ $\varepsilon(\theta,\lambda) : = S(\theta,\lambda)[\beta_{\theta}].e_{\lambda}$ \  satisfies :
$$ \alpha_{\theta}.\varepsilon(\theta,\lambda) = \lambda.\beta_{\theta}.\varepsilon(\theta,\lambda) $$
where \ $\alpha_{\theta} : = \theta(a)$ \ and \ $\beta_{\theta} : = b.\theta'(a)$.
\end{prop}

\parag{Proof} We shall prove first that there exists, for each integer \ $n$, an invertible element \ $\chi_n \in \C[\theta,\lambda][[b]]$\ such the following identity holds in \ $E_{\lambda}$ :
\begin{equation*}
b^n.e_{\lambda} = \beta_{\theta}^n.\chi_n(\theta,\lambda)[\beta_{\theta}].e_{\lambda} \tag{*}
\end{equation*}
where \ $\beta_{\theta} : = b.\theta'(a)$.\\
Remark that \ $\alpha_{\theta} : = \theta(a)$ \ and \ $\beta_{\theta}$ \ are two \ $\C-$linear endomorphisms  of \ $E_{\lambda}$, compatible with the \ \ $b-$filtration, which depend polynomially on \ $\theta $. Now the equality 
 $$b^n.\C[[b]].e_{\lambda} = \beta_{\theta}^n.\C[[\beta_{\theta}]].e_{\lambda}$$
 implies the existence and the uniqueness of \ $\chi_n$ \ for each \ $n \in \mathbb{N}$ \ for each value given to  \ $\theta$ \ and \ $\lambda$. So it is enough to show that \ $\chi_n$ \ depends polynomially on \ $(\theta,\lambda)$.  Fix an integer \ $p \gg 1$ \ and consider the basis of the vector space \ $V_p : = E_{\lambda}\big/b^p.E_{\lambda}$. We have two basis of \ $V_p$\ given by \ $b^q.e_{\lambda}$ \ and \ $\beta_{\theta}^q.e_{\lambda}$ \ for \ $ q \in [0,p-1]$ \ and they satisfy
$$ \beta^q.e_{\lambda} = b^q.e_{\lambda} + b^{q+1}.V_p \quad \forall q \in [0,p-1] .$$
So the corresponding  base change is triangular with entries equal to  $1$  on the diagonal. Let us show by induction on \ $i \in[1,p]$ \ that the coefficients of the \ $(p-i)-$th column of the corresponding matrix depends polynomially in \ $(\theta,\lambda)$. For \ $ i = 1$, this is clear. Assume it is proved for \ $i \geq 1$ \ and consider the equality
$$ \beta_{\theta}^{i-1}.e_{\lambda} = b^{i-1}.e_{\lambda} +  \sum_{j=i}^{p} \ c_{i-1,j}.b^{j}.e_{\lambda}  $$
in \ $V_p$. Then applying \ $\beta_{\theta} = b.\theta'(a)$ \ to both sides gives
$$ \beta_{\theta}^{i}.e_{\lambda} = b.\theta'(a).b^{i-1}.e_{\lambda} +  \sum_{j=i}^p \ c_{i-1,j}.b.\theta'(a).b^{j}.e_{\lambda} . $$
Now the conclusion follows from the following lemma  and the inductive hypothesis:

\begin{lemma}\label{comput.}
For each integer \ $q$ \ there exists  \ $T_q$ \ in \ $\C[\theta,\lambda][[b]]$ \ with constant term in \ $b$ \  equal to \ $1$ \ and satisfying
$$ \theta'(a).b^q.e_{\lambda} = b^q.T_q(\theta,\lambda).e_{\lambda} $$
in \ $E_{\lambda}$, where \ $e_{\lambda}$ \ is a standard generator of \ $E_{\lambda}$ \ such that \ $(a - \lambda.b).e_{\lambda} = 0 $.
\end{lemma}

\parag{Proof} As it is clear that for each \ $(\theta,\lambda)$ \ the element \ $T_q(\theta,\lambda)$ \ exists and is unique in \ $\C[[b]]$ \ the point is to prove the polynomial dependance in \ $(\theta,\lambda)$. But as \ $\theta'(a)$ \ is affine in \ $\theta$ \ and as the coefficient of \ $b^p$ \ in \ $T_q(\theta,\lambda)$ \ depends only on \ $\theta_2, \dots, \theta_p$, it is enough to prove that there exists \ $S_{p,q}(\lambda)$ \ in \ $\C[\lambda][[b]]$ \ such that
$$ a^p.b^q.e_{\lambda} = b^q.S_{p,q}(\lambda).e_{\lambda}.$$
It is an easy exercice to show that \ $ S_{p,q}(\lambda) = (\lambda+q)\dots (\lambda+q+p-1).b^p $ \  and \ 
$T_q(\theta,\lambda) = 1 + \sum_{j=2}^{\infty} \ j.\theta_j.S_{j-1,q}(\lambda) . \hfill \blacksquare$

\parag{end of the proof of the proposition \ref{rk 1 polyn. depend.}} Remark that we have
$$ \alpha_{\theta}.e_{\lambda} = \lambda.b.e_{\lambda} + \sum_{j=2}^{\infty} L_j.\theta_j.b^j.e_{\lambda} $$
where \ $L_j(\lambda) : = \lambda.(\lambda+1) \dots (\lambda+j-1) $. So using our previous statement, we may write
$$ \alpha_{\theta}.e_{\lambda} = \lambda.\beta_{\theta}.e_{\lambda} + \beta_{\theta}^2.R(\theta,\lambda)[\beta_{\theta}].e_{\lambda}$$
where \ $R \in \C[\theta,\lambda][[b]]$. Now the desired equality, with \ $S \in \C[\theta,\lambda][[b]]$ 
$$ \alpha_{\theta}.S[\beta_{\theta}].e_{\lambda} = \lambda.\beta_{\theta}.S[\beta_{\theta}].e_{\lambda} $$
is equivalent to
\begin{equation*}
 S[\beta_{\theta}].\left[ \lambda.\beta_{\theta}.e_{\lambda} + \beta_{\theta}^2.R[\beta_{\theta}].e_{\lambda}\right] + \beta_{\theta}^2.S'[\beta_{\theta}].e_{\lambda} =  \lambda.\beta_{\theta}.S[\beta_{\theta}].e_{\lambda} 
 \end{equation*}
  where \ $S'$ \ is the derivative in \ $b$ \ of \ $S \in \C[\theta,\lambda][[b]]$ ; and after simplification we find the differential equation :
 $$ S' + R.S = 0 .$$
 Let \ $\tilde{R}(\theta,\lambda)$ \ be the primitive in \ $b$ \  without constant term of \ $R(\theta,\lambda)$. Then we obtain \ $S = exp[-\tilde{R}] \in \C_0[\theta,\lambda][[b]]$ \ and its constant term (in \ $b$) \ is equal to \ $1$, so it is invertible in \ $\C[\theta,\lambda][[b]]$. $\hfill \blacksquare$\\

\parag{Proof of theorem \ref{dependance}} Let us prove first  the rank \ $1$ \ case. So let \ $e$ \  be a generator of \ $E_{\lambda_1}$. We may write \ $e = T.e_1$ \ where \ $e_1$ \ is a standard generator of \ $E_{\lambda_1}$ \ and \ $T \in \C[[b]]$ \ satisfies \ $T(0) = 1$. Then the previous lemma gives an \ $S \in \C[\theta][[b]]$ \ such that \ $S(\theta)[0] \equiv 1$ \ and such that \ $\varepsilon(\theta) : = S(\theta)[\beta_{\theta}].e_1 $ \ satisfies 
 $$ \alpha_{\theta}.\varepsilon(\theta) = \lambda_1.\beta_{\theta}.\varepsilon(\theta).$$
 But the formula \ $(^*)$ \ in the proof of the lemma \ref{rk 1 polyn. depend.} gives an element \ $U \in \C[\theta][[b]]$ \ such that \ $U(\theta)[0] \equiv 1$ \ and  satisfying for each \ $\theta$
 $$ e_1 = U(\beta_{\theta}).e .$$
 So we obtain, if \ $V \in \C[\theta][[b]]$ \ is the inverse of \ $S.U$,
 $$ \varepsilon(\theta) = (S.U)[\beta_{\theta}].e $$
and then
 $$ (\alpha_{\theta} - \lambda_1.\beta_{\theta}).(V(\theta)[\beta_{\theta}])^{-1}.e \equiv 0 $$
 concluding the rank \ $1$ \ case.
 
\parag{Remark} The coefficient of \ $b^q$ \ in \ $U$ \ does not depends on the coefficients of \ $b^r$ \ in \ $T$ \ for \ $r \geq q+1$. $\hfill \square$

\bigskip
 
 Assume now that the case of rank \ $\leq k-1$ \ is proved and consider a rank \ $k$ \ fresco \ $E$. Let \ $F_1$ \ be the rank \ $1$ \ submodule of \ $E$ \ in its principal  J-H. sequence and define \ $G : = E\big/F_1$. Now the image \ $[e]$ \ of \ $e$ \ in \ $G$ \ is a generator of \ $G$ \ with annihilator generated by
$$ Q^0 : = (a - \lambda_2.b).S_2^{-1} \dots S_{k-1}^{-1}.(a - \lambda_k.b).S_k^{-1} .$$
So the induction hypothesis gives \ $S_j(\theta) \in \C[\theta][[b]]$ \ for \ $j \in [2,k]$ \ with \ $S_j(\theta)[0] =1$ \ and such that the annihilator of \ $[\theta(e)]$ \ in \ $\theta_*(G)$ \ is generated by 
$$ Q(\theta) : = (a - \lambda_2.b).S_2(\theta)[b]^{-1} \dots (a - \lambda_k.b).S_k(\theta)[b]^{-1} .$$
This means that
$$ (\alpha_{\theta} - \lambda_2.\beta_{\theta}).S_2(\theta)[\beta_{\theta}]]^{-1} \dots (\alpha_{\theta} - \lambda_k.\beta_{\theta}).S_k(\theta)[\beta_{\theta}]]^{-1}$$
annihilates \ $[e]$ \ in \ $E\big/ F_1$ \ for each unimodular \ $\theta$.\\
Now we have \ $Q(\theta).[\theta(e)] = 0$ \ in \ $\theta_*(E\big/F_1) = \theta_*(E)\big/\theta_*(F_1)$ \  which means that \ $Q(\theta).[\theta(e)]$ \ lies in \ $\theta_*(F_1) = F_1(\theta_*(E)) \simeq E_{\lambda_1}$. So applying now the rank \ $1$ \ case gives the conclusion.

\smallskip

To prove the existence of the integer \ $m(q)$ \ it is enough to remark that in the rank \ $1$ \ case, the coefficient of  \ $b^q$ \ in \ $S_1(\theta,\lambda_1)$ \ depends neither on \ $\theta_l$ \ for \ $l \geq q+1$ \ nor on the coefficient of \ $b^r$ \ in \ $T$ \ for \ $r \geq q+1$ \ using the remark following the proof of the rank \ $1$ \ case. $\hfill \blacksquare$

\subsection{Polynomial dependance in \ $E$.}

Fix now the fundamental invariants \ $\lambda_1, \dots, \lambda_k$ \ for a  fresco. 

\begin{defn}\label{polyn. depend. E}
Consider now a complex valued  function \ $f$ \  defined on a subset \ $X$ \ of the isomorphism classes \ $\mathcal{F}(\lambda_1, \dots, \lambda_k)$ \  frescos with fundamental invariants \ $\lambda_1, \dots, \lambda_k$. We shall say that  \ $f$ \  {\bf depends polynomially on the isomorphism class of the fresco\ $E \in X $} \ if the following condition is satisfied :\\
Let \ $s$ \ be the collection of variables corresponding to the coefficients of \ $k$ \ \'elements \ $S_1, \dots, S_k$ \ in \ $\C[[b]]$ \ satisfying \ $S_j(0) = 1\quad \forall j \in [1,k]$ \ and consider for each value of \ $s$ \ the rank \ $k$ \ $[\lambda]-$primitive  fresco \ $E(s) : = \A \big/ \A.P(s) $ \ where 
 $$P(s) : =  (a - \lambda_1.b)S_1(s)^{-1} \dots (a - \lambda_k.b).S_k(s)^{-1} $$
 where \ $S_1(s), \dots, S_k(s)$ \ correspond to the given values for \ $s$.\\
 Then there exists a polynomial \ $F \in \C[s]$ \ such that for each value of  \ $s$ \ for which \ $E \in X$ \ is isomorphic to \ $\A.\big/\A.P(s)$ \ the value of \ $F(s)$ \ is equal to\ $f([E])$.
 \end{defn}

\begin{cor}\label{concrete}
 Assume that we have a family \ $E(\rho)_{\rho \in \C}$ \  of rank \ $k$ \  frescos with fundamental invariants \ $\lambda_1, \dots, \lambda_k$ \ defined by a familly of generators \ $e(\rho)_{\rho \in \C}$ \ in a given (a,b)-module \ $\mathcal{E}$, such that \ $e(\rho)$ \ polynomially depends on \ $\rho \in \C$. Let \ $ f : \mathcal{F}(\lambda_1, \dots, \lambda_k) \to \C$ \ be a function which depends polynomially on \ $[E]$. Then  the function
$$ \rho \mapsto f([E(\rho)]) $$
is a polynomial.
\end{cor}

\parag{proof} For such a family  \ $E(\rho)_{\rho \in \C}$ \  the theorem \ref{dependance} implies that there exist \ $S_1, \dots, S_k$ \ in \ $\C[\rho][[b]]$ \ such that
\begin{enumerate}[i)]
\item \ $S_j(\rho)[0] \equiv 1 \quad \forall j\in [1,k] $.
\item For each \ $\rho$ \ the annihilator of \ $e(\rho)$ \ in \ $\mathcal{E}$ \ is generated by 
$$ P(\rho) : = (a - \lambda_1.b)S_1(\rho)^{-1} \dots (a - \lambda_k.b).S_k(\rho)^{-1} $$
\end{enumerate}
This gives a polynomial  map \ $\rho \mapsto (s(\rho))$ \ and the map \ $ \rho \mapsto f([E(\rho)])$ \ is equal to \ $\rho \mapsto F(s(\rho))$ \ which is a polynomial. $\hfill \blacksquare$

\parag{Example} Let \ $(s_i)_{i \in \mathbb{N}^*}$ \ a family of algebraically independant variables and let  \ $S(s) : = 1 + \sum_{i=1}^{\infty} \ s_i.b^i \in \C[s][[b]]$. Define \ $E(s) : = \A\big/\A.(a - \lambda_1.b).S(s)^{-1}.(a - \lambda_2.b) $ \ where \ $\lambda_1 > 2$ \ is rational and \ $\lambda_2 : = \lambda_1 + p_1 -1$ \ with \ $p_1 \in \mathbb{N}^*$. Define \ $\alpha(s) : = s_{p_1}$. Then the number \ $\alpha(s)$ \ depends only of the isomorphism class of the fresco \ $E(s)$ \ and defines a function on \ $\mathcal{F}(\lambda_1, \lambda_2)$ \ which depends polynomially on \ $[E] \in \mathcal{F}(\lambda_1, \lambda_2)$.\\
This is an easy consequence of the classification of rank \ $2$ \ frescos (see [B.12]). $\hfill \square$

\section{Construction of quasi-invariant parameters.}

\subsection{Holomorphic parameters.}

Let \ $\mathcal{F}(\lambda_1, \dots, \lambda_k)$ \ be the set of all isomorphism classes of  frescos with rank \ $k$ \ and fundamental invariants \ $\lambda_1, \dots, \lambda_k$. We have a natural map
$$ I : \mathcal{G}(\lambda_1, \dots, \lambda_k) \to \mathcal{F}(\lambda_1, \dots, \lambda_k)$$
from the basis of the versal family which is given by \ $s \mapsto [\mathbb{E}^{taut}(s)]$ \ which is surjective.

\begin{defn}\label{hol. param.}
Let \ $\mathcal{F}_0$ \ be a subset of \ $\mathcal{F}(\lambda_1, \dots, \lambda_k)$ \ and let 
$$ f : \mathcal{F}_0 \to \C $$
be a function. We shall say that \ $f$ \ is a {\bf holomorphic parameter on \ $\mathcal{F}_0$} \ when there exists a locally closed analytic subset \ $\mathcal{G}_0 \subset \mathcal{G}(\lambda_1, \dots, \lambda_k)$ \  containing \ $I^{-1}(\mathcal{F}_0)$ \ and a  holomorphic  map
$$ F : \mathcal{G}_0 \to \C $$
such that for each \ $s \in \mathcal{G}_0$ \ with \ $I(s) \in \mathcal{F}_0$, we have \ $ f(I(s)) = F(s)$.
\end{defn}

We shall say that \ $f$ \ is {\bf polynomial} on \ $\mathcal{F}(\lambda_1, \dots, \lambda_k)$ \  when we may choose for \ $F$ \ a polynomial on the affine space \ $\mathcal{G}(\lambda_1, \dots, \lambda_k)$. \\

\begin{cor}[Corollary of theorem \ref{dependance}.]\label{hol. param. change}
For any change of variable \ $\theta$ \ in \ $ \C[[a]]$ \ the map
$$ \theta_* : \mathcal{F}(\lambda_1, \dots, \lambda_k) \to \mathcal{F}(\lambda_1, \dots, \lambda_k) $$
define by the change of variable \ $\theta$ \ is polynomial in the sense that for any polynomial function
$$ f : \mathcal{F}(\lambda_1, \dots, \lambda_k)  \to \C $$
the function \ $f\circ \theta_* : \mathcal{F}(\lambda_1, \dots, \lambda_k)  \to \C$ \ is again polynomial.\\
This implies that for any holomorphic parameter \ $f : \mathcal{F}_0 \to \C$ \ and any change of variable \ $\theta$ \ the function \ $f\circ \theta_* : (\theta_*)^{-1}(\mathcal{F}_0) \to \C$ \ is again a holomorphic parameter.
\end{cor}

This is a consequence of the theorem \ref{dependance} applied to the tautological family we constructed on \ $\mathcal{G}(\lambda_1, \dots, \lambda_k)$ : it gives, choosing  locally a \ $k-$admissible generator, a polynomial map \ $\tilde{\theta}$  which (locally) lift to \ $\mathcal{G}(\lambda_1, \dots, \lambda_k)$ \  the map \ $\theta$ \ on \ $\mathcal{F}(\lambda_1, \dots, \lambda_k)$. But, of course, this map is not canonical and may not gobalize.\\

Fix \ $1 \leq i \leq j \leq k $ \ and let
$$ g_{i,j} : \mathcal{F}(\lambda_1, \dots, \lambda_k) \to \mathcal{F}(\lambda_i, \dots, \lambda_j) $$
the map which sends \ $[E]$ \ to \ $[F_j\big/F_{i-1}]$ \ where \ $(F_h)_{h \in [1,k]}$ \ is the principal J-H. sequence of \ $E$. It is a consequence of the uniqueness of the principal J-H. sequence of a \ $[\lambda]-$primitive fresco that this map is well defined. Then the following obvious proposition will be used to produce inductively holomorphic parameters.

\begin{prop}\label{construct. hol. param.}
Let \ $\mathcal{F}_0$ \ be a subset of \ $\mathcal{F}(\lambda_1, \dots, \lambda_j)$ \ and
$$ f : \mathcal{F}_0 \to \C $$
be a holomorphic parameter (resp. polynomial parameter), then
$$ f\circ g_{i,j} : g_{i,j}^{-1}(\mathcal{F}_0) \to \C  $$
is a holomorphic parameter (resp. a polynomial parameter).
\end{prop}

\bigskip

\begin{defn}\label{q-invariant}
Let \ $\mathcal{F}_0$ \ be a subset of \ $\mathcal{F}(\lambda_1, \dots, \lambda_k)$ \ which is invariant by any change of variable. We shall say that a holomorphic parameter
$$ f : \mathcal{F}_0 \to \C$$
is quasi-invariant of weight \ $w \in \mathbb{N}$ \ if for any change of variable \ $\theta$ \ we have
$$ f\circ \Theta = \chi( \theta)^w.f  $$
on \ $\mathcal{F}_0$. We shall say {\bf invariant} when \ $w = 0 $.
\end{defn}

The  next  proposition is the first basic stone from which we shall build a lot of quasi-invariant holomorphic parameters (in fact polynomial).

\begin{prop}\label{theme rg 2} Assume that \ $p_1 : = \lambda_2 - \lambda_1 +1$ \ is in \ $\mathbb{N}^*$. Let \ $\alpha : \mathcal{F}(\lambda_1,\lambda_2) \to \C $ \ the map defined by the coefficient of \ $b^{p_1}$ \ in a standard presentation of a rank \ $2$ \ fresco. Then \ $\alpha$ \ is a quasi-invariant polynomial parameter with weight \ $p_1$.
\end{prop}

\parag{proof} As \ $\alpha(E) \not= 0$ \ is a necessary and suffisant condition in order that \ $[E]$ \ in \ $ \mathcal{F}(\lambda_1,\lambda_2) $ \ is a theme, the polynomial map
$$ \theta \mapsto \alpha(\theta_*(E)) $$
for a given \ $E$, which is a polynomial in \ $\theta$,  is either identically  zero when \ $E$ \ is semi-simple, either never zero when \ $E$ \ is a theme, thanks to the proposition \ref{stability}. But a polynomial on an affine space which never vanishes is constant, so we have \ $\alpha(\theta_*(E))= \alpha(E)$ \ for any unimodular \ $\theta$. It is then clear that we have a quasi-invariant polynomial parameter of weight \ $p_1$. $\hfill \blacksquare$\\

%\begin{prop}\label{ss rk 4}
%Let \ $\mathcal{F}_{ss}(\lambda_1, \dots, \lambda_4)$ \ be the subset of \ $\mathcal{F}(\lambda_1, \dots, \lambda_4)$ \ of isomorphism class of semi-simple frescos. Then the function
%$$ \rho : \mathcal{F}_{ss}(\lambda_1, \dots, \lambda_4) \to \C$$
%define in the subsection ?? is a quasi-invariant of weight \ $1$.
%\end{prop}

From the result of  [B.12] we may also build another quasi-invariant polynomial parameter.

\begin{prop}\label{alpha-inv.}
Let \ $\mathcal{F}_0(\lambda_1, \dots, \lambda_k)$ \ be the subset of \ $\mathcal{F}(\lambda_1, \dots, \lambda_k)$ \ of isomorphism classes of \ $[\lambda]-$primitive frescos \ $E$ \ such that \ $F_{k-1}$ \ and \ $E\big/F_1$ \ are semi-simple, where \ $(F_j)_{j\in [1,k]}$ \ is the principal J-H. sequence of \ $E$. Then the function
$$ \delta : \mathcal{F}_0 \to \C $$
given by the \ $\alpha-$invariant built in [B.12] section 3.2  is a quasi-invariant polynomial parameter with weight \ $p_1 + \dots + p_{k-1}$.
\end{prop}

\parag{proof} The proof is almost the same than in the previous proposition because for such a fresco \ $E$ \ we have the equivalence between \ $E$ \ semi-simple and \ $\delta(E) = 0$ : \\ 
the polynomial function
$$ \theta \mapsto \delta(\theta_*(E)) $$
is either identically zero or never vanishes,  as semi-simplicity is stable by change of variable thanks to the proposition \ref{stability}. So it has to be constant by any unimodular change of variable for any given  \ $[E] \in \mathcal{F}_0(\lambda_1, \dots, \lambda_k)$.\\
Here we have to use for \ $\mathcal{G}_0$ \ the subset defined in \ $\mathcal{G}(\lambda_1, \dots, \lambda_k)$ \ by the condition that the terms \ $F_1$ \ and \ $F_{k-1}$ \ of the principal J-H. sequence are semi-simple, which is a closed analytic (in fact polynomial) condition thanks to [B.12]. $\hfill \blacksquare$\\

So our aim  now is to prove the theorem \ref{second.parm.} which produces another basic stone for building quasi-invariant holomorphic parameters.

  \subsection{The sharp filtration.}
  
  Let \ $E$ \ be a rank \ $k$ \  fresco. Denote
  $$ 0 = F_0 \subset F_1 \subset \dots \subset F_{k-1} \subset F_k = E $$
  the principal J-H. sequence of \ $E$.
  
  \begin{defn}\label{sharp filt.}
  Define for \ $n \in \mathbb{N}$ \ and \ $h \in [0,k-1]$
  $$ \Phi_{n.k + h} : = b^n.F_{k-h} + b^{n+1}.E .$$
  The filtration of \ $E$ \ by the vector spaces  \ $(\Phi_{\nu})_{\nu \geq 0}$ \ will be called the {\bf sharp filtration of \ $E$}.
  \end{defn}

  \begin{lemma}\label{propr.}
  The sharp filtration has the following properties :
  \begin{enumerate}[i)]
  \item it is strictly decreasing,  $\Phi_0 = E$ \ and \ $\cap_{\nu \geq 0} \ \Phi_{\nu} = \{0\} $. 
  \item For \ each \ $\nu \in \mathbb{N}$ \ the dimension of the complex vector space \ $\Phi_{\nu}\big/\Phi_{\nu+1}$ \ is \ $1$ ; it is generated by the image of \ $b^n.e_{k-h}$ \ for \ $\nu = n.k + h$ \ with \ $h \in [0,k-1]$.
  \item We have \ $a.\Phi_{\nu} \subset \Phi_{\nu+1}$ \ and \ $b.\Phi_{\nu} \subset \Phi_{\nu+k}$.
  \item The sharp filtration is invariant by any change of variable.
  \end{enumerate}
  \end{lemma}  
  
  \parag{proof}  We leave to the reader as an exercice the proof of the properties  i) to iii). The property iv) is consequence of the fact that the principal J-H. sequence of \ $E$ \ is invariant by any change of variable and also the filtration \ $\big(b^m.E, m \in \mathbb{N}\big)$ \ because we defined the sharp filtration only in term of these subspaces.\\
  
  The next lemma is less obvious.
  
  \begin{lemma}\label{non triv.}
  The sharp filtration satifies also the following properties :
  \begin{enumerate}[a)]
  \item For \ $k : = rk(E)$ \ we have \ $a^k.\Phi_{\nu} \subset \Phi_{\nu+2k-1} \quad \forall \nu \in \mathbb{N}$.
  \item If \ $b.x \in \Phi_{\nu+k}$ \ then \ $x \in \Phi_{\nu} \quad \forall \nu \in \mathbb{N}$.
  \item If \ $a.x \in \Phi_{\nu+k}$ \ then \ $x \in \Phi_{\nu} \quad \forall \nu \in \mathbb{N}$.
  \end{enumerate}
  \end{lemma}
  
  \parag{proof} We leave to the reader the easy verifications of the properties b) and c). Note that they are optimal : take \ $x : = e_1$ \ the standard generator of \ $F_1$.\\
  Let \ $e_1, \dots, e_k$ \ be a standard \ $\C[[b]]-$basis of \ $E$. So \ $e : = e_k$ \ generates \ $E$ \ as a \ $\A-$module and we have \ $(a - \lambda_j.b).e_j = S_{j-1}.e_{j-1}$ \ for each \ $j \in [1,k]$ \ with the convention \ $e_0 : = 0$ \ and where \ $S_j \in \C[[b]]$ \ satisfies \ $S_j(0) = 1$ \ for \ $j \in [1,k-1]$. It will be convenient to define \ $\mu_j : = S'_j(0)$ \ for \ $j \in [1,k-1]$.\\
  We want in fact to show that we have for each \ $n \in \mathbb{N}$ \ and each \ $j \in [1,k]$ 
  $$ a^k.b^n.e_j \in \Phi_{n.k +k-j+2k-1} $$
  because \ $b^n.e_j$ \ induces a basis of \ $\Phi_{n.k+k-j}\big/\Phi_{n.k+k-j+1} = \Phi_{(n+3).k -j-1}$, for each \ $n \in \mathbb{N}$ \ and \ $j \in [1,k]$ \ and we may use the  property iii). As an easy induction on \ $k $ \ shows that we have for any integers \ $k$ \ and \ $n$
  $$ a^k.b^n = b^n.a^k + n.k.b^{n+1}.a^{k-1} + b^{n+2}.\A $$
  we see that it is enough to show that we have \ $a^k.e_j \in \Phi_{k -j+ 2k-1} = \Phi_{3k-j-1}$ \ for each \ $j \in [1,k]$ \ and \ $a^{k-1}.e_j \in \Phi_{2k-j-1}$. \\
  The second point is obvious as \ $e_j$ \ is in \ $\Phi_{k-j}$ \ thanks to the property iii).\\
  For the first point  we shall prove that \ $a^j.e_j \in \Phi_{2k-1}$ \ for each \ $j \in [1,k]$.\\
  As \ $a.e_1 = \lambda_1.b.e_1 \in \Phi_{k+k-1}$ \ because \ $e_1 \in \Phi_{k-1}$, the property is true for \ $j = 1$. Assume that it is proved for \ $j \geq 1$ \ and consider \ $a^{j+1}.e_{j+1}$. We have
  $$ a.e_{j+1} = \lambda_{j+1}.b.e_{j+1} + e_j + \mu_j.b.e_j + b^2.F_{k-j} \quad ( F_{k-j}  \subset \Phi_{k-j} )$$
  so, we obtain
  $$ a^{j+1}.e_{j+1} = a^j.e_j +a^j.( \lambda_{j+1}.b.e_{j+1} + \mu_j.b.e_j) \quad modulo \quad \Phi_{3k -j} \subset \Phi_{2k-1} .$$
  Now \ $a^j.b.e_{j+1}$ \ is in \ $\Phi_{j+k+k-j-1}$ \ and \ $ a^j.b.e_j \in \Phi_{j+k+k-j} $;  so \ $ a^{j+1}.e_{j+1} -  a^j.e_j $ \
  lies in \ $\Phi_{2k-1}$. Then \ $a^k.e_j = a^{k-j}.a^j.e_j$ \ is in \ $\Phi_{2k-1 + k-j} = \Phi_{3k-j-1}.\hfill \blacksquare$\\
  
  \subsection{Change of variable in rank 3.}
  
\begin{lemma}\label{rk 1}
Let \ $\theta(a) : = a + \theta_2.a^2 + \theta_3.a^3 + \dots $ \ an unimodular  change of variable, let \ $\alpha : = \theta(a)$ \ and \ $\beta : = b.\theta'(a)$, and fix \ $\mu > 1$ \ in \ $\mathbb{Q}$. Then there exists an unique \ $S_{\mu} \in \C[[b]]$ \ with \ $S_{\mu}(0) = 1$ \ such that \ $\varepsilon_{\mu} : = S_{\mu}(\beta).e_{\mu} $, where \ $e_{\mu}$ \ is the standard generator of \ $E_{\mu}$, satisfies
$$ (\alpha - \mu.\beta).\varepsilon_{\mu} = 0 .$$
Moreover \ $S'_{\mu}(0) = \theta_2.\mu.(\mu-1)$.
\end{lemma}

\parag{proof} We want to prove that the equation
$$ (\alpha - \mu.\beta).S_{\mu}(\beta).e_{\mu} = S_{\mu}(\beta).(\alpha - \mu.\beta).e_{\mu} + \beta^2.S_{\mu}'(\beta).e_{\mu} = 0 $$
has an unique solution such that \ $S_{\mu}(0) = 1$. But we know that \ $E_{\mu} = \C[[\beta]].e_{\mu}$ \ as the \ $\C[[\beta]]-$module. Then put \ $\alpha.e_{\mu} = \mu.\beta.T(\beta).e_{\mu}$ \ where \ $T(0) = 1$ \  \ leads to the equation
$$ \beta.S_{\mu}'(\beta) + \mu.S_{\mu}(\beta).(T(\beta)-1) = 0 $$
which has clearly an unique solution in \ $\C[[\beta]]$ \  with constant term equal to \ $1$, using the recursion formula (with the obvious notations for the coefficients of \ $S_{\mu}$ \ and \ $T$)
$$ (n+ 1).s_{n+1} = - \mu.\sum_{j=1}^n \  s_{n+1-j}.t_j  \quad \forall n \geq 0.$$
The term of degree \ $1$ \ in \ $ \beta$ \ is given by \ $s_1 = -\mu.t_1$. But
$$ \alpha.e_{\mu} = a.e_{\mu} + \theta_2.a^2.e_{\mu} + b^3.E_{\mu} $$
and so \ $\alpha.e_{\mu} = \mu.b.e_{\mu} + \theta_2.\mu.(\mu+1).b^2.e_{\mu} + b^3.E_{\mu}$. Now we have \ $\beta^3.E_{\mu} = b^3.E_{\mu}$ \ and \  $\beta.e_{\mu} = b.e_{\mu} + 2.\theta_2.\mu.b^2.e_{\mu} + b^3.E_{\mu} $ \ which implies \ $b.e_{\mu} = \beta.e_{\mu} - 2.\theta_2.\mu.\beta^2.e_{\mu} + \beta^3.E_{\mu}  $ \ and also \ $b^2.e_{\mu} = \beta^2.e_{\mu} + \beta^3.E_{\mu}$. So finally
$$ \alpha.e_{\mu} = \mu.\big[\beta.e_{\mu} - 2.\theta_2.\mu.\beta^2.e_{\mu} \big] + \theta_2.\mu.(\mu+1).\beta^2.e_{\mu} + \beta^3.E_{\mu} $$
and so \ $t_1 =  \theta_2.(1-\mu)$ \ and \ $s_1 =  \theta_2.\mu.(\mu-1)$. $\hfill \blacksquare$

  \begin{prop}\label{ordre 3 rk 3}
  Let \ $E$ \ be a rank \ $3$ \  fresco with fundamental invariants \ $\lambda_1, \lambda_2, \lambda_3$. We assume that \ $p_1 : = \lambda_2 - \lambda_1 + 1$ \ and \ $p_2 : = \lambda_3 - \lambda_2 + 1$ \ are different from \ $1$.  Assume that \ $E$ \ is generated by \ $e$ \ such that the element in \ $\A$ 
   $$ P : = (a - \lambda_1.b).S_1^{-1}.(a - \lambda_2.b).(1 + v.b^{p_2})^{-1}.(a - \lambda_3.b)$$
  generates the annihilator of \ $e$, where \ $S_1(0) = 1$, and where the complex number \ $v$ \ is \ $0$ \ when \ $p_2$ \ is not an integer \ $\geq 2$. Let \ $\theta(a ) = a + a^3.\tau(a)$ \ be a change of variable. Then there exists a generator \ $\varepsilon$ \ of \ $\theta_*(E)$ \ with annihilator generated by 
   $$(a - \lambda_1.b).T_1^{-1}.(a - \lambda_2.b).(1 + v.b^{p_2})^{-1}.(a - \lambda_3.b)$$
    with \ $T_1(0) = 1$ \ and \ $T_1'(0) = S_1'(0)$.
  \end{prop}
  
  \parag{Remarks} \begin{enumerate}
  \item In fact the only assumption really made in the previous lemma is the fact that \ $p_1 $ \ and \ $p_2 $ \ are different from \ $1$, because any rank \ $3$ \ fresco with this condition admit a generator with such an annihilator, thanks to the classification of rank \ $2$ \ frescos.
  \item For \ $p_1 = 1$ \ the number \ $S_1'(0)$ \ is the parameter of the rank \ $2$ \ fresco \ $F_2$, so the assertion is obvious.
  \item When \ $p_1 \not= 1$ \ and  the coefficient \ $S_1'(0)$ \ is determined by the isomorphism class of \ $E$, then  the classification of rank \ $3$ \  frescos shows that the isomorphism class of \ $E$ \ is completely  determined by the parameters of \ $F_2$ \ and \ $E\big/F_1$ \ and the coefficient of \ $b$ \ in \ $S_1$. As the parameters of \ $F_2$ \ and \ $E\big/F_1$ \ are invariant by any unimodular change of variable (i.e. \ $\theta_1 = 1$), the previous lemma shows that in these cases we have \ $\theta_*(E) \simeq E$ \ for such a change of variable. This will reduce  mainly our study of the changes of variable  to the special case of the form \ $\theta(a) = a + \tau.a^2$ \ for  rank \ $3$ \ frescos.
  \end{enumerate}
  
  \parag{proof} Define \ $e_3 : = e, e_2 : = (1 +v.b^{p_2})^{-1}.(a - \lambda_3.b).e$ \ and 
   $$e_1 : = S_1^{-1}.(a - \lambda_2.b).(1 + v.b^{p_2})^{-1}.(a - \lambda_3.b).e  S_1^{-1}.(a - \lambda_2.b).e_2 .$$
    Then we shall construct a corresponding \ $\C[[\beta]]-$basis \ $\varepsilon_1, \varepsilon_2, \varepsilon_3$ \ such that
  \begin{align*}
  & (\alpha - \lambda_1.\beta).\varepsilon_1 = 0, \\
  & (\alpha - \lambda_2.\beta).\varepsilon_2 = T_1(\beta).\varepsilon_1 \quad {\rm and} \\
  & (\alpha - \lambda_3.\beta).\varepsilon_3 = (1 + v.\beta^{p_2}).\varepsilon_2 .
  \end{align*}
  The first step is easy, thanks to the lemma \ref{rk 1} ; up to a non zero constant we must choose
  \begin{equation*}
  \varepsilon_1 : = S_{\lambda_1}(\beta).e_1 . \tag{1}
  \end{equation*}
   As the principal J-H. sequence of \ $E$ \ is invariant by any change of variable, we must have \ $\varepsilon_2$ \ which is a generator of \ $F_2$. But, up to \ $F_1$ \ and a non zero  multiplicative constant, \ $\varepsilon_2$ \ has to coincide with \ $\varepsilon_2^0 : = S_{\lambda_2}(\beta).e_2 $; so we shall define \ $\varepsilon_2 : = \varepsilon_2^0 + W(\beta).\varepsilon_1$. We shall define \ $U$ \ by the relation \ $(\alpha - \lambda_2.\beta).\varepsilon_2^0 = U(\beta).\varepsilon_1$ \ where we normalize \ $\varepsilon_1$ \ and \ $e_1, e_2,e_3$ \ by the same non zero constant if necessary to have \ $U(0) = 1$.\\
   Now, with the same argument, we must have
   $$\varepsilon_3 : = \rho.S_{\lambda_3}(\beta).e_3 + X(\beta).\varepsilon_2^0 + Y(\beta).\varepsilon_1$$
   where \ $\rho $ \ is in \ $\C^*$. Now the last equation \ $ (\alpha - \lambda_3.\beta).\varepsilon_3 = (1 + v.\beta^{p_2}).\varepsilon_2$ \ gives
   \begin{align*}
   & \rho.(\alpha - \lambda_3.\beta).S_{\lambda_3}(\beta).e_3 + (\lambda_2-\lambda_3).\beta.X(\beta).\varepsilon_2^0 + \beta^2.X'(\beta).\varepsilon_2^0 + X(\beta).U(\beta).\varepsilon_1 \\
   &  + (\lambda_1-\lambda_3).\beta.Y(\beta).\varepsilon_1 + \beta^2.Y'(\beta).\varepsilon_1 =  (1 + v.\beta^{p_2}).\varepsilon_2^0 + (1 + v.\beta^{p_2}).W(\beta).\varepsilon_1 . \tag{@}
   \end{align*}
  As we have \ $a^3.E + b.a^2.E \subset \Phi_5 $ \ thanks to lemma \ref{non triv.} and \ $S_{\lambda_3}(\beta) = 1 + \beta^2.\C[[\beta]]$ \ using \ $\theta_2 = 0$ \ and lemma \ref{rk 1}, when we compute \ $(\alpha - \lambda_3.\beta).S_{\lambda_3}(\beta).e_3 $ \ modulo \ $\Phi_5$, we find only \ $(a - \lambda_3.b).e_3 = e_2 \quad {\rm modulo} \quad \Phi_5$. But we have also \ $\varepsilon_2^0 = e_2$ \ modulo  $\Phi_5$,  and \ $\beta.\varepsilon_1 \in \Phi_5$, so our equation above reduced modulo \ $\Phi_5$ \ to 
  $$ \rho.e_2 - (p_2-1).X(0).b.e_2 + X(0).\varepsilon_1 = e_2 + W(0).\varepsilon_1 \quad {\rm modulo} \quad \Phi_5 .$$
  So we must have \ $\rho = 1, X(0) = W(0)$ \ and \ $X(0) = 0$ \ as \ $p_1 \not= 1$, looking successively in \ $\Phi_1\big/\Phi_2, \Phi_2\big/\Phi_3$ \ and \ $\Phi_4\big/\Phi_5$. Then we obtain that \ $W(0) = 0$ \ and then
  $$(\alpha - \lambda_2.\beta).\varepsilon_2 = U(\beta).\varepsilon_1 + \beta^2.F_1 \subset \Phi_7 .$$
   So we have \ $T_1'(0) = U'(0)$. But now
  \begin{align*}
  & (\alpha - \lambda_2.\beta).\varepsilon_2 = S_{\lambda_2}(\beta).(\alpha - \lambda_2.\beta).e_2 + \beta^2.F_2 \subset \Phi_7 \\
  & \qquad = (a - \lambda_2.b).e_2 + \Phi_6  \\
  & \qquad = S_1(b).e_1 + \Phi_6 = \varepsilon_1 + U'(0).\beta.\varepsilon_1 + \Phi_6
  \end{align*}
  where we uses the fact that \ $a^3.F_2$ \ and \ $b.a^2.F_2$ \ are in \ $\Phi_6$. So as 
  $$e_1 = \varepsilon_1 \quad {\rm modulo} \quad \Phi_8$$
  we conclude that \ $T_1'(0) = U'(0) = S_1'(0)$.\\
  Of course, we know "a priori" that such a \ $\C[[\beta]]-$basis \ $\varepsilon_1, \varepsilon_2, \varepsilon_3$ \  exists thanks to the classification of rank \ $2$ frescos and the invariance by an unimodular change of variable of the parameter of a rank \ $2$ \ $[\lambda]-$primitive fresco. $\hfill \blacksquare$
  
  \begin{prop}\label{chgt var. rk. 3}
   Let \ $E$ \ be a rank \ $3$ \  fresco with fundamental invariants \ $\lambda_1, \lambda_2, \lambda_3$. We assume that \ $p_1 : = \lambda_2 - \lambda_1 + 1$ \ and \ $p_2 : = \lambda_3 - \lambda_2 + 1$ \ are different from \ $1$. Assume that \ $E$ \ is generated by \ $e$ \ such that 
   $$ P : = (a - \lambda_1.b).S^{-1}.(a - \lambda_2.b).(1 + v.b^{p_2})^{-1}.(a - \lambda_3.b)$$
  generates the annihilator of \ $e$, where \ $S(0) = 1$ \ and where the complex number \ $v$ \ is \ $0$ \ when \ $p_2$ \ is not an integer \ $\geq 2$. Let \ $\theta(a ) = a + \tau.a^2$ \ with \ $\tau \in \C$. Then there exists a generator \ $\varepsilon$ \ of \ $\theta_*(E)$ \ with annihilator generated by 
   $$(a - \lambda_1.b).T_1^{-1}.(a - \lambda_2.b).(1 + v.b^{p_2})^{-1}.(a - \lambda_3.b)$$
    with \ $T(0) = 1$ \ and \ $T'(0) = S'(0) - \tau.L$ \ where \ $L : =   (p_1-1).(p_1+p_2 -1) $.
  \end{prop}
  
  Note that the result is also true for \ $p_1 = 1$ \ as \ $S_1'(0)$ \ is the parameter of \ $F_2$ \ which is invariant by any unimodular change of variable.\\
    
  \parag{proof} We shall follow the same line than in the proof of the previous proposition. Now we have \ $\alpha : = a + \tau.a^2$ \ and \ $\beta : = b + 2\tau.b.a $. Define
  \begin{align*}
  & s_i : = \tau.\lambda_i.(\lambda_i - 1) \quad i = 1,2,3 ; \\
  & \varepsilon_1 : = S_{\lambda_1}(\beta).e_1 , \quad \varepsilon^0_2 : = S_{\lambda_2}(\beta).e_2  ;\\
  & \varepsilon_3 : = \rho.S_{\lambda_3}(\beta).e_3 + X(\beta).\varepsilon^0_2 + Y(\beta).\varepsilon_1 .
    \end{align*} 
   We  we want to define \ $\varepsilon_2 : = \varepsilon^0_2 + W(\beta).\varepsilon_1 $ \ in order to have 
    \begin{equation*}
    (\alpha - \lambda_3.\beta).\varepsilon_3 = (1 + v.\beta^{p_2}).\varepsilon_2 . \tag{@@}
    \end{equation*}
    We also define \ $U \in \C[[\beta]]$ \ by the equation \ $(\alpha - \lambda_2.\beta).\varepsilon_2  = U(\beta).\varepsilon_1$ \ with \ $U(0) = 1$ \ (after a normalization if necessary of the initial basis \ $e_1, e_2, e_3$). Then
  \begin{align*}
  & \varepsilon_1 = e_1 + s_1.b.e_1 + \Phi_8 \\
  & \varepsilon_2^0 = e_2 + s_2.(b.e_2 + 2\tau.b.e_1) +  \Phi_7
 \end{align*}
 Now the  equation \ $(@ @)$ \ gives modulo \ $\Phi_6$ \ using
 \begin{align*}
 & S_{\lambda_3}(\beta).(\alpha - \lambda_3.\beta).e_3 = (1 + s_3.\beta).\big[(a - \lambda_3.b) + \tau.(a - 2\lambda_3.b).a\big].e_3  + \Phi_6 \\
 & \qquad = (1 + s_3.\beta).\big[ e_2  + \tau.(a - 2\lambda_3.b)(\lambda_3.b.e_3 + e_2)\big] + \Phi_6 \\
 & \qquad =  e_2 + s_3.e_2 + \tau\big[\lambda_3.b.e_2 + \lambda_2.b.e_2 + e_1 + S'(0).b.e_1  - 2\lambda_3.b.e_2\big]  + \Phi_6
 \end{align*}
 we have, from the relation \ $ (\alpha - \lambda_3.\beta).\varepsilon_3 = (1 + v.\beta^{p_2}).(\varepsilon^0_2 + W(\beta).\varepsilon_1)$
  \begin{align*}
 &  \rho.e_2 + \rho.\tau\big[ (\lambda_2 - \lambda_3).b.e_2 + e_1 + S'(0).b.e_1 \big] -(p_2-1).X(0).\beta.\varepsilon_2^0 + X(0).\varepsilon_1 +  \\
 & \quad  + (X.U)'(0).\beta.\varepsilon_1 - (p_1+p_2-2).\beta.Y(0).\varepsilon_1 =  \varepsilon^0_2 + W(0).\varepsilon_1 + W'(0).\beta.\varepsilon_1 + \Phi_6.
\end{align*}
 Using  the basis \ $e_1, e_2, e_3$ \ this gives
 \begin{align*}
 & \rho.e_2  - e_2 = 0 \quad {\rm in} \quad \Phi_1\big/\Phi_2 \quad {\rm so} \quad \rho = 1 \\
 &  \tau.e_1 + X(0).e_1 - W(0).e_1 = 0 \quad {\rm in} \quad \Phi_2\big/\Phi_3 \quad  {\rm so} \\
 &  W(0) = X(0) + \tau  \\
  & {\rm and \  modulo} \ \Phi_5  \quad {\rm we \ obtain} \\
 & s_3 +  \tau.(\lambda_2 - \lambda_3) - (p_2-1).X(0) - s_2 = 0 \quad   {\rm and \ so} \\
 &   (p_2-1).X(0)  =  - \tau \big[ (p_2-1) + \lambda_2.(\lambda_2-1) - \lambda_3.(\lambda_3-1)\big] 
  \end{align*}
  This implies 
   $$W(0) = \tau.(2\lambda_2 + p_2 - 2)$$
   Now we are looking for the value of  \ $T'(0)$ \ where \ $T$ \ is defined by the relation
   $$ (\alpha - \lambda_2.\beta).\varepsilon_2 = T(\beta).\varepsilon_1 .$$
   As we have
   \begin{align*}
   & (\alpha - \lambda_2.\beta).\varepsilon_2^0 = U(\beta).\varepsilon_1 = \varepsilon_1 + U'(0).\beta.\varepsilon_1 \quad {\rm modulo} \quad  \Phi_6  \\
   & (\alpha - \lambda_2.\beta).W(\beta).\varepsilon_1 = (\lambda_1-\lambda_2).\beta.W(0).\varepsilon_1 \quad {\rm modulo} \quad  \Phi_6 
   \end{align*}
   we obtain that
   $$ T'(0) = U'(0) - (p_1-1).W(0) .$$
   So to conclude, it is enough to compute \ $U'(0)$. As \ $\varepsilon^0_2 = S_{\lambda_2}(\beta).e_2$ \ we have
   \begin{align*}
   & (\alpha - \lambda_2.\beta).\varepsilon_2^0 = (1 + s_2.\beta).(\alpha - \lambda_2.\beta).e_2 \quad {\rm modulo} \quad \Phi_6 = \beta^2.E = b^2.E \\
   & \quad = (1 + s_2.\beta).\big[ a - \lambda_2.b + \tau.(a - 2\lambda_2.b).a\big].e_2  \quad {\rm modulo} \quad \Phi_6 \\
& \quad  = (1 + s_2.\beta).\big[ S.e_1 + \tau.(a - 2\lambda_2.b).(\lambda_2.b.e_2 + (1 + S'(0).b).e_1)\big]\quad {\rm modulo} \quad \Phi_6 \\
& \quad  = (1 + s_2.\beta).\big[ e_1 + S'(0).b.e_1 + \tau.(\lambda_1-\lambda_2).b.e_1 \big] \quad {\rm modulo} \quad \Phi_6 \\
& \quad  = e_1 + s_2.b.e_1+ \big[ S'(0) - \tau.(p_1-1)\big].b.e_1 \quad {\rm modulo} \quad \Phi_6\\
& \quad = \varepsilon_1 + \big[S'(0) -\tau.(p_1-1) + s_2 - s_1\big].\beta.\varepsilon_1 \quad {\rm modulo} \quad \Phi_6
\end{align*}
So we have \ $U'(0) = S'(0) - \tau.(p_1-1) + s_2 - s_1$ \ and 
$$ T'(0) = S'(0) - \tau.(p_1-1) + \tau.(p_1 -1).(2\lambda_1 + p_1 - 2) - \tau.(p_1 -1).(2\lambda_2 + p_2 - 2)
 .$$
This gives \ $L = (p_1-1).(p_1+p_2 -1).  \hfill \blacksquare$
  
  \bigskip

  \subsection{Constructions of quasi-invariant holomorphic \\ parameters.}

 \begin{thm}[The new parameter in rank 3.]\label{second.parm.}
Let \ $E$ \ be a rank \ $3$ \  fresco with fundamental invariants \ $\lambda_1, p_1,p_2$.
Let \ $S_1, S_2 \in \C[[b]]$ \ with \ $S_1(0) = 1$, $ S_2(0) = 1$,  such that  \ $E$ \ is isomorphic to 
 $$\A\big/\A.(a - \lambda_1.b).S_1^{-1}.(a - \lambda_2.b).S_2^{-1}.(a - \lambda_3.b) .$$
 Then the number
  $$\gamma(E) : = (p_2-1).S_1'(0) - (p_1-1).S_2'(0) $$
  depends only of the isomorphism class of \ $E$. For any unimodular change of  variable \ $\theta(a) = a + \tau.a^2 + \dots \in \C[[a]]$ \ we have 
  $$ \gamma(\theta_*(E)) = \gamma(E) - (p_1-1)(p_2-1)(p_1+p_2-1).\tau .$$
  \end{thm}
  
  \parag{Proof} Remark first that the classification of rank \ $3$ \ $[\lambda]-$primitive frescos which is given in the appendix implies immediately the theorem when we assume that \ $S_2 = (1 + \beta.b^{p_2})$ \ where \ $\beta$ \ is the parameter of \ $E\big/F_1$ : 
  with  this assumption, if \ $p_2 \not= 1$ \ we have \ $S_2'(0) = 0$ \ and then \ $\gamma(E) = (p_2-1).\gamma$, where \ $\gamma = S_1'(0)$. So the result is clear  thanks also to the propositions \ref{ordre 3 rk 3} and  \ref{chgt var. rk. 3}. For \ $p_2 = 1$ \ the number \ $S_2'(0)$ \ is the parameter \ $\beta$ \ of \ $E\big/F_1$, so the result is again clear.\\
  So the only point to prove is that for a general choice for \ $S_2$ \ and with \ $p_2 \not= 1$, \ $\gamma(E)$ \ coincides with the invariant \ $\gamma$ \ which appear in the classification of the appendix and in the proposition  \ref{chgt var. rk. 3}. Using the basis \ $e_1, e_2,e_3$ \ of \ $E$ \ deduced from the generator \ $e$ \ and satisfying
  \begin{align*}
  & (a - \lambda_3.b).e_3 = S_2.e_2  \quad {\rm with} \quad e_3 : = e \\
  & (a - \lambda_2.b).e_2 = S_1.e_1 \\
  & (a - \lambda_1.b).e_1 = 0 
  \end{align*}
  then we look for \ $\varepsilon_2 : = e_2 + T.e_1, \varepsilon_3 : = e_3 + X.e_2 + Y.e_1 $ \ in order to have 
  \begin{align*}
  & (a - \lambda_2.b).\varepsilon_2 = R_1.e_1 , {\rm with} \quad R_1(0) = 1, \\
  & (a - \lambda_3.b).\varepsilon_3 = (1 + \beta.b^{p_2}).\varepsilon_2 .
  \end{align*}
  Then this last relation implies
  \begin{align*}
  & (a - \lambda_3.b).\varepsilon_3 = S_2.e_2 + (\lambda_2-\lambda_3).b.X.e_2 + b^2.X'.e_2 + \\
  & \qquad  + X.S_1.e_1 + (\lambda_1-\lambda_3).b.Y.e_1 = b^2.Y'.e_1 = (1 + \beta.b^{p_2}).(e_2 + T.e_1) 
  \end{align*}
  So we have to consider the equations
  \begin{align*}
  & b^2.X' - (p_2-1).b.X = (1 + \beta.b^{p_2}) - S_2 \quad {\rm and} \tag{@}\\
  & X.S_1 - (p_1+p_2-2).b.Y + b^2.Y' = (1 + \beta.b^{p_2}).T   \tag{@@}
  \end{align*}
  for which we know that a solution exists for a good choice of \ $T \in \C[[b]]$.\\
  But in fact  we are only  interested in the computation of \ $R'_1(0)$ \ with the assumption that \ $p_2 \not= 1$. We have
  $$ (a - \lambda_2.b).(e_2 + T.e_1) = S_1.e_1 + (\lambda_1-\lambda_2).b.T.e_1 + b^2.T'.e_1$$
  which gives \ $R'_1(0) = S'_1(0) -  (p_1-1).T(0)$. So the computation of \ $T(0)$ \ is enough for our purpose. The equation \ $(@@)$ \ gives \ $X(0) = T(0)$ \ as \ $S_1(0) = 1$. The equation \ $(@)$ \ gives \ 
  $-(p_2-1).X(0) = -S_2'(0)$ \ for \ $p_2 \not= 1$ \ and we conlude then that
  $$ R_1'(0) = S_1'(0) - \frac{p_1-1}{p_2-1}.S_2'(0) \quad {\rm for} \quad p_2 \not= 1.$$
  This completes the proof. $\hfill \blacksquare$
  
  \begin{thm}\label{inv. param.}
  Let \ $\mathcal{S}(\lambda_1, \dots, \lambda_k)$ \ the set of isomorphism class of frescos with fundamental parameters \ $\lambda_1, \dots, \lambda_k$. Define the rational numbers \ $p_j : = \lambda_{j+1} - \lambda_j+1$ \ for \ $j \in [1,k-1]$. For \ $1 \leq i < j \leq k-2$ \ define the function
  $$ \pi_{i,j} : \mathcal{S}(\lambda_1, \dots, \lambda_k) \to \C $$
  by the formula
  $$ \pi_{i,j}(E) = (p_j-1).(p_{j+1}-1).(p_j+p_{j+1}-1).\gamma_i(E) - (p_i-1).(p_{i+1}-1).(p_i+p_{i+1}-1).\gamma_j(E) $$
  where \ $\gamma_h(E) : = \gamma(F_{h+2}\big/F_{h-1})$ \ for \ $h \in [1,k-2]$, with \ $F_h$ \ the \ $h-$th term in the principal J-H. of \ $E$.
  This function is a holomorphic parameter which is quasi-invariant of weigth \ $1$ \ by change of variables.
  \end{thm}
  
	  \parag{proof} This is an easy consequence of the theorem \ref{second.parm.} and the covariance of the principal J-H. sequence by change of variables. $\hfill \blacksquare$
	  
	\parag{Remark} Il we have a holomorphic family \ $\mathbb{E}$ \ of  frescos parametrized by a connected reduced complex space, then for each \ $[\lambda]$ \ we have a holomorphic family of \  $[\lambda]-$primitive frescos corresponding to the \ $[\lambda]-$primitive parts of the \ $\mathbb{E}(x)[\lambda]$ \ (thanks to proposition \ref{primitive part} ) and the theorem \ref{hol. ss part} gives a Zariski dense open set \ $X'$ \ in \ $X$ \ on which the semi-simple part of the holomorphic family \ $\mathbb{E}[\lambda]$ \ is holomorphic. Then we may apply the previous theorem to this family and get on \ $X'$ \ holomorphic functions which only depends on the the isomorphism classes of \ $\mathbb{E}(x)$ \ modulo an unimodular change of variables. \\
	These functions cannot be obtained, in general, without using the theorem on the semi-simple part, because, even for a \ $[\lambda]-$primitive fresco, the semi-simple part is not given, in general, by a term of the principal J-H. sequence.

\newpage

\section{Appendix : Classification of rank 3   frescos.}

\subsection{The main computation.}
We fix a rank \ $3$ \ frescos with fundamental invariants \ $\lambda_1,\lambda_2,\lambda_3$. We assume that inside a class modulo \ $\mathbb{Z}$ \ the \ $\lambda_i$ \ are such that \ $\lambda_j+j$ \ is increasing. In the case where a class modulo 1 has exactly two members we shall assume that they are \ $\lambda_2$ \ and \ $\lambda_3$. We fix a \ $\C[[b]]-$basis \ $e_1, e_2,e_3$ \ such that 
$$(a - \lambda_3.b).e_3 = S_2.e_2 \quad (a - \lambda_2.b).e_2 = S_1.e_1 \quad {\rm and} \quad (a -  \lambda_1.b).e_1 = 0 $$
where \ $S_2 = 1 + \beta.b^{p_2} $ \ when \ $\lambda_3 = \lambda_2 + p_2 -1$ \ with \ $p_2 \in \mathbb{N}^*$ \ and \ $\beta \in \C$ \ is the parameter of \ $E\big/F_1, F_1: = \C[[b]].e_1 \simeq E_{\lambda_1}$ \ 
and where \ $S_2 \equiv 1$ \ in all other cases. We assume \ $S_1(0) = 1$ \ and we want to give a unique normal form of the type  \ $\Sigma_1 = 1 + \gamma.b + \alpha.b^{p_1} + \delta.b^{p_1+p_2}$ \ to \ $S_1$, of course as simple as possible, up to an isomorphism of \ $E$ \ (so in fact by  a change of \ $\C[[b]]-$basis).\\
So we look for a new basis \ $\varepsilon_1, \varepsilon_2, \varepsilon_3$ \  with the following properties :
$$ (a - \lambda_3.b).\varepsilon_3 = S_2.\varepsilon_2 \quad (a - \lambda_2.b).\varepsilon_2 = \Sigma_1.\varepsilon_1 \quad {\rm and} \quad  (a -  \lambda_1.b).\varepsilon_1 = 0 .$$
We shall look to this new basis by looking for \ $X,Y,T \in \C[[b]]$ \ such that
$$ \varepsilon_3 : = e_3 + X.e_2 + Y.e_1 \quad  \varepsilon_2 : = e_2 + T.e_1 \quad {\rm and}  \quad  \varepsilon_1 = e_1 .$$
The relation \ $ (a - \lambda_3.b).\varepsilon_3 = S_2.\varepsilon_2 $ \  leads to the following two equations
\begin{align*}
& (\lambda_2-\lambda_3).b.X + b^2.X' = 0  \tag{@} \\
&  (\lambda_1 - \lambda_3).b.Y + b^2.Y' = S_2.T - X.S_1 . \tag{@@}
\end{align*}
Assume  \ $\lambda_3 - \lambda_2 \not\in \mathbb{N}$ ; the first equation gives \ $X = 0$ \ and the second equation becomes \ $b^2.Y' - (\lambda_3-\lambda_1).b.Y = S_2.T$. This forces \ $T(0) = 0$ \  as \ $S_2(0) = 1$ \ and defining \ $T = b.\tilde{T}$ \ the second equation becomes
\begin{equation*}
b.Y' - (\lambda_3-\lambda_1).Y = S_2.\tilde{T} \tag{@@bis}
\end{equation*}
Remark that  \ $\lambda_3 - \lambda_1$ \ is not in \ $\mathbb{N}$ \ with our convention and so  for any choice of \ $\tilde{T}$ \ there exists an unique solution in \ $Y$. So the equation \ $(@@bis)$ \ imposes no condition on \ $\tilde{T}$ \ in this case. \\

Now, when we have find a solution for \ $X,Y,T$ \ we want to choose \ $\Sigma_1$ \ in order to have
\begin{equation*}
(a - \lambda_2.b).\varepsilon_2 = \Sigma_1.\varepsilon_1 \quad {\rm which \ gives} \quad 
 (\lambda_1-\lambda_2).b.T + b^2.T' =\Sigma_1 - S_1  \tag{@@@}
\end{equation*}
Note that if \ $T(0) = 0$ \ this implies the condition \ $\Sigma_1 - S_1 \in b^2.\C[[b]]$. This means that in this case our change of basis keep the coefficient in \ $b$ \ fix in \ $S_1$.

\parag{Remark} If we want to keep the normal form for \ $E\big/F_1$ :
$$E\big/F_1  \simeq \A\big/\A.(a - \lambda_2.b).S_2^{-1}.(a - \lambda_3.b) $$
and also the condition \ $(a - \lambda_2.b).\varepsilon_2 \in Ker\,(a - \lambda_1.b) \setminus\{0\}$ \ we may always normalize by a non zero constant our initial basis in order  to have \ $\varepsilon_1 = e_1$ \ and then\ $ \varepsilon_2 = e_2 + T.e_1$ \ is necessary because we assume \ $S_1(0) = \Sigma_1(0) = 1$. Then the condition \ $S_2(0) = 1$ \ forces the fact that \ $ \varepsilon_3 : = e_3 + X.e_2 + Y.e_1 $. So the form we choose for the new \ $\C[[b]]-$basis does not reduce the generality. $\hfill \square$\\

\subsection{The different cases.}

\subsubsection*{The case  \ $\lambda_2 -\lambda_1 \not\in \mathbb{Z}$.}
 With our convention, this implies \ $\lambda_3 - \lambda_1 \not\in \mathbb{Z}$ \ so we have an unique solution \ $Y$ \ for any choice of \ $\tilde{T}$ \ for the equation \ $(@@bis)$. Then there exists an unique solution \ $\tilde{T}$ \ for the equation \ $(@@@)$ \ with \ $\Sigma_1 : = 1 + \gamma.b$ \ where \ $\gamma : = S'_1(0)$. In this case we have an unique \ $\gamma \in \C$ \ such that
 \begin{equation*}
 E \simeq \A\big/\A.(a - \lambda_1.b).(1 + \gamma.b)^{-1}.(a - \lambda_2.b).S_2^{-1}.(a - \lambda_3.b). \tag{1}
 \end{equation*}
 
\parag{Notation} Remark that with our convention, in all other case the fresco \ $E$ \ is \ $[\lambda]-$primitive. So \ $p_1$ \ and \ $ p_2 $ \ are natural integers.

\subsubsection*{The case where \ $p_1\geq 2, p_2 \geq 2$ \ with \ $F_2$ \ a theme \ ($\alpha \not= 0$).} We have \ $X = \sigma.b^{p_2-1}$ \ for some \ $\sigma \in \C$ \ and the equation \ $(@@)$ \ forces again  \ $T(0) = 0$. Now to have a solution to the equation
\begin{equation*}
b.Y' - (p_1+p_2-2).Y = S_2.\tilde{T} - \sigma.b^{p_2-2}.S_1 \tag{@@ter}
\end{equation*}
we must have no term in \ $b^{p_1+p_2-2}$ \ in the right handside of \ $(@@ter)$. Let \ $\alpha \not= 0$ \ be the parameter of \ $F_2$ \ which is also the coefficient of \ $b^{p_1}$ \ in \ $S_1$. Then to have a solution we must choose \ $\sigma$ \ in order that \ $\beta.t_{p_1-1} = \sigma.\alpha$ \ where \ $ t_{p_1-1}$ \ is the coefficient of \ $b^{p_1-1}$ \ in \ $T$. This shows that we may choose such a \ $\sigma$ \ for every choice of \ $T$. So the only constraint on \ $T$ \ is \ $T(0) = 0$ \ and we may choose \ $\Sigma = 1 + \gamma.b + \alpha.b^{p_1}$, where \ $\gamma$ \ is the coefficient of \ $b$ \ in \ $S_1$. We conclude that in this case we have
\begin{equation*}
 E \simeq \A\big/\A.(a - \lambda_1.b).(1 + \gamma.b + \alpha.b^{p_1})^{-1}.(a - \lambda_2.b).(1 + \beta.b^{p_2})^{-1}.(a - \lambda_3.b) \tag{2}
 \end{equation*}
where \ $\alpha \not= 0$ \ and \ $\beta$ \ are the parameters respectively of \ $F_2$ \ and \ $E\big/F_1$ \ and where \ $\gamma$ \ is a complex number uniquely determined by the isomorphism class of \ $E$.

\subsubsection*{The case where \ $p_1\geq 2, p_2 \geq 2$ \ with \ $F_2$ \ semi-simple \  ($\alpha = 0$).} We have to add to the previous computation the condition  \ $\beta.t_{p_1-1} = 0 $. Assume first that \ $\beta \not= 0$. So to solve the equation 
\begin{equation*}
b.\tilde{T}' - (p_1-2).\tilde{T} = (\Sigma_1 - S_1)\big/b^2  , \tag{@@@bis} 
\end{equation*}
 as \ $S_1$ \ has no term in \ $b^{p_1}$ \ from our assumption, we may choose \ $\Sigma_1 = 1 + \gamma.b$ \ where \ $\gamma$ \ is the coefficient of \ $b$ \ in \ $S_1$. For \ $\beta = 0$ \  the equation \ $(@@bis)$ \ does not imply a condition on \ $\tilde{T}$ \ and we may make the same choice. Then we find independantely of the value of \ $\beta$ \  the same normal form than in the previous case (but \ $\alpha = 0$ \ by assumption !).
 \begin{equation*}
  E \simeq \A\big/\A.(a - \lambda_1.b).(1 + \gamma.b)^{-1}.(a - \lambda_2.b).(1 + \beta.b^{p_2})^{-1}.(a - \lambda_3.b) . \tag{3}
  \end{equation*}

\subsubsection*{The case \ $p_1 \geq 2, p_2 = 1$.} 
In the previous computation we now must choose \ $\sigma = T(0) = t_0$ \ and to have a solution in \ $(@@)$ \ we must ask  to \ $S_2.T - \sigma.S_1$ \ to have no term in \ $b^{p_1}$ \ as \ $p_1 + p_2 - 1 = p_1$. This  forces \ $ t_{p_1} + \beta.t_{p_1-1} = t_0.\alpha$. Now \ $T$ \ is solution of \ $(@@@)$ \ which is now 
\begin{equation*}
b.T' - (p_1-1).T = (\Sigma_1- S_1)\big/b .\tag{@@@ter}
\end{equation*}
This  implies, as \ $\Sigma_1 : = 1 + \gamma.b + \alpha.b^{p_1} + \delta.b^{p_1+1} $,  the following conditions :
\begin{equation*}
 -(p_1-1).t_0 = \gamma - S_1'(0)  \qquad  \beta.t_{p_1-1} = t_0.\alpha - t_{p_1} \qquad   t_{p_1} = \delta - s_{p_1+1}
\end{equation*}
where \ $\alpha$ \ and  \ $s_{p_1+1}$ \ are the coefficient of \ $b^{p_1}$ \ and  \ $b^{p_1+1}$ \ in \ $S_1$. If \ $\beta \not= 0$, as we may choose the coefficient \ $t_{p_1-1}$ \ as we want for a solution of \ $(@@@ter)$, we may take \ $\gamma   = \delta = 0 $ . Then choose \ $\sigma = t_0 = S_1'(0)/(p_1-1)$ \ and we find  
\begin{equation*}
 E \simeq \A\big/\A.(a - \lambda_1.b).(1 + \alpha.b^{p_1})^{-1}.(a - \lambda_2.b).(1 + \beta.b)^{-1}.(a - \lambda_2.b) \qquad \beta \not= 0  . \tag{4}
 \end{equation*}
But for \ $\beta = 0$ \ the previous relations forces \ $t_{p_1} = \delta - s_{p_1+1} = (S'(0)- \gamma)\big/(p_1-1)$. So we may choose \ $\delta = 0$ \ and then \ $\gamma$ \ is determined uniquely. Then there exists an unique \ $\gamma \in \C$ \ such that we have 
\begin{equation*}
  E \simeq \A\big/\A.(a - \lambda_1.b).(1 + \gamma.b + \alpha.b^{p_1})^{-1}.(a - \lambda_2.b).(a - \lambda_2.b) \quad {\rm case} \quad \beta = 0 . \tag{4'}
  \end{equation*}

\parag{Remark}  For \ $\alpha \not= 0$ \ and \ $\beta = 0$ \ we may also choose  \ $\gamma = 0$ \ and  find an unique \ $\delta \in \C$ \ such that 
$$ E \simeq \A\big/\A.(a - \lambda_1.b).(1 + \alpha.b^{p_1} + \delta.b^{p_1+1})^{-1}.(a - \lambda_2.b).(a -\lambda_2.b)$$
where \ $\gamma$ \ and \ $\delta$ \ are linked by the relation \ $\alpha.\gamma = (p_1-1).\delta$.

\subsubsection*{The case \ $p_1 \geq 2, p_2 = 0$.} 
The equation \ $(@)$ \ gives \ $X = 0$ \ as \ $\lambda_3 - \lambda_2 = -1$ \ and then \ $T = b.\tilde{T}$ \ and equation \ $(@@bis)$ \ gives, as \ $S_2 = 1$,
$$ b.Y' - (p_1-2).Y = \tilde{T}.$$
 So we must have no \ $b^{p_1-1}$ \ term in \ $T$ \ to find a solution \ $Y$.
Then the equation \ $(@@@)$ \ leads to choose \ $\Sigma_1 = 1 + \gamma.b + \alpha.b^{p_1} $ \ where  \ $\gamma : = S_1'(0)$ \ and \ $\alpha$ \  is the parameter of \ $F_2$. The previous computation implies the uniqueness of \ $\gamma$. Then in this case 
\begin{equation*}
E \simeq \A\big/\A.(a - \lambda_1.b).(1 + \gamma.b +  \alpha.b^{p_1})^{-1}.(a - \lambda_2.b).(a - (\lambda_2-1).b) . \tag{5}
\end{equation*}

\subsubsection*{The case \ $p_1=1, p_2 \geq 1$.} 
The equation \ $(@)$ \ gives \ $X = \sigma.b^{p_2-1}$ \ and this implies \ $T = b.\tilde{T}$ \ for \ $p_2 \geq 2$ \ and \ $T(0) = \sigma$ \ for \ $p_2 = 1$.\\
 Assume first \ $p_2 \geq 2$. Then we have from \ $(@@bis)$ 
  $$ b.Y' - (p_2-1).Y = S_2.\tilde{T} - \sigma.b^{p_2-2}.S_1 .$$
  So to have a solution we need that the right handside has no term in \ $b^{p_2-1}$. If \ $\alpha$,  the parameter of \ $F_2$,  which is the coefficient of \ $b$ \ in \ $S_1$, is not zero, there is a unique  choice of \ $\sigma$ \  to satisfy this condition. So the conditions on \ $T$ \ are \ $T(0) = 0$ \ and \ $b^2.T' = \Sigma_1-S_1$. Then we may choose \ $\Sigma = 1 + \alpha.b$ \ and obtain the normal form
  \begin{equation*}
   \A\big/\A.(a - \lambda_1.b)(1 + \alpha.b)^{-1}(a - \lambda_1.b).(1 + \beta.b^{p_2})^{-1}.(a - \lambda_3.b) \quad {\rm for} \quad \alpha \not= 0,  p_2 \geq 2. \tag{6}
   \end{equation*}
  For \ $\alpha = 0$, we have the extra condition on \ $T$ \ that \ it has no \ $b^{p_2}$ \ term which forces that \ $\Sigma_1 - S_1$ \ has no \ $b^{p_2+1}$ \ term. So we must choose \ $\Sigma_1 = 1 + \delta.b^{p_2+1}$ \ where \ $\delta$ \ is the coefficient of \ $b^{p_2+1}$ \ in \ $S_1$. Then \ $E$ \ is isomorphic to
  \begin{equation*}
   \A\big/\A.(a - \lambda_1.b).(1 + \delta.b^{p_2+1})^{-1}(a - \lambda_1.b).(1 + \beta.b^{p_2})^{-1}.(a - \lambda_3.b) \quad\alpha = 0, p_2 \geq 2 . \tag{6'}
   \end{equation*}
  
  In the case \ $p_2 = 1$ \ we choose \ $T(0) = \sigma$ \ and \ $(@@bis)$ \ is now 
   $$b^2.Y' = S_2.T - T(0).S_1 $$
   which have a solution if and only if \ $T'(0) + \beta.T(0) = \alpha.T(0) $ \ where \ $\beta$ \ is the parameter of \ $E\big/F_1$. The equation \ $(@@@)$ \ is then \ $b^2.T' = \Sigma_1 - S_1$. If \ $\beta \not= \alpha$ \ we may choose \ $T'(0) $ \ as we want and deduce \ $T(0) = \sigma$. So \ $\Sigma_1 = 1 + \alpha.b $.\\
   If \ $\beta = \alpha$ \ we must choose \ $\Sigma_1 = 1 + \alpha.b + \delta.b^2$ \ where \ $\delta$ \ is the coefficient of \ $b^2$ \ in \ $S_1$. So we have normal forms
   \begin{align*}
   & \A\big/\A.(a - \lambda_1.b).(1 + \alpha.b)^{-1}.(a - \lambda_1.b).(1 +\beta.b)^{-1}.(a - \lambda_3.b) \quad {\rm for} \quad \beta \not = \alpha \tag{6''} \\
   & \A\big/\A.(a - \lambda_1.b).(1 + \alpha.b + \delta.b^2)^{-1}.(a - \lambda_1.b).(1 + \alpha.b)^{-1}.(a - \lambda_1.b) \quad \beta = \alpha. \tag{6'''}
   \end{align*}

 \subsubsection*{The case \ $p_1 = 1, p_2 = 0$.} 
 The  equation \ $(@)$ \ gives \ $X = 0$ \ and this implies \ $T = b.\tilde{T}$. We have \ $\lambda_3-\lambda_1 = -1$ \  the equation \ $(@@bis)$ \  gives no condition on\ $\tilde{T}$. To solve the equation \ $(@@@bis)$ \ with \ $p_1 = 1$ \  there is again no condition on \ $\tilde{T}$, so we have 
 \begin{equation*}
   E \simeq \A\big/\A.(a - \lambda_1.b).(1 + \gamma.b )^{-1}.(a - \lambda_1.b).(a - (\lambda_1-1).b) \quad {\rm in \  this \ case}. \tag{7}
   \end{equation*}
 
\subsubsection*{The case \ $p_1 = 0$.} 
Assume first \ $p_2 \geq 2$. Again we have \ $X = \sigma.b^{p_2-1}$ \ and \ $T = b.\tilde{T}$. Now in equation \ $(@@bis)$ \ we have \ $\lambda_3 - \lambda_1 = p_2-2 \geq 0$, so,  to have a solution \ $Y$ \ we must have no term in \ $b^{p_2-2}$ \ in \ $\tilde{T}$. This implies via equation \ $(@@@bis)$ \ that we have no \ $b^{p_2}$ \ in \ $\Sigma_1 - S_1$. Then we have to  choose \ $\Sigma_1 = 1 + \gamma.b + \delta.b^{p_2}$ \ with \ $\gamma = S_1'(0)$ \ and \ $\delta$ \ the coefficient of \ $b^{p_2}$ \ in \ $S_1$. Then in this case we have the normal form
 \begin{equation*}
  \A\big/\A.(a - \lambda_1.b).(1 + \gamma.b + \delta.b^{p_2})^{-1}.(a - (\lambda_1-1).b).(1 + \beta.b^{p_2})^{-1}.(a - \lambda_3.b)\quad p_2 \geq 2  . \tag{8}
  \end{equation*}
 
 \bigskip
 
For \ $p_2 = 0$ \ we have \ $X = 0$ \ and again \ $T = b.\tilde{T}$. Then there is no  condition on \ $\tilde{T}$ \ to solve \ $(@@bis)$ \ and we conclude that we may choose \ $\Sigma_1 = 1 + \gamma.b$.\\
 For \ $p_2 = 1$ \ we must choose \ $T(0) = \sigma = X$ \ and we want to solve 
 $$ b^2.Y' + b.Y = S_2.T - T(0).S_1 $$
 which has a solution for any choice of \ $T$. Then the equation \ $(@@@)$ \ is 
 $$ b.T' + T = (\Sigma_1 - S_1)\big/b $$
 and choosing \ $\sigma  = - S_1'(0)$ \ we may have \ $\Sigma = 1$ \ in this case.\\
 This may be summarized as follows
 \begin{align*}
 &  p_2 = 1 \quad 
  E \simeq \A\big/\A.(a - \lambda_1.b).(a - (\lambda_1-1).b).(1 + \beta.b)^{-1}.(a - (\lambda_1-1).b) \tag{8'} \\
& p_2 = 0 \quad 
 E \simeq \A\big/\A.(a - \lambda_1.b).(1 + \gamma.b)^{-1}.(a - (\lambda_1-1).b).(a - (\lambda_1-2).b) \tag{8''}
\end{align*}
 
 \parag{Remarks} 
 \begin{enumerate}
 \item The only choices of \ $p_1, p_2$ \ where the number \ $\gamma$ \  is not determined by the isomorphism class of the  fresco \ $E$ \ with fundamental invariants \ $(\lambda_1, p_1,p_2)$ \  satisfy \ $p_2 = 1$. So for \ $p_2 \not= 1$ \ the function defined on \ $\mathcal{F}(\lambda_1, \lambda_2, \lambda_3)$ \ by this number \ $\gamma$ \ is a non trivial holomorphic parameter.\\
  Note that the proposition \ref{chgt var. rk. 3} implies that for \ $p_1\not= 1$ \ and \ $p_1 + p_2 \not= 1$, with the changes of variable \ $\theta(a) = a + \tau.a^2$,  the number \ $\gamma(\theta_*(E))$ \ takes all complex values (with fix \ $\alpha$ \ and \ $\beta$) when \ $\tau$ \ describes \ $\C$.
  \item Note that for \ $p_1 \not\in \mathbb{N}$ \ and \ $p_1+p_2 \in \mathbb{N}^*$ \ we may change the order on \ $\mathbb{Q}\big/\mathbb{Z}$ \ to have new fundamental invariants \ $\mu_1 = \lambda_2 + 1, \mu_2 = \lambda_1 - 1, \mu_3 = \lambda_3$ \ so we have \ $p'_1 = - p_1$ \ and \ $p'_2 = p_1+ p_2 -1 $ \ where \ $\mu_2 = \mu_1 + p'_1 - 1, \mu_3 = \mu_2 + p'_2 - 1$. So we reach with this new order our initial hypothesis in the case where  a class modulo \ $\mathbb{Z}$ \ contains exactly two of the fundamental invariants.
  \end{enumerate}
 
 \newpage
 
  \section{Bibliography.}

\begin{itemize}

\item{[Br.70]} Brieskorn, E. {\it Die Monodromie der Isolierten Singularit{\"a}ten von Hyperfl{\"a}chen}, Manuscripta Math. 2 (1970), p. 103-161.\\

\item {[B.93]}  Barlet, D. {\it Th\'eorie des (a,b)-modules I}, Univ. Ser. Math. Plenum (1993), p.1-43.\\

\item {[B.08]}Barlet, D. {\it Sur certaines singularit\'es d'hypersurfaces II}, J. Alg. Geom. 17 (2008), p. 199-254.\\

\item {[B.09]}  Barlet, D. {\it P\'eriodes \'evanescentes et (a,b)-modules monog\`enes}, Bollettino U.M.I (9) II (2009), p.651-697.\\

\item {[B.10]} Barlet, D. {\it Le th\`eme d'une p\'eriode \'evanescente}, preprint de l'Institut E. Cartan (Nancy) 2009  $n^0 33$,  57 pages, math. arXiv. 1110 1353 (english version).\\

\item{[B.12]}  Barlet, D. {\it Asymptotics of a vanishing period : General existence theorem and basic properties of frescos}  preprint of the Institut E. Cartan (Nancy) 2012  $n^0 1$ 38 pages, math. arXiv. 1201 2757.\\

\item{[Kz.11]} Karwasz, P. {\it Hermitian (a,b)-modules and Saito' higher residue pairings"} math. arXiv. 1104 1505.\\

 \item{[K.76]} Kashiwara, M. {\it b-function and holonomic systems}, Inv. Math. 38 (1976) p. 33-53.\\

\item{[M.74]} Malgrange, B. {\it Int\'egrale asymptotique et monodromie}, Ann. Sc. Ec. Norm. Sup. 7 (1974), p.405-430.\\

\end{itemize}

\end{document}